\documentclass[a4paper,fleqn]{cas-sc}

\usepackage[numbers,sort&compress]{natbib}

\usepackage{accents}	
\usepackage{stmaryrd}   
\usepackage{subcaption} 
\usepackage{bm} 
\usepackage{xfrac} 
\usepackage{booktabs} 

\usepackage{float} 

\usepackage[nameinlink,noabbrev]{cleveref} 

\usepackage{tikz}

\def\tsc#1{\csdef{#1}{\textsc{\lowercase{#1}}\xspace}}
\tsc{BB}
\tsc{AR}
\tsc{DK}
\tsc{GG}

\newdefinition{remark}{Remark}

\newcommand{\norm}[1]{\left\lVert#1\right\rVert}

\newcommand\stateL[1]{\mathbf{#1}}
\newcommand\stateG[1]{\boldsymbol #1}

\newcommand{\Nabla} {\vec{\nabla}}
\newcommand\acclrvec[1]{\accentset{\,\leftrightarrow}{#1}}	
\newcommand{\blocktensor}[1]{\acclrvec{{\mathbf #1}}}

\newcommand{\FV}{{\mathrm{FV}}}
\newcommand{\numfluxb}[1]{\hat{\mathbf{#1}} }

\def\NN{\mathcal{N}} 

\def\cc{\hat{\mathtt{c}}} 
\def\mm{\mathtt{m}} 
\def\NN{\mathcal{N}} 
\def\dd{\hat{\mathtt{d}}} 

\begin{document}
\let\WriteBookmarks\relax
\def\floatpagepagefraction{1}
\def\textpagefraction{.001}

\shorttitle{Invariant-domain-preserving mortar flux for LGL-DGSEM}    

\shortauthors{Bolm, Rueda-Ramírez, Kuzmin, Gassner} 

\title [mode = title]{Invariant-domain-preserving limiting with Adaptive Mesh Refinement for Legendre--Gauss--Lobatto Discontinuous Galerkin Spectral Element Methods}  



%

\author[1]{Benjamin Bolm}[type=editor,
orcid=0009-0009-6055-8472]
\cormark[1]
\ead{benjamin.bolm@uni-koeln.de}

\affiliation[1]{organization={Department of Mathematics and Computer Science, University of Cologne},
            addressline={Weyertal 86-90}, 
            city={Cologne},
            postcode={50931}, 
            state={NRW},
            country={Germany}}

\author[2]{Andrés M. Rueda-Ramírez}[
orcid=0000-0001-6557-9162
]
\ead{am.rueda@upm.es}
\affiliation[2]{organization={ETSIAE-UPM-School of Aeronautics, Universidad Politécnica de Madrid (UPM)},
            city={Madrid},
            country={Spain}}

\author[3]{Dmitri Kuzmin}[
orcid=0000-0003-0084-9752
]
\ead{kuzmin@math.uni-dortmund.de}
\affiliation[3]{organization={Institute of Applied Mathematics (LS III), TU Dortmund University},
            addressline={Vogelpothsweg 87}, 
            city={Dortmund},
            postcode={D-44227}, 
            state={NRW},
            country={Germany}}

\author[1, 4]{Gregor J. Gassner}[
orcid=0000-0002-1752-1158
]
\ead{ggassner@uni-koeln.de}
\affiliation[4]{organization={Center for Data and Simulation Science, University of Cologne},
            addressline={Weyertal 86-90}, 
            city={Cologne},
            postcode={50931}, 
            state={NRW},
            country={Germany}}

\begin{abstract}
We present an invariant-domain-preserving (IDP) treatment of nonconforming interfaces for Legendre--Gauss--Lobatto Discontinuous Galerkin Spectral Element Methods (LGL-DGSEM) with adaptive mesh refinement (AMR) on Cartesian meshes. 
The proposed methodology extends recently developed convex limiting and graph-viscosity frameworks for DGSEM to meshes containing hanging nodes. 

Starting from a conservative mortar formulation, we derive low-order interface fluxes that satisfy the requirements of invariant-domain-preserving discretizations. 
To avoid the excessive diffusion associated with fully connected mortar couplings, a sparsification strategy based on LGL subcell characteristic functions is introduced, yielding compact interface stencils. 
The resulting mortar fluxes remain conservative, reduce to the standard conforming formulation on matching interfaces, and naturally fit into graph-viscosity-based low-order schemes used for convex limiting. 

The proposed construction provides the missing ingredient required to combine high-order DGSEM discretizations, invariant-domain-preserving limiting, and adaptive mesh refinement within a unified framework for nonlinear hyperbolic conservation laws.
We provide numerical verifications of the properties of the proposed scheme and run challenging simulations that require positivity limiting and shock-capturing.
\end{abstract}




\maketitle

\section{Introduction}

High-order discontinuous Galerkin (DG) methods have become a popular framework for the numerical approximation of hyperbolic conservation laws due to their excellent accuracy, geometric flexibility, and suitability for modern high-performance computing architectures.
Among these methods, the Discontinuous Galerkin Spectral Element Method (DGSEM) based on Legendre--Gauss--Lobatto (LGL) collocation points is particularly attractive because of its tensor-product structure, diagonal mass matrix, and efficient implementation on modern hardware.
Moreover, split-form DGSEM discretizations possess favorable nonlinear stability properties and can be constructed to satisfy discrete entropy conservation or entropy stability for a broad class of hyperbolic systems \cite{Fisher2013a,gassner2016split,ranocha2021efficient}.

Despite these advantages, high-order DG methods may generate spurious oscillations in the vicinity of discontinuities or under-resolved solution features.
Such oscillations can violate essential physical constraints, including positivity of density and pressure for the compressible Euler equations or, more generally, invariant domain properties associated with the underlying system of conservation laws.
Preserving these constraints is crucial for the robustness and reliability of simulations involving shocks, strong rarefactions, low-density regions, or under-resolved turbulence.

The development of nonlinear stabilization techniques for high-order finite element and DG methods has a long history.
A particularly influential class of approaches is provided by algebraic flux correction (AFC) and flux-corrected transport (FCT) methods, which combine a robust low-order discretization with a high-order target scheme through carefully designed antidiffusive fluxes and limiters \cite{zalesak1979, selmin1987b, kuzmin2010a,kuzmin2012,lohmann2017,lohner2008}.
Over the last two decades, AFC methodologies have been extended to continuous and discontinuous finite element discretizations of hyperbolic conservation laws and have demonstrated their ability to enforce positivity, maximum principles, and other physically motivated bounds while retaining high-order accuracy in smooth regions \cite{kuzmin2023}.

Recently, Guermond, Popov, and collaborators developed a rigorous invariant-domain-preserving (IDP) framework based on graph viscosity and convex limiting \cite{Guermond2016,guermond2018,guermond2019}.
In this family of AFC approaches for hyperbolic problems, a provably invariant-domain-preserving low-order scheme is combined with a high-order approximation through limited flux corrections.
The resulting methodology provides strong theoretical guarantees, including preservation of convex invariant sets and compatibility with entropy inequalities.
These ideas have subsequently been incorporated into high-order DG discretizations \cite{Pazner2020, hajduk2021,hajduk2020,RUEDARAMIREZ2022, ruedaramirez2024} and have established a systematic pathway towards robust, invariant-domain-preserving methods for nonlinear hyperbolic systems.

For nodal DG discretizations, several recent developments have demonstrated that sparse low-order operators are essential for obtaining accurate and efficient convex limiting schemes \cite{Pazner2020, RUEDARAMIREZ2022}.
Fully connected graph-viscosity operators lead to excessive numerical diffusion, especially for high polynomial degrees.
To address this issue, sparse invariant-domain-preserving DG discretizations based on compact subcell stencils were proposed in the context of Bernstein \cite{hajduk2021,hajduk2020} and discontinuous Galerkin finite element methods \cite{RUEDARAMIREZ2022, ruedaramirez2024}.
These approaches combine the robustness of low-order subcell discretizations with the accuracy of high-order DG methods and provide an effective foundation for modern convex limiting strategies.

While these developments have significantly improved the robustness of DGSEM discretizations, large-scale simulations additionally require adaptive mesh refinement (AMR) to efficiently resolve localized structures while controlling computational cost.
In DGSEM, the use of AMR introduces nonconforming interfaces with hanging nodes between coarse and fine elements.
The standard treatment of such interfaces relies on mortar methods and projection operators, which preserve conservation and high-order accuracy across nonmatching interfaces \cite{KOPRIVA1996475}.
These mortar techniques have become a standard ingredient of adaptive spectral element discretizations \cite{ranocha2022adaptive}.
However, the corresponding interface operators are generally not designed to satisfy invariant-domain-preserving properties and therefore do not directly fit into existing graph-viscosity and convex limiting frameworks.

The objective of this work is to bridge this gap.
We develop a conservative and invariant-domain-preserving mortar formulation for nonconforming interfaces in LGL-DGSEM discretizations with adaptive mesh refinement.
Starting from a conservative all-to-all mortar coupling, we derive a low-order interface discretization that is compatible with graph-viscosity-based convex limiting.
To avoid the excessive diffusion associated with dense mortar couplings, we introduce a sparsification strategy based on characteristic functions of LGL subcells, yielding compact interface stencils whose size is independent of the polynomial degree.
The resulting mortar fluxes remain conservative, naturally reduce to the standard conforming-interface formulation in the absence of hanging nodes, and fit seamlessly into existing IDP DGSEM frameworks.

The main contribution of this paper is therefore the construction and analysis of sparse invariant-domain-preserving mortar fluxes for nonconforming DGSEM discretizations.
These fluxes provide the missing low-order interface building block required for combining high-order DGSEM, convex limiting, and adaptive mesh refinement within a unified framework for nonlinear hyperbolic conservation laws.

The paper is organized as follows. In \Cref{sec:numerical_methods}, we briefly review the LGL-DGSEM, derive the stable low-order mortar formulation with sparse connectivity stencil, and prove its invariant-domain-preserving property. We summarize the limiting approach at the nonconforming interfaces and briefly describe a positivity-preserving approach we used when transferring the solution between coarse and fine meshes. In \Cref{sec:numerical_results}, we validate the derived method and use it to perform challenging simulations of the compressible Euler equations. Finally, we conclude in \Cref{sec:conclusion}.
All the implementations in the paper are carried out using the open-source code \texttt{Trixi.jl}~\cite{ranocha2022adaptive, schlottkelakemper2021purely, schlottkelakemper2025trixi}.

\section{Numerical Methods}\label{sec:numerical_methods}
The following derivation of the Discontinuous Galerkin Spectral Element Method (DGSEM) and its notation is based on \cite{ruedaramirez2024}. Consider the system of hyperbolic conservation laws of the form
\begin{equation}\label{eq:conservationlaw}
    \frac{\partial \stateL{u}}{\partial t} + \Nabla \cdot \blocktensor{f}(\stateL{u}) = \stateL{0}, \quad \text{in}\; \Omega \times \mathbb{R}^+,
\end{equation}
where $\Omega \subseteq \mathbb{R}^d$ is the spatial computational domain with space dimensions $d\in\{1, 2, 3\}$. The vector $\stateL{u}(\vec{x}, t)\in \mathbb{R}^{n_\text{eq}}$, depending on the spatial location $\vec{x}$ and time $t$, contains conserved quantities. The flux function $\blocktensor{f}:\mathcal{G}\to\mathbb{R}^{d\times n_\text{eq}}$ depends on $\stateL{u}:\bar\Omega \times \mathbb{R}^+_0 \to \mathcal{G}$, where the set $\mathcal{G}\subseteq\mathbb{R}^{n_\text{eq}}$ is called an invariant domain if $\mathcal{G}$ is convex and $\stateL{u}(\vec{x}, t)\in \mathcal{G}$ for all $(\vec{x}, t)\in\bar\Omega\times \mathbb{R}^+_0$.
The system in \eqref{eq:conservationlaw} comes with an initial condition $\stateL{u}(\cdot, 0) = \stateL{u}_0$, and suitable boundary conditions on the boundary $\partial \Omega$ of the spatial domain.

\subsection{The Discontinuous Galerkin Spectral Element Method}
We divide the spatial domain $\Omega$ into $K$ non-overlapping elements $\mathcal{T} = \{\Omega_1, \dots, \Omega_K\}$, and approximate the solution within each element by a polynomial of degree $N$. This yields a DG approximation $\stateL{u}$ which may be discontinuous at element interfaces and lies in the following space:
\begin{equation}
    \mathcal{V}^N = \lbrace u \in L_2(\Omega) \colon u \rvert_{\Omega_e} \in \mathcal{P}^N(\Omega_e) \, \forall \Omega_e \in \mathcal{T} \rbrace.
\end{equation}

To derive a semi-discrete update formula, we focus on the solution of \eqref{eq:conservationlaw} within an element $\Omega_e$ and multiply the equation by an arbitrary polynomial test function $\stateG{\phi}\in(\mathcal{P}^N(\Omega_e))^{n_\text{eq}}$. Integrating over the element $\Omega_e$ and using integration by parts moves the derivation from the solution to the test function and yields the weak formulation of the system of conservation laws
\begin{equation}\label{eq:weakform}
    \int_{\Omega_e} \frac{\partial \stateL{u}}{\partial t} \cdot \stateG{\phi} \mathrm{d}\vec{x} - \int_{\Omega_e} \blocktensor{f}(\stateL{u}) \cdot \Nabla \stateG{\phi} \mathrm{d}\vec{x} + \oint_{\partial \Omega_e} \stateG{\phi} \cdot \numfluxb{f} \mathrm{d}s = \stateL{0}.
\end{equation}

As mentioned before, the approximation is not necessarily unique on element interfaces, yielding the flux $\numfluxb{f}\approx \blocktensor{f} \cdot \vec{n}\rvert_{\partial\Omega_e}$ with outside-pointing normal vector $\vec{n}$ to be not unique on $\partial \Omega_e$ as well. We use an approximate Riemann solver that receives the state values from both sides of the interface and returns a numerical flux.

The approximated polynomial solution $\stateL{u}$ within the element $\Omega_e$ is represented using a Lagrange basis with $(N+1)^d$ Legendre--Gauss--Lobatto (LGL) interpolation points. Multidimensional LGL approximations are constructed using a tensor-product approach. To use the same polynomial basis functions on every element, we introduce mappings $X^e:\tilde\Omega=[-1,1]^d \to \Omega_e$ that transform the elements from reference space to physical space ($\xi \mapsto \vec{x}$ for $\xi\in\tilde\Omega$ and $\vec{x}=X^e(\xi)\in\Omega_e$). In the present work, we will focus on Cartesian meshes so the transformation consists only of shifting and scaling.

The one-dimensional Lagrange basis functions of the reference element $\tilde\Omega=[-1,1]$ are denoted by $\{\varphi_i\}_{i=0}^N$. Since those functions form a basis, it is sufficient to use their vector-valued counterparts instead of the arbitrary test function  $\stateG{\phi}$ in \eqref{eq:weakform}.
Moreover, for numerical integration, we apply a quadrature rule with the same LGL collocation nodes on the reference element $\tilde\Omega = [-1,1]^d$. The 1D quadrature weights in $[-1,1]$ are
\begin{equation}
    \omega_i = \int_{-1}^1 \varphi_i d\xi.
\end{equation}
In this work, we refer to nodes with their local nodal tensor-product index, i.e., in two spatial dimensions, the node $ij$ with $i, j\in\{0, \dots, N\}$ is the $i$th local degree of freedom in $\xi_1$-direction and the $j$th degree of freedom in $\xi_2$-direction on the reference element.

To maintain clarity while still working with the actual element interfaces, we will introduce the method in two spatial dimensions ($d = 2$).
With the described approach and some manipulations, the semi-discrete evolution equation of a two-dimensional LGL-DGSEM discretization of \eqref{eq:conservationlaw} for node $ij$ in element $\Omega_e$ can be written as \cite{ranocha2021efficient,RUEDARAMIREZ2022}
\begin{equation}\label{eq:evolutionformula}
\begin{split}
    J \omega_i \omega_j \dot{\stateL{u}}_{ij}
    =&\; \omega_j \left(\stateL{F}_{ij}^{\mathrm{Vol(1)}} + \delta_{i0} \numfluxb{f}_{(0j,L)} - \delta_{iN} \numfluxb{f}_{(Nj, R)} \right)
    +\; \omega_i \left(\stateL{F}_{ij}^{\mathrm{Vol(2)}} + \delta_{j0} \numfluxb{f}_{(i0,L)} - \delta_{jN} \numfluxb{f}_{(iN,R)} \right),
\end{split}
\end{equation}
where $\dot{\stateL{u}} = \frac{\partial \stateL{u}}{\partial t}$, $J=\frac{\Delta x \Delta y}{4}$ is the determinant of the geometric mapping Jacobian from the reference space $\tilde\Omega$ to the physical space of $\Omega_e$, $\delta$ is Kronecker's delta function, and $\Delta x$ and $\Delta y$ are the dimensions of the DG element.
The presented ansatz is the so-called nodal collocation variant of the DG method.

In this work, we use a split-form DG discretization of the divergence term of \eqref{eq:conservationlaw} based on two-point fluxes. See, e.g., \cite{gassner2016split,Fisher2013a,ranocha2021efficient}.
In this framework, $\stateL{F}_{ij}^{\mathrm{Vol}(m)}$ denotes the volume integral part in the $m$th direction of the reference space.
For instance, the $\xi_1$ component of the volume-integral term reads
\begin{equation}
   \stateL{F}_{ij}^{\mathrm{Vol}(1)} = -\sum_{k=0}^N \bar{S}_{ik} \stateL{f}^{1*}_{(ij, kj)},
\end{equation}
where we use the skew-symmetric volume integral matrix $\bar{S}=2\bar{Q}-\bar{B}$, which is a function of the summation-by-parts (SBP) derivative matrix $\bar{Q}_{ij} = \omega_i \varphi_j'(\xi_i)$ and the so-called boundary matrix $\bar{B}=\text{diag}(-1,0, \dots, 0, 1)$.
The flux $\stateL{f}^{1*}_{(ij,kj)} := \stateL{f}^{1*}(\stateL{u}_{ij},\stateL{u}_{kj}, \vec{n}_1)$ is the two-point contravariant numerical volume flux between nodes $ij$ and $kj$, typical for the present split-form DGSEM (see for instance \cite{RUEDARAMIREZ2022}). The flux is evaluated in the direction of the normal vector $\vec{n}_1$, which, in this Cartesian case, corresponds to $(\frac{\Delta y}{2}, 0)^T$.
The fluxes in the $y$ direction, $\stateL{f}^{2*}_{(\cdot,\cdot)}$, are computed analogously.

More relevant for this work are the interface contributions with the numerical fluxes $\numfluxb{f}$ in \eqref{eq:evolutionformula}. The subscripts $L$ and $R$ indicate the direction of the flux to the left/bottom (negative directions) and to the right/top (positive directions), respectively. We consider the flux at an interface $S\in\partial \Omega_e$ in positive $\xi_1$-direction connecting element $\Omega_e$ with its right neighbor. The flux reads as
\begin{equation}
    \numfluxb{f}_{(Nj,R)} = \numfluxb{f}(\stateL{u}_{Nj}, \stateL{u}^+, \vec{n}_S)
\end{equation}
and is computed using the inner value $\stateL{u}_{Nj}$ of node $Nj$ in $\Omega_e$, the outer value $\stateL{u}^+$ from the right neighbor and the associated outward-pointing normal vector $\vec{n}_S$, which simplifies for Cartesian meshes to a scaled first standard unit vector and is constant for all nodes at $S$.

All diagonal-norm SBP discretizations of conservation laws -- including the present split-form LGL-DGSEM -- can be rewritten in the so-called \textit{flux-differencing} formulation \cite{Fisher2013a}:
\begin{equation}\begin{split}\label{eq:fluxdifferencingform}
    J\dot{\stateL{u}}_{ij} = \frac{1}{\omega_i} \left( \numfluxb{f}_{((i-1)j,ij)} - \numfluxb{f}_{(ij,(i+1)j)} \right)
    + \frac{1}{\omega_j} \left(\numfluxb{f}_{(i(j-1),ij)} - \numfluxb{f}_{(ij,i(j+1))} \right), \qquad \forall i, j=0, \dots, N,
\end{split}\end{equation}
where the indices $-1$ and $N+1$ refer to outer states. The symmetric and consistent fluxes $\numfluxb{f}_{(ij,k\ell)}=\numfluxb{f}_{(k\ell,ij)}$ in $\xi_1$ direction read as follows:
\begin{equation}\begin{split}
\numfluxb{f}_{((-1)j,0j)}  &= \numfluxb{f}_{(0j,L)},\\
\numfluxb{f}_{(ij,(i+1)j)} &= \sum_{l=0}^i \sum_{k=0}^N \bar{S}_{lk} \stateL{f}^{1*}_{(lj,kj)} \qquad \forall i=0, \ldots, N-1,\\
\numfluxb{f}_{(Nj,(N+1)j)} &= \numfluxb{f}_{(Nj,R)},
\end{split}\end{equation}
and the fluxes in the $\xi_2$ direction are obtained analogously.

The general idea of nodal convex limiting is now to combine a high-order accurate DG scheme with a compatible low-order FV method,
\begin{equation}\begin{split}\label{eq:FVmethod}
    J\dot{\stateL{u}}_{ij}^{\FV} = \frac{1}{\omega_i} \left( \numfluxb{f}^{\FV}_{((i-1)j,ij)} - \numfluxb{f}^{\FV}_{(ij,(i+1)j)} \right)
    + \frac{1}{\omega_j} \left(\numfluxb{f}^{\FV}_{(i(j-1),ij)} - \numfluxb{f}^{\FV}_{(ij,i(j+1))} \right), \qquad \forall i, j=0, \dots, N,
\end{split}\end{equation}
to ensure invariant-domain-preserving (IDP) properties.

The index notation from above, where nodes (in 2 dimensions) are indicated by 2 local indices $ij, \, i,j=0,\dots, N$, perfectly fits to the locality of DG and it is relatively similar to the actual implementation. Nevertheless, it has also disadvantages. For instance, looping through those nodes without loss of generality always requires 2 variable indices.
Because of that, we will use a global index notation as well. In the global notation, each node is indicated by one unique global index $i$.

An equivalent global index notation for \eqref{eq:FVmethod}, for the particular case of LLF flux, reads (see \cite{KUZMIN20044915, kuzmin2010a, kuzmin2020, RUEDARAMIREZ2022}),
\begin{equation}\label{eq:LowOrderScheme_globalIndices}
    \mm_i \dot{\stateL{u}}^{\FV}_{i} = -\sum_{j \in \NN(i)} \cc_{(i,j)} \cdot \blocktensor{f}_{j} + \sum_{j \in \NN(i)} \dd_{(i,j)} \left( \stateL{u}_j - \stateL{u}_i \right),
\end{equation}
where $\mm_i$ is the mass matrix entry that corresponds to degree of freedom $i$, $\dd_{(i,j)} := \norm{\cc_{(i,j)}} \lambda_{(i,j)}$ is the so-called graph viscosity coefficient, $\lambda_{(i,j)}$ is the upper bound for the maximum wave speed between nodes $i$ and $j$, $\cc_{(i,j)}$ is an entry of a discrete differentiation operator, and $\NN(i)$ contains all indices of neighbor degrees of freedom of the low-order stencil.

For the particular case of an LGL node distribution and affine geometrical mappings, the discrete differentiation operator for a node can be defined for each coordinate direction using the two-dimensional indexing $i\rightarrow ij$ as \cite{RUEDARAMIREZ2022},
\begin{equation}
    \cc_{(ij,mj)} := \left(\frac{\omega_j \Delta y}{4}, \, 0\right)^T
    \qquad
    \cc_{(ij,im)} := \left(0, \, \frac{\omega_i \Delta x}{4}\right)^T.
\end{equation}
It is also straightforward to define the set of neighbor nodes, $\NN(i)$, using 2D notation: $ij\leftarrow i$.
For an inner DOF, $\NN(ij)$ contains the global indices of the nodes that correspond to nodes (in local 2D notation) $(i-1)j, (i+1)j, i(j-1), i(j+1)$.
For a degree of freedom at the element boundary, $\NN(ij)$ also contains the global indices of the degrees of freedom across the boundary, with which node $ij$ is connected.
For instance, if $i$ ($\rightarrow ij$) is the right-most node on an interface of a Cartesian grid, as sketched in \Cref{fig:SketchConformingInterface}, then $\NN(ij)$ contains the global indices of the nodes that correspond to nodes (in local 2D notation) $\NN(ij) =\{(0j)^+, (i-1)j, i(j-1), i(j+1)\}$, where $(0j)^+$ corresponds to the $0j$ index of the neighbor element on the right.

We will use both notation systems throughout this work. It will be clear from the context and the number of indices used which one is being used.

\NewDocumentCommand{\conforminginterface}{ O{1.2} m }{%
  \centering
  \begin{tikzpicture}[scale=#1,
    declare function={
      x2UR(\x)=\x*0.5+3/2+\shift;
      y2UR(\y)=0.5+0.5*\y;
    }]

    \pgfdeclarelayer{background}
    \pgfsetlayers{background,main}

    \useasboundingbox (-1.2, -1.2) rectangle (3.4, 1.2);

    \def\shift{0.12}
    \def\a{0.4472135955} 

    \def\b{0.6666666667} 


    \draw[thick] (-1,-1) rectangle (1,1);

    \foreach \x in {-\b,0,\b}
    \draw[dashed,gray] (\x,-1) -- (\x,1);

    \foreach \y in {-\b,0,\b}
    \draw[dashed,gray] (-1,\y) -- (1,\y);

    \foreach \x in {-1,-\a,\a,1}
    \foreach \y in {-1,-\a,\a,1}
    \fill (\x,\y) circle (0.04);


    \draw[thick] (1+\shift,-1) rectangle (3+\shift,1);

    \foreach \x in {-\b, 0, \b}
    \draw[dashed,gray] (2+\shift+\x,-1) -- (2+\shift+\x,1);

    \foreach \y in {-\b,0,\b}
    \draw[dashed,gray] (1+\shift,\y) -- (3+\shift,\y);

    \foreach \x in {1+\shift,2+\shift-\a,2+\shift+\a,3+\shift}
    \foreach \y in {-1,-\a,\a,1}
    \fill (\x,\y) circle (0.04);


    \ifnum#2=2
    \fill[red] (1,\a) circle (0.04);
    \foreach \p in {
        (1,1),
        (\a,\a),
        (1,{-\a}),
        ({1+\shift},{\a})
    }{
        \fill[red] \p circle (0.04);
        \draw[red,thick] (1,\a) -- \p;
    }
    \draw[black, line width=0.6pt] (1, \a) circle (0.04);
    \fi
\end{tikzpicture}
}
\begin{figure}[pos=htbp]
    \begin{minipage}{0.65\textwidth}
        \centering
        \conforminginterface[1.7]{2}
    \end{minipage}
    \begin{minipage}[c]{0.3\textwidth}
    \begin{tikzpicture}
          \matrix [draw=black!40, align=center, rounded corners=1pt] at (0,0) {
          \fill (0.45,0.45) circle (0.07);
          \node[anchor=west, font=\small] at (1,0.45) {LGL nodes};

          \draw[black, thick] (0,0) -- (0.9,0);
          \node[anchor=west, font=\small] at (1,0) {Element boundary};

          \draw[dashed, gray] (0,-0.45) -- (0.9,-0.45);
          \node[anchor=west, font=\small] at (1,-0.45) {Subcell boundary};

          \draw[red, thick] (0,-0.9) -- (0.9,-0.9);
          \node[anchor=west, font=\small] at (1,-0.9) {\small Node connections};
          \\
          };
    \end{tikzpicture}
    \end{minipage}
    \caption{Sketch of a conforming interface for a polynomial degree of $N=3$.}
    \label{fig:SketchConformingInterface}
\end{figure}
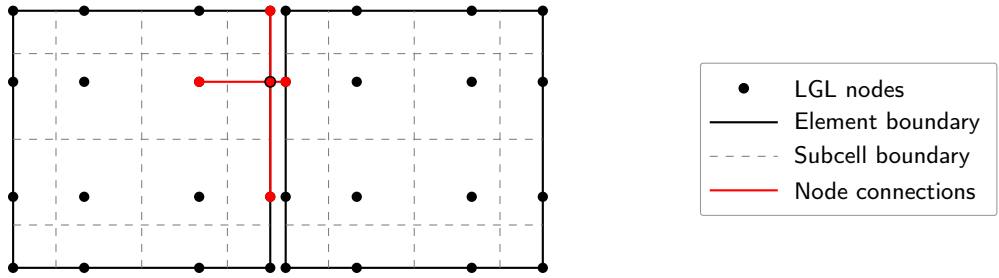

Conforming meshes allow building a subcell LGL grid where every node has $2^d$ neighboring nodes based on the flux-differencing formulation. This holds for inner nodes as well as for nodes on interfaces. For that purpose, one high-order and one low-order version of the flux are constructed and then blended together. This is explained in detail for DGSEM for instance in \cite{RUEDARAMIREZ2022,ruedaramirez2024}. In the following section, we extend the idea of convex limiting to treat nonconforming meshes with hanging nodes. This enables locally refined meshes and adaptive mesh refinement (AMR).

\subsection{Invariant-domain-preserving (IDP) low-order fluxes at nonconforming interfaces}
The standard DG approach for interface fluxes on nonconforming interfaces with hanging nodes -- we call these interfaces mortars -- is to use an $L_2$-projection of the polynomial solution of the large element to the resolution of the smaller ones \cite{KOPRIVA1996475}.
This approach is high-order accurate but not invariant-domain-preserving (IDP), which can cause problems when handling challenging simulation setups. In order to still allow interfaces with hanging nodes and therefore benefit from the efficiency of AMR but without compromising stability, another approach is needed. The goal of this section is to construct such an IDP scheme to compute the interface fluxes at mortars.

\begin{figure}[pos=htbp]
    \begin{minipage}{0.65\linewidth}
    \centering
    \begin{tikzpicture}[scale=1.7,
        declare function={
          x2UR(\x)=\x*0.5+3/2+\shift;
          y2UR(\y)=0.5+0.5*\y;
        }]
        \pgfdeclarelayer{background}
        \pgfsetlayers{background,main}

        \useasboundingbox (-1.1, -1.2) rectangle (2.3, 1.2);

        \def\shift{0.4}
        \def\a{0.4472135955} 

        \def\b{0.6666666667} 


        \draw[thick] (-1,-1) rectangle (1,1);

        \foreach \x in {-\b,0,\b}
        \draw[dashed,gray] (\x,-1) -- (\x,1);

        \foreach \y in {-\b,0,\b}
        \draw[dashed,gray] (-1,\y) -- (1,\y);

        \foreach \x in {-1,-\a,\a,1}
        \foreach \y in {-1,-\a,\a,1}
        \fill (\x,\y) circle (0.04);


        \draw[thick] ({1+\shift},0) rectangle ({2+\shift},1);

        \foreach \x in {-\b,0,\b}
        \draw[dashed,gray] ({\x*0.5+3/2+\shift},0) -- ({\x*0.5+3/2+\shift},1);

        \foreach \y in {-\b,0,\b}
        \draw[dashed,gray] ({1+\shift},\y*0.5+0.5) -- ({2+\shift},\y*0.5+0.5);

        \foreach \x in {-1,-\a,\a,1}
        \foreach \y in {-1,-\a,\a,1}
        \fill ({\shift+\x*0.5+3/2},{0.5+0.5*\y}) circle (0.04);


        \draw[thick] ({1+\shift},-1) rectangle ({2+\shift},0);

        \foreach \x in {-\b,0,\b}
        \draw[dashed,gray] ({\x*0.5+3/2+\shift},-1) -- ({\x*0.5+3/2+\shift},0);

        \foreach \y in {-\b,0,\b}
        \draw[dashed,gray] ({1+\shift},\y*0.5-0.5) -- ({2+\shift},\y*0.5-0.5);

        \foreach \x in {-1,-\a,\a,1}
        \foreach \y in {-1,-\a,\a,1}
        \fill ({\shift+\x*0.5+3/2},{-0.5+0.5*\y}) circle (0.04);


        \draw[thick, blue] (1+0.5*\shift, -1) -- (1+0.5*\shift, 1);
        \node[anchor=south, font=\normalsize, blue] at (1+0.5*\shift,1) {$S$};


        \foreach \y in {-1,-\a,\a,1}
        \fill[green] ({1},\y) circle (0.04);

        \node[anchor=south, font=\scriptsize, green] at (1-0.1, 1) {$\NN_S^-$};

        \foreach \k in {-0.5, 0.5}
        \foreach \y in {-1,-\a,\a,1}
        \fill[orange] ({1+\shift},{\k+0.5*\y}) circle (0.04);

        \node[anchor=south, font=\scriptsize, orange] at (1+\shift+0.1, 1) {$\NN_S^+$};
    \end{tikzpicture}
    \end{minipage}
    \begin{minipage}[c]{0.34\linewidth}
    \begin{tikzpicture}
          \matrix [draw=black!40, align=center, rounded corners=1pt] at (0,0) {
              \fill (0.45,0.45) circle (0.07);
              \node[anchor=west, font=\scriptsize] at (1,0.45) {LGL nodes};

              \draw[black, thick] (0,0) -- (0.9,0);
              \node[anchor=west, font=\scriptsize] at (1,0) {Element boundary};

              \draw[dashed, gray] (0,-0.45) -- (0.9,-0.45);
              \node[anchor=west, font=\scriptsize] at (1,-0.45) {Subcell boundary};

              \\
          };
    \end{tikzpicture}
    \end{minipage}
    \caption{Sketch of a nonconforming interface for a polynomial degree of $N=3$.}
    \label{fig:SketchNonConformingInterface}
\end{figure}
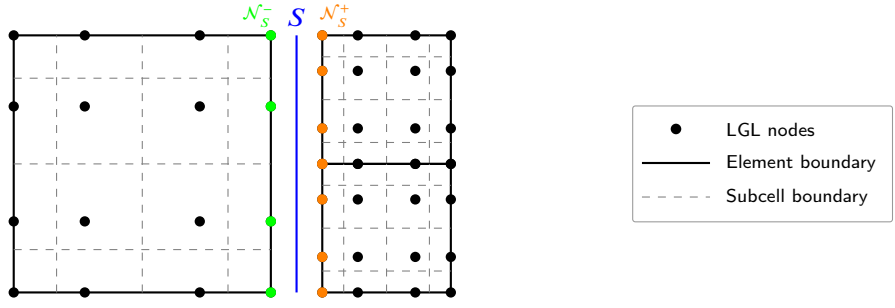

We now consider a nonconforming interface $S\in\partial \Omega_e$ with hanging nodes. In this paper, we focus on nonconforming interfaces where the connected elements use the same polynomial degree and differ only by one level of refinement.
In two spatial dimensions, a nonconforming interface is therefore always a 2-to-1 element connection.
Such an interface is sketched in \Cref{fig:SketchNonConformingInterface}.

We no longer have matching nodes on both sides of $S$.
Therefore, we will distinguish between nodes on both sides of the interface and denote the index sets of nodes on the internal and external sides of interface $S$ by $\NN_S^-$ and $\NN_S^+$, respectively.
For the nonconforming interface shown in \Cref{fig:SketchNonConformingInterface}, $\NN_S^-$ contains the global indices of nodes corresponding to nodes (in local 2D indices) $Nj$ with $j=0,\dots, N$. $\NN_S^+$ contains the global indices of nodes corresponding to $0j$ with $j=0, \dots, N$ of both small elements on the right.

Due to the different sizes of the elements we define the weights of the collocated physical-space LGL quadrature on both sides of $S$
\begin{equation}\label{eq:weightsWithPhi}
    \omega_{i, S}^- = \int_{S} \varphi_i^- \mathrm{d}s, \qquad
    \omega_{j, S}^+ = \int_S \varphi_j^+ \mathrm{d}s, \qquad  i \in\NN_S^-, \, j\in\NN_S^+,
\end{equation}
where $\{\varphi_i^-\}_{i\in\NN_S^-}$ are the basis functions in the left element and $\{\varphi_i^+\}_{i\in\NN_{S}^+}$ in the right elements across $S$.
For tensor-product LGL elements, $\varphi_i$ are 1D Lagrange basis functions.
Note that the integrals are computed exactly using LGL quadrature with $N+1$ nodes.
Since the basis functions on both sides of $S$ sum to unity, we have
\begin{equation}
    \sum_{i\in\NN^-_S} \omega^-_{i,S} = \sum_{j\in\NN^+_S} \omega^+_{j,S}.
\end{equation}

At first, we assume that every node is connected to all nodes on the other side of the interface as in the $L_2$ mortars by Kopriva~\cite{KOPRIVA1996475}. For these all-to-all mortars, we have
\begin{equation}\begin{split}
    &\NN^{+}_S \subset \NN(i) \,\qquad \forall i\in\NN_S^-,\qquad
    \NN^{-}_S \subset \NN(j) \qquad \forall j\in\NN_S^+.
\end{split}\end{equation}

We also distinguish between the outward-pointing normal vector $\vec{n}_S^-$ on the left side of $S$ and $\vec{n}_S^+$ on the right side. They are unit vectors since all physical face measures are contained in the weights. In the case with a Cartesian mesh, the normal vectors are constant on each side of $S$ and $\vec{n}_S^+=-\vec{n}_S^-$.

We now proceed to construct a new mortar flux, which is designed to be conservative and IDP.
The interface flux contribution $\stateL{F}_{i}$ from node $i\in\NN_S^-$ at interface $S$ can be computed with the surface integral \eqref{eq:weakform},
\begin{equation}\label{eq:surface_integral_S}
    \stateL{F}^-_i = \int_{S, N} \varphi_i^- \numfluxb{f}(\stateL{u}^-, \stateL{u}^+, \vec{n}_S^-) \mathrm{d}s,
\end{equation}
with the numerical flux $\numfluxb{f}$, and inner and outer states $\stateL{u}^-$ and $\stateL{u}^+$ with respect to $i$ and $S$. Note that the integrals are considered as numerical integrals using LGL quadrature with $N+1$ nodes, as indicated by the subscript $N$. The quadrature rule is not exact anymore.

The interface contribution of node $j\in\NN_S^+$ on the opposing side reads as,
\begin{equation}
    \stateL{F}_{j}^{+} = \int_{S, N} \varphi_j^+ \numfluxb{f}(\stateL{u}^-, \stateL{u}^+, \vec{n}_S^+) \mathrm{d}s.
\end{equation}

We will construct the mortar flux using nodal quadrature based on a local Lax--Friedrichs type flux containing a central part and a graph viscosity term.
Let
\begin{equation}
    \numfluxb{f}_S^- = \sum_{k \in\NN_S^-} \blocktensor{f}_{k} \cdot \vec{n}_S^- \,\varphi_{k}^-, \qquad
    \numfluxb{f}_S^+ = \sum_{j\in\NN_S^+} \blocktensor{f}_{j} \cdot \vec{n}_S^- \,\varphi_j^+
\end{equation}
be left and right fluxes in a common direction.
Then, the central component of the total LLF flux can be defined as
\begin{equation}
    \numfluxb{f}(\stateL{u}^-, \stateL{u}^+, \vec{n}_S^-) = \frac{1}{2} \left(\numfluxb{f}_S^- + \numfluxb{f}_S^+\right).
\end{equation}
We insert this approach into formula \eqref{eq:surface_integral_S} to obtain an initial all-to-all central (dissipation-free) nodal flux contribution,
\begin{equation}\label{eq:MortarFluxFirstStep}
    \stateL{F}_{i}^-
    = \frac{1}{2} \sum_{k\in\NN_S^-} \blocktensor{f}_{k} \cdot \vec{n}_S^- \int_{S, N} \varphi_i^- \varphi_k^- \mathrm{d}s
    + \frac{1}{2} \sum_{j\in\NN_S^+} \blocktensor{f}_{j} \cdot \vec{n}_S^- \int_{S, N} \varphi_i^- \varphi_j^+ \mathrm{d}s.
\end{equation}

Note that $\int_{S, N} \varphi_i^- \varphi_k^- \mathrm{d}s = \delta_{ik} \omega_{i, S}^-$. Thus, \eqref{eq:MortarFluxFirstStep} simplifies to
\begin{equation}\begin{split}\label{eq:interfaceflux_nonlocal}
    \stateL{F}_{i}^-
    &= \frac{1}{2} \blocktensor{f}_{i} \cdot \vec{n}_S^- \,\omega_{i, S}^- + \frac{1}{2} \sum_{j\in\NN_S^+} \blocktensor{f}_{j} \cdot \vec{n}_S^- \int_{S, N} \varphi_i^- \varphi_j^+ \mathrm{d}s
    = \frac{1}{2} \sum_{j\in\NN_S^+} \omega_{(i, j)}^S \left(\blocktensor{f}_{i} + \blocktensor{f}_{j} \right) \cdot \vec{n}_S^-,
\end{split}\end{equation}
where we introduced the weights
\begin{equation}\label{eq:nonlocal_weight}
    \omega_{(i, j)}^S := \int_{S, N} \varphi_i^- \varphi_j^+ \mathrm{d}s, \qquad i\in\NN_S^-, j\in\NN_S^+
\end{equation}
and used the fact that
\begin{equation}\label{eq:sum_weights}
    \sum_{j\in\NN_S^+} \omega_{(i,j)}^S = \sum_{j\in\NN_S^+} \int_{S, N} \varphi_i^- \varphi_j^+ \mathrm{d}s = \int_S \varphi_i^- \mathrm{d}s = \omega_{i, S}^-,
\end{equation}
due to $\sum_{j\in\NN_S^+}\varphi_j^+ \equiv 1$.

Equivalently, for node $j\in\NN_S^+$, we obtain
\begin{equation}\begin{split}\label{eq:interfaceflux_nonlocal_right}
    \stateL{F}_{j}^{+}
    &= \frac{1}{2} \blocktensor{f}_{j} \cdot \vec{n}_S^+ \,\omega_{j, S}^+ + \frac{1}{2} \sum_{i\in\NN_S^-} \blocktensor{f}_{i} \cdot \vec{n}_S^+ \int_{S, N} \varphi_i^- \varphi_j^+ \mathrm{d}s
    = \frac{1}{2} \sum_{i\in\NN_S^-} \omega_{(i, j)}^S \left(\blocktensor{f}_{i} + \blocktensor{f}_{j}\right) \cdot \vec{n}_S^+.
\end{split}\end{equation}
On a Cartesian mesh, we have constant normal vectors on the left and right side of the interface $S$. Hence, all fluxes in both formulas are computed in the same direction, which simplifies the method significantly.

The constructed fluxes \eqref{eq:interfaceflux_nonlocal} and \eqref{eq:interfaceflux_nonlocal_right} satisfy the conservation relation
\begin{equation}\label{eq:conservationProperty}
    \sum_{i\in\NN_S^-} \stateL{F}_{i}^-
    = -\sum_{j\in\NN_S^+} \stateL{F}_{j}^{+},
\end{equation}
due to $\vec{n}_S^+ = - \vec{n}_S^-$.

Equation \eqref{eq:interfaceflux_nonlocal} defines an all-to-all approach for the nodal interface flux containing only the central part of the LLF without any dissipation. At this point, we can add a graph viscosity term to add dissipation and ensure the IDP property. It was first observed in \cite{anderson2017} that this yields a very diffusive behavior for high-order Bernstein- and Gauss-Lobatto-FE methods. In particular, \cite[Fig. 2]{anderson2017} shows that this behavior is even more pronounced for Gauss-Lobatto than for Bernstein basis polynomials.

One alternative approach is to use a sparsified version of the method. The idea of designing a sparse low-order scheme where the stencil does not grow with the polynomial degree was explored in \cite{lohmann2017,kuzmin2020subcell} for continuous Galerkin and in \cite{hajduk2021} for DG methods using Bernstein polynomials, but also applied to DGSEM in \cite{Pazner2020, RUEDARAMIREZ2022}.

\subsubsection*{Sparsification}
With the approach described above, every node is connected to every node on the other side of the interface using the weight $\omega_{(i, j)}^S$. This lack of sparsity is a potential disadvantage.

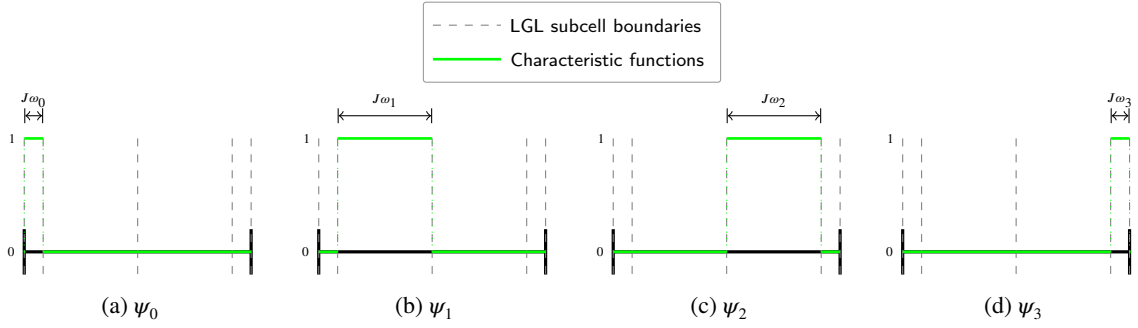
\begin{figure}[pos=htbp]  
  \centering
  \def\xA{-1}
  \def\xB{-5/6}
  \def\xC{0}
  \def\xD{5/6}
  \def\xE{1}

  \def\colorA{green}
  \def\colorB{blue}
  \def\colorC{gray}
  \def\colorD{yellow}

  \tikzset{
    base/.style={black, very thick},
    subcells/.style={gray, dashed}
  }
  \NewDocumentCommand{\basisfunctions}{m m m m}{%
      \begin{tikzpicture}[scale=1.5]

      \draw[base] (\xA, 0) node[black, left] {\tiny$0$} -- (\xE,0);
      \draw[base] (\xA, 0.2) -- (\xA, -0.2);
      \draw[base] (\xE, 0.2) -- (\xE, -0.2);

      \draw[subcells] (\xA, 1) node[black, left] {\tiny$1$} -- (\xA, -0.2);
      \draw[subcells] (\xB, 1) -- (\xB, -0.2);
      \draw[subcells] (\xC, 1) -- (\xC, -0.2);
      \draw[subcells] (\xD, 1) -- (\xD, -0.2);
      \draw[subcells] (\xE, 1) -- (\xE, -0.2);

      \draw[#2, line width=1pt] (\xA, 0) -- (#3, 0);
      \draw[#2, line width=1pt] (#3, 1) -- (#4, 1);
      \draw[#2, line width=1pt] (#4, 0) -- (\xE, 0);

      \draw[dotted, #2] (#3, 0) -- (#3, 1);
      \draw[dotted, #2] (#4, 0) -- (#4, 1);

      \draw[|<->|, line width=0.4pt] (#3, 1.2) -- (#4, 1.2) node[midway, above] {\tiny$J\omega_{#1}$};
    \end{tikzpicture}
    }

    \begin{subfigure}[c]{\textwidth}
    \centering
    \begin{tikzpicture}
    \matrix [draw=black!40, align=center, rounded corners=1pt] at (0,0) {
        \draw[subcells] (0,0) -- (0.9,0) node[anchor=west, black, font=\scriptsize] {LGL subcell boundaries};
        \draw[green, line width=1pt] (0,-0.45) -- (0.9,-0.45) node[anchor=west, black, font=\scriptsize] {Characteristic functions}; \\
        };
    \end{tikzpicture}
\end{subfigure}\\
\begin{subfigure}{0.23\textwidth}
    \centering
    \basisfunctions{0}{green}{\xA}{\xB}
    \caption{$\psi_0$}
\end{subfigure}
\begin{subfigure}{0.23\textwidth}
    \centering
    \basisfunctions{1}{green}{\xB}{\xC}
    \caption{$\psi_1$}
\end{subfigure}
\begin{subfigure}{0.23\textwidth}
    \centering
    \basisfunctions{2}{green}{\xC}{\xD}
    \caption{$\psi_2$}
\end{subfigure}
\begin{subfigure}{0.23\textwidth}
    \centering
    \basisfunctions{3}{green}{\xD}{\xE}
    \caption{$\psi_3$}
\end{subfigure}
\caption{Sketches of characteristic functions of LGL subcells for $N=3$, i.e., $N+1=4$ LGL nodes.}
\label{fig:characteristic_functions}
\end{figure}

To use a more local approach, we replace $\omega_{(i, j)}^S$ from \eqref{eq:nonlocal_weight} with a local weight $\tilde{\omega}_{(i, j)}^S$, which is only nonzero if the nodes are close together. To achieve that, we define
\begin{equation}\label{eq:local_weights}
    \tilde{\omega}_{(i, j)}^S := \int \psi_i^- \psi_j^+ \mathrm{d}s,
\end{equation}
where $\psi_i^-$ and $\psi_j^+$ are the characteristic functions of the LGL subcells around the nodes $i\in\NN_S^-$ and $j\in\NN_S^+$ as sketched in \Cref{fig:characteristic_functions}. These piecewise-constant basis functions take the value of $1$ within the associated subcell and $0$ elsewhere.
These integrals are, in turn, exact integrals, which, however, are easy to evaluate due to the form of the characteristic functions.

We also denote these functions by low-order basis functions.
They are a natural choice for the low-order variant since we also use piecewise constant approximations on the LGL subgrid for the low-order method in the subcell limiting approach in the volume integral \cite{hennemann2021, rueda2021entropy, RR2021}.

Due to the nonnegativity of the characteristic functions, the local weights themselves are nonnegative as well.
Moreover, similar to $\varphi^-$ and $\varphi^+$, we have
\begin{equation}
    \sum_{i\in\NN_S^-} \psi_i^- \equiv 1, \qquad
    \sum_{j\in\NN_S^+} \psi_j^+ \equiv 1,
\end{equation}
which yields
\begin{equation}\begin{split}\label{eq:sum_local_omega}
    \sum_{i\in\NN_S^-} \tilde\omega_{(i, j)}^S &= \int_S \psi_j^+ \mathrm{d}s = \omega_{j, S}^+, \qquad \forall j\in\NN_S^+,\\
    \sum_{j\in\NN_S^+} \tilde\omega_{(i, j)}^S &= \int_{S} \psi_i^- \mathrm{d}s = \omega_{i, S}^-, \qquad \forall i\in\NN_{S}^-,
\end{split}\end{equation}
due to the size of the LGL subcells.

Applying the sparsification approach to \eqref{eq:interfaceflux_nonlocal} yields a local (sparsified) version of the nodal flux for node $i\in\NN_S^-$:
\begin{equation}\begin{split}\label{eq:interfaceflux_central_local}
    \bar{\stateL{F}}_{i}^-
     &= \frac{1}{2} \sum_{j\in \NN_S^+} \tilde{\omega}_{(i, j)}^S \left(\blocktensor{f}_{i}
    + \blocktensor{f}_{j} \right) \cdot \vec{n}_S^-\\
\end{split}\end{equation}
and for the node $j\in\NN_S^+$:
\begin{equation}\begin{split}
    \bar{\stateL{F}}_{j}^{+}
    &= \frac{1}{2} \sum_{i\in \NN_S^-} \tilde{\omega}_{(i, j)}^S \left(\blocktensor{f}_{i}
    + \blocktensor{f}_{j} \right) \cdot \vec{n}_S^+.
\end{split}\end{equation}

These fluxes also satisfy the conservation relation \eqref{eq:conservationProperty}.

Finally, we add the LLF graph viscosity term to ensure the IDP properties. This gives us the low-order nodal interface fluxes
\begin{equation}\label{eq:interfaceflux_local}
    \tilde{\stateL{F}}_{i}^-
    = \frac{1}{2} \sum_{j\in\NN_S^+} \tilde\omega_{(i, j)}^S \left( \left(\blocktensor{f}_{i} + \blocktensor{f}_{j} \right) \cdot \vec{n}_S^-
    - \lambda_{(i, j)} (\stateL{u}_{j} - \stateL{u}_{i})\right), \quad i\in\NN_S^-,
\end{equation}
and
\begin{equation}\begin{split}
    \tilde{\stateL{F}}_{j}^+
    &= \frac{1}{2} \sum_{i\in \NN_S^-} \tilde{\omega}_{(i, j)}^S \left(\left(\blocktensor{f}_{i} + \blocktensor{f}_{j} \right) \cdot \vec{n}_S^+ - \lambda_{(i, j)} (\stateL{u}_{i} - \stateL{u}_{j})\right),\quad j\in\NN_S^+,
\end{split}\end{equation}
where $\lambda_{(i, j)}$ is the upper bound for the maximum wave speed between the left node $i\in\NN_S^-$ and the right node $j\in\NN_S^+$.

Adding the graph viscosity term preserves the conservation relation \eqref{eq:conservationProperty}.
Moreover, the constructed local weights are very natural since when applied to a conforming interface without hanging nodes the flux reduces to the standard conforming-interface flux. In that case, every node is again only connected to its associated direct neighboring node on the other side of the interface due to $\psi_i^+ = \psi_i^- \, \forall i$ and therefore
\begin{equation}\label{eq:mortarfluxforconforminginterface}
    \tilde\omega_{(i, j)}^S = \int_S \psi_i^- \psi_j^+ \mathrm{d}s
    = \delta_{ij} \int_S \psi_i^- \mathrm{d}s
    = \delta_{ij} \omega_{i, S}^-.
\end{equation}

\begin{remark}
The constructed mortar flux \eqref{eq:interfaceflux_local} can also be written in local 2D index notation. Assume an interface $S$ in positive $\xi_1$ direction with one large element on the left side as sketched in \Cref{fig:SketchNonConformingInterface}. The global indices $i\in\NN_S^-$ corresponds to the local indices $Nk$ with $k=0,\dots,N$. For nodes on the right, the global indices $j\in\NN_S^+$ correspond to the local indices $(0\ell)^{e}$ with $\ell=0,\dots, N$ in both small neighboring elements ($e=1, 2$). The local weight $\tilde\omega_{(i, j)}^S$ corresponds to $\tilde\omega_{(Nk, (0\ell)^{e})}^S$. So, \eqref{eq:interfaceflux_local} can be written as follows,
\begin{equation}
    \tilde{\stateL{F}}^-_{Nk}
    = \frac{1}{2} \sum_{e=1}^2 \sum_{\ell=0}^N \tilde\omega_{(Nk, (0\ell)^{e})}^S \left( \left(\blocktensor{f}_{Nk} + \blocktensor{f}_{(0\ell)^{e}} \right) \cdot \vec{n}_S^-
    - \lambda_{(Nk, (0\ell)^{e})} (\stateL{u}_{(0\ell)^e} - \stateL{u}_{Nk})\right).
\end{equation}
\end{remark}

\begin{remark}\label{re:Equivalence_mortarFlux_cdForm}
    Equation~\eqref{eq:interfaceflux_local} can be written as
    \begin{equation}\label{eq:MortarFlux_cd}
        \tilde{\stateL{F}}_i^- = \sum_{j \in \NN_S^+} \cc_{(i,j)} \cdot \left(\blocktensor{f}_{i} + \blocktensor{f}_{j}\right)
        - \sum_{j \in \NN_S^+} \dd_{(i,j)} \left(\stateL{u}_j - \stateL{u}_i \right),
    \end{equation}
    with
    \begin{equation}\label{eq:c_nonconforming}
        \cc_{(i,j)} := \frac{\tilde{\omega}^S_{(i,j)} \vec{n}^-_S}{2}
    \end{equation}
    and again $\dd_{(i,j)} := \norm{\cc_{(i,j)}} \lambda_{(i,j)}$ for $i \in \NN^-_S, \, j \in \NN^+_S$.

    The net contribution of the right-face mortar flux to the semi-discrete residual is
    \begin{equation}\label{eq:MortarFlux_contribution_cd}
        -\tilde{\stateL{F}}_i^- = -\sum_{j \in \NN_S^+} \cc_{(i,j)} \cdot \left(\blocktensor{f}_{i} + \blocktensor{f}_{j}\right)
        + \sum_{j \in \NN_S^+} \dd_{(i,j)} \left(\stateL{u}_j - \stateL{u}_i \right).
    \end{equation}

    All terms involving $\blocktensor{f}_i$ can be canceled out from the sum because of the discrete metric identities,
    \begin{equation}\label{eq:metric_identity}
        \sum_{j \in \NN(i)} \cc_{(i,j)} = \vec{0}.
    \end{equation}
    These identities naturally extend to the nonconforming case, since
    \begin{equation}
        \sum_{j \in \NN_S^+} \cc_{(i,j)}
        = \sum_{j \in \NN_S^+} \frac{1}{2} \tilde{\omega}^S_{(i,j)} \vec{n}^-_S
        = \frac{1}{2} \omega_{i, S}^- \vec{n}^-_S,
    \end{equation}
    where \eqref{eq:sum_local_omega} has been used. This relation then reduces the argument to the conforming case, for instance in \cite[eq. (36)]{RUEDARAMIREZ2022}.
    Consequently, \eqref{eq:MortarFlux_contribution_cd} fits into the low-order update formulation given in \eqref{eq:LowOrderScheme_globalIndices}.
\end{remark}

With the derivation above, we have constructed a stable mortar flux. In the next section, we will prove that this formulation yields an invariant-domain-preserving scheme in the case of nonconforming interfaces.

\subsubsection*{Invariant-domain-preserving (IDP) property}
To prove the IDP property of the constructed surface fluxes, we introduce the so-called auxiliary \textit{bar-states} associated with every pair of nodes $i$ and $j$:
\begin{equation}
    \bar{\stateL{u}}_{(i, j)} := \frac{\stateL{u}_{i} + \stateL{u}_{j}}{2} - \frac{\cc_{(i,j)} \cdot \left(\blocktensor{f}_{j} - \blocktensor{f}_{i}\right)}{2 \dd_{(i, j)}}.
\end{equation}

These bar-states play a central role in the construction of invariant-domain-preserving schemes. They may be interpreted as local intermediate states associated with a one-dimensional Riemann problem with the initial states $\stateL{u}_i$ and $\stateL{u}_j$. The expression corresponds to the mean state generated by a local Lax--Friedrichs, or equivalently HLL-type~\cite{hll1983} approximate Riemann solver.
Consequently, if the nodal states $\stateL{u}_{i}$ and $\stateL{u}_{j}$ belong to a convex invariant set, then the corresponding bar-state $\bar{\stateL{u}}_{(i,j)}$ belongs to the same set. This property was exploited in \cite{Guermond2016} for continuous multi-linear finite element discretizations and follows from the convexity of the invariant set together with the Riemann-solver interpretation of the bar-state. Hence, properties such as positivity preservation and the satisfaction of admissibility constraints can be inferred locally from the analysis of the bar-states.

Since the constructed mortar flux fits the form of the semi-discrete low-order update~\eqref{eq:LowOrderScheme_globalIndices} due to \Cref{re:Equivalence_mortarFlux_cdForm}, the IDP property of the constructed scheme can be shown by rewriting \eqref{eq:LowOrderScheme_globalIndices} into the so-called \textit{bar-state} form, where only a combination of different \textit{bar-states} is used.

Due to the metric identities~\eqref{eq:metric_identity}, we can add
$0=\sum_{j \in \NN(i)} \cc_{(i,j)} \cdot \blocktensor{f}_i$ to show that
\begin{equation}\begin{split}
    \mm_i \dot{\stateL{u}}_i^{\FV}
    &= -\sum_{j \in \NN(i)} \cc_{(i,j)} \cdot \blocktensor{f}_{j}
    + \sum_{j \in \NN(i)} \dd_{(i,j)} \left(\stateL{u}_j - \stateL{u}_i \right)\\
    &= -\sum_{j \in \NN(i)} \cc_{(i,j)} \cdot \left(\blocktensor{f}_{j} - \blocktensor{f}_i\right)
    + \sum_{j \in \NN(i)} \dd_{(i,j)} \left(\stateL{u}_j - \stateL{u}_i \right).
\end{split}\end{equation}
Small rearrangements yield the bar-state form
\begin{equation}\begin{split}\label{eq:barstateform}
    \mm_i \dot{\stateL{u}}_i^{\FV}
    &= -\sum_{j \in \NN(i)} \cc_{(i,j)} \cdot \left(\blocktensor{f}_{j} - \blocktensor{f}_i\right)
    + \sum_{j \in \NN(i)} \dd_{(i,j)} \left(\stateL{u}_j - \stateL{u}_i \right)\\
    &= -\sum_{j \in \NN(i)} \cc_{(i,j)} \cdot \left(\blocktensor{f}_{j} - \blocktensor{f}_i\right)
    + \sum_{j \in \NN(i)} \dd_{(i,j)} \left(\stateL{u}_j + \stateL{u}_i - 2\stateL{u}_i\right)\\
    &= \sum_{j \in \NN(i)} \left[\dd_{(i,j)} \left(\stateL{u}_j + \stateL{u}_i \right) - \cc_{(i,j)} \cdot \left(\blocktensor{f}_{j} - \blocktensor{f}_i\right)\right]
    - \sum_{j \in \NN(i)} 2\dd_{(i,j)} \stateL{u}_i\\
    &= \sum_{j \in \NN(i)} 2\dd_{(i,j)} \left[\frac{\stateL{u}_j + \stateL{u}_i}{2} - \frac{\cc_{(i,j)} \cdot \left(\blocktensor{f}_{j} - \blocktensor{f}_i\right)}{2\dd_{(i,j)}}\right]
    - \sum_{j \in \NN(i)} 2\dd_{(i,j)} \stateL{u}_i\\
    &= \sum_{j \in \NN(i)} 2\dd_{(i,j)} \left(\bar{\stateL{u}}_{(i,j)} - \stateL{u}_i\right).
\end{split}\end{equation}
The semi-discretization contains the node value itself $\stateL{u}_{i} = \bar{\stateL{u}}_{(i, i)}$ and all bar-states between the node and its neighboring nodes $\NN(i)$.

In the conforming case, each node has exactly 4 neighboring nodes (in 2D). The total of five bar-states form the typical 5-point stencil. A node at a nonconforming interface is connected to multiple nodes across the interface based on the sparse weights $\tilde\omega$. Note that due to the nature of the piecewise constant characteristic functions, the local weight $\tilde{\omega}_{(i,j)}^S$ of two nodes is nonzero only if the corresponding LGL subcells overlap. As a result, most of the nodes on the opposite side of the interface do not contribute to $\dot{\stateL{u}}$ due to $\tilde{\omega}_{(i,j)} = 0$ yielding $\dd_{(i,j)} = 0$. \Cref{fig:SketchNonConformingInterface_Connections} sketches active node connections for different nodes across a nonconforming interface. An exemplary computation of the local weights is given in \Cref{sec:local_weights}.
\NewDocumentCommand{\nonconforminginterface}{ O{1.6} m }{%
  \centering
  \begin{tikzpicture}[scale=#1,
    declare function={
      x2UR(\x)=\x*0.5+3/2+\shift;
      y2UR(\y)=0.5+0.5*\y;
    }]

    \pgfdeclarelayer{background}
    \pgfsetlayers{background,main}

    \useasboundingbox (-1.2, -1.2) rectangle (2.4, 1.2);

    \def\shift{0.12}
    \def\a{0.4472135955} 

    \def\b{0.6666666667} 


    \draw[thick] (-1,-1) rectangle (1,1);

    \foreach \x in {-\b,0,\b}
    \draw[dashed,gray] (\x,-1) -- (\x,1);

    \foreach \y in {-\b,0,\b}
    \draw[dashed,gray] (-1,\y) -- (1,\y);

    \foreach \x in {-1,-\a,\a,1}
    \foreach \y in {-1,-\a,\a,1}
    \fill (\x,\y) circle (0.04);


    \draw[thick] ({1+\shift},0) rectangle ({2+\shift},1);

    \foreach \x in {-\b,0,\b}
    \draw[dashed,gray] ({\x*0.5+3/2+\shift},0) -- ({\x*0.5+3/2+\shift},1);

    \foreach \y in {-\b,0,\b}
    \draw[dashed,gray] ({1+\shift},\y*0.5+0.5) -- ({2+\shift},\y*0.5+0.5);

    \foreach \x in {-1,-\a,\a,1}
    \foreach \y in {-1,-\a,\a,1}
    \fill ({\shift+\x*0.5+3/2},{0.5+0.5*\y}) circle (0.04);


    \draw[thick] ({1+\shift},-1) rectangle ({2+\shift},0);

    \foreach \x in {-\b,0,\b}
    \draw[dashed,gray] ({\x*0.5+3/2+\shift},-1) -- ({\x*0.5+3/2+\shift},0);

    \foreach \y in {-\b,0,\b}
    \draw[dashed,gray] ({1+\shift},\y*0.5-0.5) -- ({2+\shift},\y*0.5-0.5);

    \foreach \x in {-1,-\a,\a,1}
    \foreach \y in {-1,-\a,\a,1}
    \fill ({\shift+\x*0.5+3/2},{-0.5+0.5*\y}) circle (0.04);

    \ifnum#2=1
    \begin{pgfonlayer}{background}
        \fill[blue!15] (\b,\b) rectangle (1,1);
    \end{pgfonlayer}
    \fill[red] (1,1) circle (0.04);
    \foreach \x/\y in {
        {1+\shift}/1,
        \a/1,
        1/\a,
        {1+\shift}/{0.5*\a+0.5}
    }{
        \fill[red] (\x,\y) circle (0.04);
        \draw[red,thick] (1,1) -- (\x,\y);
    }
    \draw[red,thick] (1,1) -- (1,1+\shift);
    \draw[black, line width=0.6pt] (1, 1) circle (0.04);

    \else\ifnum#2=2
    \begin{pgfonlayer}{background}
        \fill[blue!15] (\b,0) rectangle (1,\b);
    \end{pgfonlayer}
    \fill[red] (1,\a) circle (0.04);
    \foreach \p in {
        (1,1),
        (\a,\a),
        (1,{-\a}),
        ({1+\shift},{0.5*\a+0.5}),
        ({1+\shift},{0.5*(-\a)+0.5}),
        ({1+\shift},0)
    }{
        \fill[red] \p circle (0.04);
        \draw[red,thick] (1,\a) -- \p;
    }
    \draw[black, line width=0.6pt] (1, \a) circle (0.04);

    \else\ifnum#2=3
    \begin{pgfonlayer}{background}
        \fill[blue!15] (1+\shift,1) rectangle ({x2UR(-\b)},{y2UR(\b)});
    \end{pgfonlayer}
    \fill[red] (1+\shift,1) circle (0.04);
    \foreach \p in {
        (1,1),
        ({x2UR(-1)},{y2UR(\a)}),
        ({x2UR(-\a)},{y2UR(1)})
    }{
        \fill[red] \p circle (0.04);
        \draw[red,thick] ({x2UR(-1)},{y2UR(1)}) -- \p;
    }
    \draw[red,thick] (1+\shift,1) -- (1+\shift,1+\shift);
    \draw[black, line width=0.6pt] (1+\shift,1) circle (0.04);

    \else\ifnum#2=4
    \begin{pgfonlayer}{background}
        \fill[blue!15] ({x2UR(-1)},{y2UR(\b)}) rectangle ({x2UR(-\b)},{y2UR(0)});
    \end{pgfonlayer}
    \fill[red] ({x2UR(-1)},{y2UR(\a)}) circle (0.04);
    \foreach \p in {
        (1,1),
        (1,\a),
        ({x2UR(-1)},1),
        ({x2UR(-\a)},{y2UR(\a)}),
        ({x2UR(-1)},{y2UR(-\a)})
    }{
        \fill[red] \p circle (0.04);
        \draw[red,thick] ({x2UR(-1)},{y2UR(\a)}) -- \p;
    }
    \draw[black, line width=0.6pt] ({x2UR(-1)},{y2UR(\a)}) circle (0.04);

    \else\ifnum#2=5
    \begin{pgfonlayer}{background}
        \fill[blue!15] ({x2UR(-1)},{y2UR(0)}) rectangle ({x2UR(-\b)},{y2UR(-\b)});
    \end{pgfonlayer}
    \fill[red] ({x2UR(-1)},{y2UR(-\a)}) circle (0.04);
    \foreach \p in {
        (1,\a),
        ({x2UR(-1)},{y2UR(\a)}),
        ({x2UR(-1)},{y2UR(-1)}),
        ({x2UR(-\a)},{y2UR(-\a)})
    }{
        \fill[red] \p circle (0.04);
        \draw[red,thick] ({x2UR(-1)},{y2UR(-\a)})-- \p;
    }
    \draw[black, line width=0.6pt] ({x2UR(-1)},{y2UR(-\a)}) circle (0.04);

    \else\ifnum#2=6
    \begin{pgfonlayer}{background}
        \fill[blue!15] ({x2UR(-1)},{y2UR(-\b)}) rectangle ({x2UR(-\b)},{y2UR(-1)});
    \end{pgfonlayer}
    \fill[red] ({x2UR(-1)},{y2UR(-1)}) circle (0.04);
    \foreach \p in {
        (1,\a),
        ({x2UR(-1)},{y2UR(-\a)}),
        ({x2UR(-\a)},{y2UR(-1)})
    }{
        \fill[red] \p circle (0.04);
        \draw[red,thick] ({x2UR(-1)},{y2UR(-1)})-- \p;
    }
    \draw[black, line width=0.6pt] ({x2UR(-1)},{y2UR(-1)}) circle (0.04);
    \fi\fi\fi\fi\fi\fi
\end{tikzpicture}
}
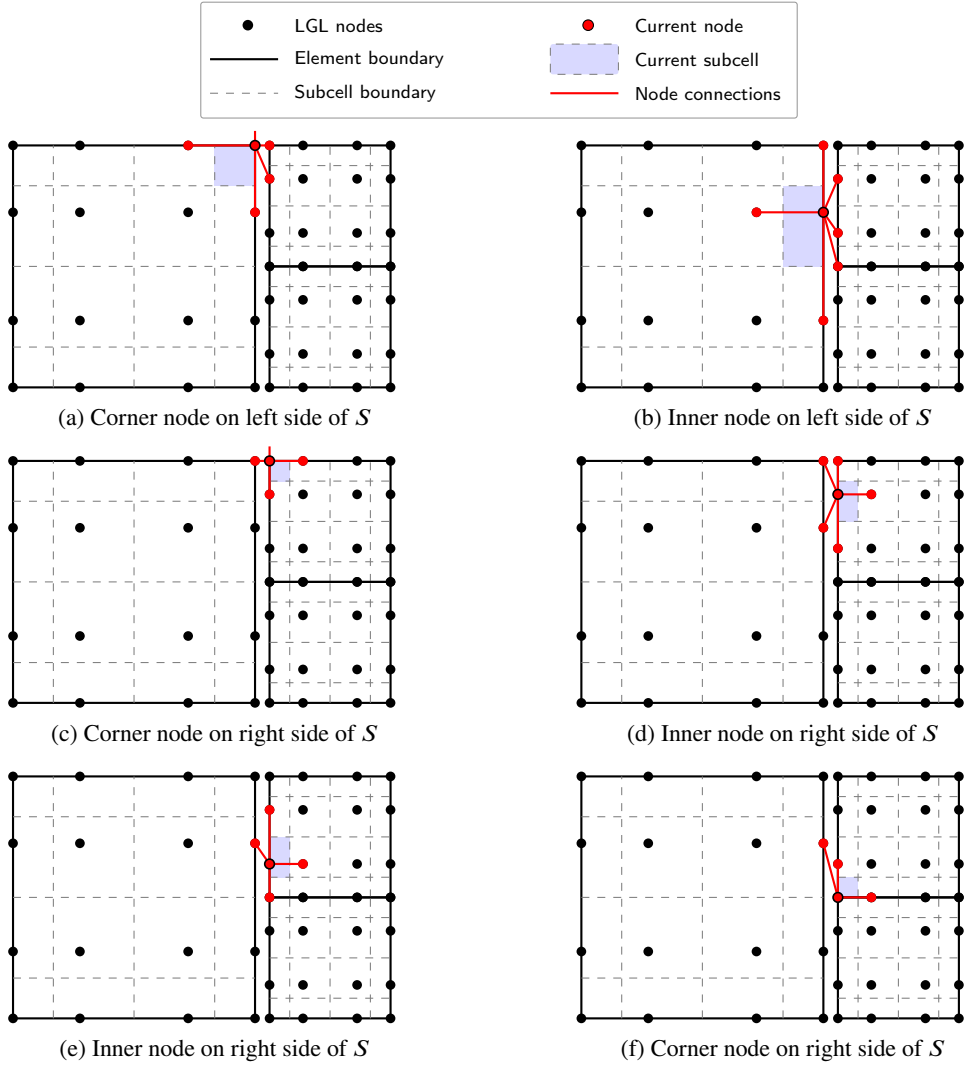
\begin{figure}[pos=htbp]
    \centering
    \begin{tikzpicture}[scale=1.5]
        \matrix [draw=black!40, align=center, rounded corners=1pt] at (0,0) {
        \fill (0.5,0.45) circle (0.07);
        \node[anchor=west, font=\scriptsize] at (1,0.45) {LGL nodes};

        \draw[draw=black, fill=red] (5,0.45) circle (0.07);
        \node[anchor=west, font=\scriptsize] at (5.5,0.45) {Current node};

        \draw[black, thick] (0,0) -- (0.9,0);
        \node[anchor=west, font=\scriptsize] at (1,0) {Element boundary};

        \draw[dashed, gray, fill=blue!15] (4.5, -0.2) rectangle (5.4, 0.2);
        \node[anchor=west, font=\scriptsize] at (5.5,0) {Current subcell};

        \draw[dashed, gray] (0,-0.45) -- (0.9,-0.45);
        \node[anchor=west, font=\scriptsize] at (1,-0.45) {Subcell boundary};

        \draw[red, thick] (4.5,-0.45) -- (5.4,-0.45);
        \node[anchor=west, font=\scriptsize] at (5.5,-0.45) {Node connections};\\
        };
    \end{tikzpicture}\\
    \begin{subfigure}{0.45\textwidth}
        \centering
        \nonconforminginterface{1}
        \vspace{-3mm}
        \caption{Corner node on left side of $S$}
    \end{subfigure}
    \begin{subfigure}{0.45\textwidth}
        \centering
        \nonconforminginterface{2}
        \vspace{-3mm}
        \caption{Inner node on left side of $S$}
    \end{subfigure}
    \begin{subfigure}{0.45\textwidth}
        \centering
        \nonconforminginterface{3}
        \vspace{-3mm}
        \caption{Corner node on right side of $S$}
    \end{subfigure}
    \begin{subfigure}{0.45\textwidth}
        \centering
        \nonconforminginterface{4}
        \vspace{-3mm}
        \caption{Inner node on right side of $S$}
    \end{subfigure}
    \begin{subfigure}{0.45\textwidth}
        \centering
        \nonconforminginterface{5}
        \vspace{-3mm}
        \caption{Inner node on right side of $S$}
    \end{subfigure}
    \begin{subfigure}{0.45\textwidth}
        \centering
        \nonconforminginterface{6}
        \vspace{-3mm}
        \caption{Corner node on right side of $S$}
    \end{subfigure}
    \caption{Illustration of node connections using the sparse approach for a polynomial degree of $N=3$.}
    \label{fig:SketchNonConformingInterface_Connections}
\end{figure}

Using a forward Euler method for the time integration yields a fully discrete version of the update scheme,
\begin{equation}\label{eq:convex_combination}
    \stateL{u}_{i}^{(n+1)}
    = \left(1 - \frac{\Delta t}{\mm_i} \sum_{j \in \NN(i)} 2\dd_{(i,j)}\right) \, \stateL{u}_i^{(n)}
    + \frac{\Delta t}{\mm_i} \sum_{j \in \NN(i)} 2\dd_{(i,j)}\, \bar{\stateL{u}}_{(i,j)}^{(n)}.
\end{equation}

To finish the proof of invariant domain preservation, we need the new solution to be a convex combination of the old solution and the bar-states. This holds under the following time step restriction,
\begin{equation}\label{eq:timestep_restriction}
    \frac{\Delta t}{\mm_i} \sum_{j \in \NN(i)} 2\dd_{(i,j)} \leq 1
    \quad \Longleftrightarrow \quad
    \Delta t \leq \frac{\mm_i}{\sum_{j \in \NN(i)} 2\dd_{(i,j)}}.
\end{equation}
We refer to this as the IDP time step restriction.

With that, the proof is complete.

\begin{remark}
    The proof uses a forward Euler step for the time integration. In order to achieve a higher order time integration, we use a strong-stability preserving (SSP) Runge--Kutta (RK) method which can be written as a convex combination of forward Euler steps.
\end{remark}

\begin{remark}
    The time step restriction \eqref{eq:timestep_restriction} generalizes the one used before in the conforming case \cite{kuzmin2022limiter} to nonconforming meshes using the new definition of $\dd_{(i,j)}$.
\end{remark}

\begin{remark}
Equation \eqref{eq:timestep_restriction} can also be written with 2D local indices. Assuming that node $i\rightarrow ij$ lies on one interface in positive $\xi_1$-direction (with corresponding 2D local indices $Nj$) we have,
\begin{equation}\begin{split}\label{eq:time_restriction_2dlocalindices}
    \frac{\Delta t}{J}
    \left( \frac{1}{\omega_{N}} \left(\lambda_{((N-1)j, Nj)} + \tilde{\lambda}_{(Nj, (Nj)^+)} \right)
    + \frac{1}{\omega_{j}} \left(\lambda_{(N(j-1), Nj)} + \lambda_{(Nj, N(j+1))} \right) \right)
    \leq 1,
\end{split}\end{equation}
where
\begin{equation}
    \tilde{\lambda}_{(Nj, (Nj)^+)} := \frac{1}{\omega_{j, S}^-} \sum_{\ell\in\NN_S^+} \tilde\omega_{(j, \ell)} \lambda_{(Nj, 0\ell)}.
\end{equation}
This formulation naturally extends to nodes on more than one interface. \eqref{eq:time_restriction_2dlocalindices} uses the same notation as in \cite[eq. (44)]{RUEDARAMIREZ2022} making it clear that it is a generalization of the previously used time step restriction for nonconforming meshes.
\end{remark}

\begin{remark}
    Unfortunately, the IDP time step restriction is more restrictive than the typical stability CFL condition of the low-order method, which scales as
    \begin{equation}\label{eq:CFLcondition_lowordermethod}
         \Delta t \lesssim \min_i
         \left(
         \min_{j\in\{1,\dots,d\}} \left(
         \frac{\Delta x^j_i}{\lambda^{j}_i}
         \right) \right),
    \end{equation}
    where $\Delta x_i^j$ is the size of the subcell $i$ in direction $j\in\{1,\dots,d\}$.
    While \eqref{eq:timestep_restriction} is required to have a provable IDP scheme, in practice, we can use an approach to circumvent the strict IDP CFL condition and still ensure set bounds by avoiding the bar-state bounds completely. This approach was already described and used in \cite{ruedaramirez2024}. Instead, we can use bounds that are computed based on the (stable) low-order solution after the forward Euler step but before the correction stage. Such a strategy is especially useful in pure positivity-preserving simulations that do not make use of local bounds to reduce oscillations, since the required minimum bounds do not use the bar-states anyway.
    Nevertheless, not being provably IDP can cause stability issues, which in practice often can be avoided by decreasing the CFL number. In total, despite the smaller CFL number, this approach can reduce the computing time significantly in some situations.
    We compare the results and time step sizes of both CFL conditions in \Cref{sec:sedovblast}.
\end{remark}

\subsection{Mortar Limiting}\label{sec:MortarLimiting}
In the previous section, we constructed a provably invariant-domain-preserving approach to compute the mortar flux $\tilde{\stateL{F}}$. Due to its low-order accuracy, we will use it only where it is absolutely needed to stabilize the simulation and to ensure compliance with certain properties. Elsewhere, we want to use the standard high-order accurate mortar flux based on the $L_2$ projection~\cite{KOPRIVA1996475}, which we will denote $\stateL{F}^{L_2}$.

The general idea is to blend the stable but low-order accurate and the high-order accurate fluxes at the mortar with a limiting factor $\alpha_S\in[0,1]$ as follows:
\begin{equation}\begin{split}\label{eq:BlendingFluxes_alpha}
    \stateL{F}^{\text{Lim}} &= \alpha_S \tilde{\stateL{F}} + \left(1 - \alpha_S\right) \stateL{F}^{L_2}
\end{split}\end{equation}

So, $\alpha_S=0$ implies that we only use the high-order flux $\stateL{F}^{L_2}$, while we only use the low-order flux $\tilde{\stateL{F}}$ if $\alpha_S=1$. Note that we limit on the mortar level with one limiting factor per mortar to preserve conservation.
In the present work, we use an \textit{a posteriori} flux-corrected transport-type (FCT) limiting approach in the volume integral as well as at the mortars \cite{RUEDARAMIREZ2022}. First, the pure low-order fluxes are used. Subsequently, in a correction stage, stable limiting factors are computed and then a stability-preserving portion of the high-order correction is applied - first in the volume integral and then at the mortars. With this approach the resulting flux can be written as
\begin{equation}
    \stateL{F}^\text{Lim} = \tilde{\stateL{F}} + (1 - \alpha_S) (\stateL{F}^{L_2} - \tilde{\stateL{F}}).
\end{equation}
In this work, we consider two sets of properties that we enforce. First, the setup that uses the least amount of limiting only ensures positivity of density and pressure. We follow \cite{RUEDARAMIREZ2022} and apply slightly stricter non-negative lower bounds depending on the low-order solution,
\begin{equation}\label{eq:Bounds_Positivity}
    \rho_i \geq \beta \rho_i^\FV, \qquad
    p_i \geq \beta p_i^\FV
\end{equation}
with $\beta=0.1$. We will refer to this as ``positivity limiting''.

Secondly, to suppress oscillations, we additionally apply local minimum and/or maximum principles for the density and entropy.
Typical bounds are
\begin{equation}\label{eq:Bounds_Barstates}
    \min_{j\in\NN(i)} \bar\rho_{(i,j)} \leq \rho_i \leq \max_{j\in\NN(i)} \bar\rho_{(i,j)}, \qquad
    \min_{j\in \NN(i)} \eta(\bar{\stateL{u}}_{(i,j)}) \leq \eta(\stateL{u}_i),
\end{equation}
where $\eta = \frac{p}{\gamma - 1} \rho^{-\gamma}$ is a modified specific entropy introduced by Guermond et al.~\cite{guermond2019}. To ensure positivity, we always apply the positivity limiting as well.
We call this local limiting due to the locality of bounds.

These minimum and maximum bounds are applied in the limiting for the volume integral and in the mortars.
For the computation of the limiting factors in the volume integral we follow the procedure described in \cite{RUEDARAMIREZ2022}. The limiting factors $\alpha_S$ at the mortars can be computed similarly.
\begin{remark}
    In practice, there is a significant difference between limiting at the mortars and in the volume integral. In the volume integral, all differences between high- and low-order fluxes, i.e., the so-called anti-diffusive fluxes, for one node are summed up when computing the limiting factor for conservative variables (see \cite[eq. (28)]{RUEDARAMIREZ2022}) using a Zalesak-type limiter~\cite{zalesak1979}. This is not possible at the mortars, since for each mortar the limiting factor is computed independently for each mortar.
    For each node $i$, the admissible correction budget is divided among all mortar contributions incident on $i$. If $n_i$ denotes the number of active mortar contributions at node $i$, we limit each mortar correction by scaling the anti-diffusive flux by $n_i$, thereby simulating the effect as if each mortar contributed this flux individually.
    This guarantees that the sum of all accepted mortar corrections remains within the original nodal admissible interval.
    A similar approach was used in \cite[eq. (29)]{RUEDARAMIREZ2022} to limit nonlinear variables using the factor $\Gamma=2d$.
\end{remark}

\subsection{Stabilized solution transfer during refinement and coarsening}\label{sec:LimitingTransfer}
Even if the time integration update using the constructed mortar fluxes and a stable limiting factor is invariant-domain-preserving, stability issues can occur when the solution needs to be transferred from one mesh to another.
A common practice uses interpolation in the refinement step and an $L_2$ projection for the coarsening step~\cite{ranocha2022adaptive}.
These transfer operators are not IDP and can produce non-physical values as well.
In our simulations, applying a positivity-preserving limiter was sufficient to stabilize the transfer step. We used a Zhang--Shu type scaling limiter~\cite{zhang2011} after every refinement and coarsening step on affected elements. The limiter rescales the conservative nodal states toward a conservative mean state when necessary.

The variables used to detect non-admissible states can be chosen flexibly. In this work, an element is marked as troubled if the density or pressure falls below a prescribed global threshold, e.g., $10^{-10}$, at least at one node.
A scaling factor is then chosen such that the corrected conservative states satisfy the prescribed lower bounds for density and pressure at all nodes.
To ensure that the reference state used by the scaling limiter is admissible, we use the mean state of the parent element before refinement as the reference state for the affected child elements.
During coarsening, the scaling limiter is applied using the mean state obtained from the coarsening projection.

We do not claim that this transfer operator is invariant-domain preserving in the strict sense. Rather, it is a conservative positivity-stabilizing post-processing step that was sufficient to prevent non-admissible states in the simulations considered here.

\section{Numerical results}\label{sec:numerical_results}
In this section, we assess the convergence and robustness of the proposed mortar treatment for the compressible Euler equations of gas dynamics with a heat capacity ratio of $\gamma=1.4$.

In all cases, we use the local Lax--Friedrichs (Rusanov) numerical flux for the compatible robust low-order subcell scheme in the volume integral as well as for the surface fluxes in the high-order DG method. Unless stated otherwise, we use the entropy-conserving and kinetic energy preserving flux of Ranocha \cite{ranocha2018generalised} for the two-point volume numerical flux of the split-form DGSEM method.

We use the third-order explicit strong stability preserving (SSP) Runge--Kutta (RK) method as numerical time integration method \cite{SHU1988439}.

As mentioned before in \Cref{sec:MortarLimiting}, we distinguish between two types of limiting. Positivity limiting enforces admissibility of density and pressure. Local limiting additionally imposes local bounds for density and entropy to suppress oscillations for numerical admissibility.

In some of the following examples, we illustrate the amount of limiting in the volume integral and at the mortars. As in \cite{RUEDARAMIREZ2022}, we use a weighted average of the limiting factors across all elements and mortars, respectively. $\alpha=1$ denotes pure low-order mortar fluxes, while $\alpha=0$ denotes the high-order flux.

As mentioned before, all results are conducted using the code framework \texttt{Trixi.jl}~\cite{ranocha2022adaptive, schlottkelakemper2021purely, schlottkelakemper2025trixi}.

\subsection{Advection of a Density Wave}
To test the convergence property of a scheme using the presented low-order mortar flux, we simulate the advection of a density wave with initial condition
\begin{equation}
\label{eq:entropywave2D}
\rho(\vec{x}) = 2+A\,\sin(2\,\pi\,(x+y)),\quad A=0.98,\quad \vec{v} = (0.1,0.2,0),\quad p=20,
\end{equation}
in the periodic domain $\Omega=[-1,1]^2$ until $T=1$. To get a nonconforming mesh we manually add a refined box in $[-0.5, 0]\times [-0.5, 0.5]$.

Our mesh implementation is based on a quad-tree approach~\cite{BursteddeWilcoxGhattas11}, where the leaves represent the elements. The coarsest mesh contains only one element, i.e., one node. It has refinement level 0. Every element can be divided into $2^d$ child elements individually by appending $2^d$ child nodes to the tree node. The new elements are one refinement level higher than the parent element. In total, a mesh can have at most $2^n$ elements per dimension of refinement level $n$.

\begin{figure}[pos=H]
    \centering
    \includegraphics[trim=450 0 450 0, clip, width=0.5\linewidth]{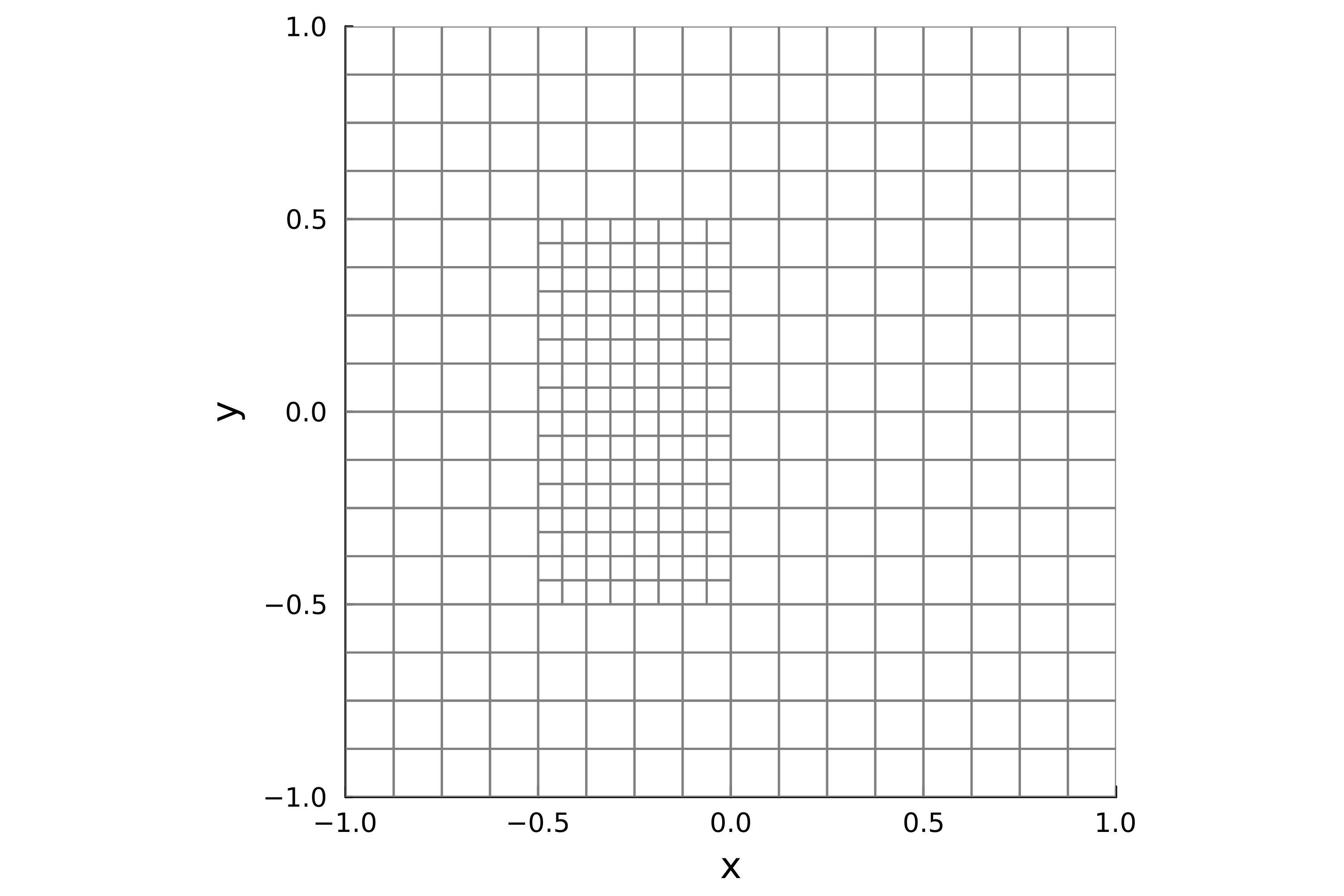}
    \caption{Nonconforming mesh with base refinement level 4 used for the convergence tests.}
    \label{fig:RefinedMesh_Box}
\end{figure}

The mesh used for the convergence test is divided into elements of two sizes. While most of the elements have one refinement level, the so-called base level, the elements in the refined box are one level higher.
An exemplary mesh with base refinement level 4 (corresponding to $2^4=16$ elements per dimension) is shown in \Cref{fig:RefinedMesh_Box}.

We quantify the $L_2$ error of the solution as the mesh is further refined while the refined box is refined once more in every iteration. We perform convergence tests with polynomial degrees $N=3$ and $N=4$. We use $\text{CFL}=0.95$ with the IDP time step restriction~\eqref{eq:timestep_restriction} in the following simulations.
For the first simulation we enabled positivity limiting for density and pressure.

Because the density and pressure remain far from the admissibility positivity bounds, no limiting is expected and therefore a pure high-order convergence rate. The results with polynomial degree of $3$ are shown in \Cref{table:ConvergenceTestPositivity}.
\begin{table}[pos=htbp]
    \centering
    \resizebox{0.85\textwidth}{!}{%
    \begin{tabular}{c|cc|cc|cc|cc}
    \toprule
     & \multicolumn{2}{c}{\(\rho\)}
     & \multicolumn{2}{c}{\(\rho v_1\)}
     & \multicolumn{2}{c}{\(\rho v_2\)}
     & \multicolumn{2}{c}{\(\rho e\)} \\
    level & error & EOC & error & EOC & error & EOC & error & EOC \\
    \midrule
    2(+1) & $7.87 \times 10^{-2}$ & -    & $7.87 \times 10^{-3}$ & -    & $1.57 \times 10^{-2}$ & -    & $1.97 \times 10^{-3}$ & -    \\
    3(+1) & $9.09 \times 10^{-3}$ & 3.11 & $9.09 \times 10^{-4}$ & 3.11 & $1.82 \times 10^{-3}$ & 3.11 & $2.27 \times 10^{-4}$ & 3.11 \\
    4(+1) & $6.67 \times 10^{-4}$ & 3.77 & $6.67 \times 10^{-5}$ & 3.77 & $1.33 \times 10^{-4}$ & 3.77 & $1.67 \times 10^{-5}$ & 3.77 \\
    5(+1) & $2.19 \times 10^{-5}$ & 4.93 & $2.19 \times 10^{-6}$ & 4.93 & $4.38 \times 10^{-6}$ & 4.93 & $5.48 \times 10^{-7}$ & 4.93 \\
    6(+1) & $1.09 \times 10^{-6}$ & 4.32 & $1.09 \times 10^{-7}$ & 4.32 & $2.19 \times 10^{-7}$ & 4.32 & $2.73 \times 10^{-8}$ & 4.32 \\
    7(+1) & $7.21 \times 10^{-8}$ & 3.92 & $7.21 \times 10^{-9}$ & 3.92 & $1.44 \times 10^{-8}$ & 3.92 & $1.80 \times 10^{-9}$ & 3.92 \\
    \bottomrule
    \end{tabular}}
    \caption{$L_2$ errors and EOC for convergence test with positivity limiting. Each mesh consists of elements of the given refinement level with an additionally refined box. The polynomial degree is $3$.}
    \label{table:ConvergenceTestPositivity}
\end{table}

As expected, the results are purely high-order with an experimental order of convergence (EOC) of $N+1=4$. In fact, the results without any enabled limiting are identical. To test the IDP mortar fluxes, we do the convergence test once more with enabled local limiting.
The results are shown in \Cref{table:ConvergenceTestLocalLimiting}.

\begin{table}[pos=htbp]
    \centering
    \resizebox{0.85\textwidth}{!}{%
    \begin{tabular}{c|cc|cc|cc|cc}
    \toprule
     & \multicolumn{2}{c|}{\(\rho\)}
     & \multicolumn{2}{c|}{\(\rho v_1\)}
     & \multicolumn{2}{c|}{\(\rho v_2\)}
     & \multicolumn{2}{c}{\(\rho e\)} \\
    level & error & EOC & error & EOC & error & EOC & error & EOC \\
    \midrule
    2(+1) & $5.43 \times 10^{-1}$ & -    & $5.43 \times 10^{-2}$ & -    & $1.09 \times 10^{-1}$ & -    & $1.36 \times 10^{-2}$ & -    \\
    3(+1) & $1.69 \times 10^{-1}$ & 1.68 & $1.69 \times 10^{-2}$ & 1.68 & $3.38 \times 10^{-2}$ & 1.68 & $4.22 \times 10^{-3}$ & 1.68 \\
    4(+1) & $6.86 \times 10^{-2}$ & 1.30 & $6.86 \times 10^{-3}$ & 1.30 & $1.37 \times 10^{-2}$ & 1.30 & $1.72 \times 10^{-3}$ & 1.30 \\
    5(+1) & $3.37 \times 10^{-2}$ & 1.03 & $3.37 \times 10^{-3}$ & 1.03 & $6.73 \times 10^{-3}$ & 1.03 & $8.42 \times 10^{-4}$ & 1.03 \\
    6(+1) & $1.72 \times 10^{-2}$ & 0.97 & $1.72 \times 10^{-3}$ & 0.97 & $3.44 \times 10^{-3}$ & 0.97 & $4.31 \times 10^{-4}$ & 0.97 \\
    7(+1) & $8.80 \times 10^{-3}$ & 0.97 & $8.80 \times 10^{-4}$ & 0.97 & $1.76 \times 10^{-3}$ & 0.97 & $2.20 \times 10^{-4}$ & 0.97 \\
    \bottomrule
    \end{tabular}}
    \caption{$L_2$ errors and EOC for convergence test with local limiting. Each mesh consists of elements of the given refinement level with an additionally refined box. The polynomial degree is $3$.}
    \label{table:ConvergenceTestLocalLimiting}
\end{table}

In this case, the experimental order of convergence is about $1$. This was expected due to the use of local non-oscillatory bounds, which lead to a higher amount of the low-order scheme.

The results of the convergence tests with a polynomial degree of $4$ are shown in \Cref{table:ConvergenceTest_polydeg4}.
\begin{table}[pos=htbp]
    \begin{subtable}{\textwidth}
    \centering
    \resizebox{0.85\textwidth}{!}{%
    \begin{tabular}{c|cc|cc|cc|cc}
    \toprule
     & \multicolumn{2}{c}{\(\rho\)}
     & \multicolumn{2}{c}{\(\rho v_1\)}
     & \multicolumn{2}{c}{\(\rho v_2\)}
     & \multicolumn{2}{c}{\(\rho e\)} \\
    level & error & EOC & error & EOC & error & EOC & error & EOC \\
    \midrule
    2(+1) & $1.45 \times 10^{-2}$ & -    & $1.45 \times 10^{-3}$ & -    & $2.90 \times 10^{-3}$ & -    & $3.63 \times 10^{-4}$ & -    \\
    3(+1) & $1.02 \times 10^{-3}$ & 3.82 & $1.02 \times 10^{-4}$ & 3.82 & $2.05 \times 10^{-4}$ & 3.82 & $2.56 \times 10^{-5}$ & 3.82 \\
    4(+1) & $2.82 \times 10^{-5}$ & 5.18 & $2.82 \times 10^{-6}$ & 5.18 & $5.64 \times 10^{-6}$ & 5.18 & $7.05 \times 10^{-7}$ & 5.18 \\
    5(+1) & $1.28 \times 10^{-6}$ & 4.46 & $1.28 \times 10^{-7}$ & 4.46 & $2.56 \times 10^{-7}$ & 4.46 & $3.20 \times 10^{-8}$ & 4.46 \\
    6(+1) & $6.96 \times 10^{-8}$ & 4.20 & $6.96 \times 10^{-9}$ & 4.20 & $1.39 \times 10^{-8}$ & 4.20 & $1.74 \times 10^{-9}$ & 4.20 \\
    7(+1) & $3.69 \times 10^{-9}$ & 4.24 & $3.69 \times 10^{-10}$ & 4.24 & $7.38 \times 10^{-10}$ & 4.24 & $9.22 \times 10^{-11}$ & 4.24 \\
    \bottomrule
    \end{tabular}}
    \subcaption{Simulation with positivity limiting.}
    \end{subtable}

    \begin{subtable}{\textwidth}
    \centering
    \resizebox{0.85\textwidth}{!}{%
    \begin{tabular}{c|cc|cc|cc|cc}
    \toprule
     & \multicolumn{2}{c|}{\(\rho\)}
     & \multicolumn{2}{c|}{\(\rho v_1\)}
     & \multicolumn{2}{c|}{\(\rho v_2\)}
     & \multicolumn{2}{c}{\(\rho e\)} \\
    level & error & EOC & error & EOC & error & EOC & error & EOC \\
    \midrule
    2(+1) & $2.94 \times 10^{-1}$ & -    & $2.94 \times 10^{-2}$ & -    & $5.88 \times 10^{-2}$ & -    & $7.35 \times 10^{-3}$ & -    \\
    3(+1) & $1.06 \times 10^{-1}$ & 1.47 & $1.06 \times 10^{-2}$ & 1.47 & $2.12 \times 10^{-2}$ & 1.47 & $2.65 \times 10^{-3}$ & 1.47 \\
    4(+1) & $3.88 \times 10^{-2}$ & 1.45 & $3.88 \times 10^{-3}$ & 1.45 & $7.77 \times 10^{-3}$ & 1.45 & $9.71 \times 10^{-4}$ & 1.45 \\
    5(+1) & $1.67 \times 10^{-2}$ & 1.21 & $1.67 \times 10^{-3}$ & 1.21 & $3.35 \times 10^{-3}$ & 1.21 & $4.19 \times 10^{-4}$ & 1.21 \\
    6(+1) & $8.96 \times 10^{-3}$ & 0.90 & $8.96 \times 10^{-4}$ & 0.90 & $1.79 \times 10^{-3}$ & 0.90 & $2.24 \times 10^{-4}$ & 0.90 \\
    7(+1) & $4.89 \times 10^{-3}$ & 0.87 & $4.89 \times 10^{-4}$ & 0.87 & $9.78 \times 10^{-4}$ & 0.87 & $1.22 \times 10^{-4}$ & 0.87 \\
    \bottomrule
    \end{tabular}}
    \subcaption{Simulation with local limiting.}
    \end{subtable}
    \caption{$L_2$ errors and EOC for convergence tests. Each mesh consists of elements of the given refinement level with an additionally refined box. The polynomial degree is $4$.}
    \label{table:ConvergenceTest_polydeg4}
\end{table}

\subsection{Smooth Isentropic Flow}
We can test the convergence property of our scheme with a smooth solution that actually requires positivity limiting. We use the two-dimensional version of the isentropic flow convergence test which was introduced by Cheng and Shu \cite{cheng2014positivity} with the initial condition
\begin{align}
    \rho(x') = 1 + 0.9999999 \sin(\pi x'),
    ~~
    v_1 = v_2 = 0,
    ~~
    p(x') = \rho(x')^{\gamma},
\end{align}
where we choose $\gamma = 3$. To obtain a 2D version we use $x'$ to apply a rotation with the angle $\theta = 45^\circ$,
\begin{equation}
x' = x \cos(\theta) + y \sin(\theta).
\end{equation}
To get a periodic initial condition we use the domain $[-\sqrt{2}, \sqrt{2}]^2$.

The important differences to the prior setup are the reduced minimum density and pressure, as in \cite{vilar2019posteriori}. The minimum density at the beginning of the simulation is $\rho (x'_0,t=0) = 10^{-7}$ and the minimum pressure is $p (x'_0,t=0) = 10^{-21}$ at both $x'_{0} = -0.5$ and $x'_{0} = 1.5$.

The flow is isentropic and $p=\rho^\gamma$ is preserved during the simulation. This causes the compressible Euler equations to reduce to two Burgers' equations. The analytical solution is given by \cite{bacigaluppi2023posteriori},
\begin{equation}\begin{split}
    \rho(x',t) &= \frac{1}{2} \left( \rho_0(x_1) + \rho_0(x_2) \right),\qquad
    v_{\theta}(x',t) = \sqrt{3} \left( \rho(x', t) - \rho_0(x_1) \right), \\
    p(x',t) &= \rho(x', t)^{\gamma},
\end{split}\end{equation}
where for each coordinate $x'$ and time $t$, $x_1$ and $x_2$ are solutions of the nonlinear equations
\begin{equation}\begin{split}\label{eq:isentropicflow_nonlinearequations}
    x' + \sqrt{3} \rho_0(x_1) t - x_1 &= 0, \\
    x' - \sqrt{3} \rho_0(x_2) t - x_2 &= 0.
\end{split}\end{equation}
The velocity components are simply obtained as
$$\vec{v}(\vec{x},t) = (v_1, v_2) = (v_{\theta} \cos(\theta), v_{\theta} \sin(\theta)).$$

\begin{figure}[pos=htbp]
    \centering
    \begin{subfigure}{0.45\textwidth}
        \includegraphics[width=\textwidth]{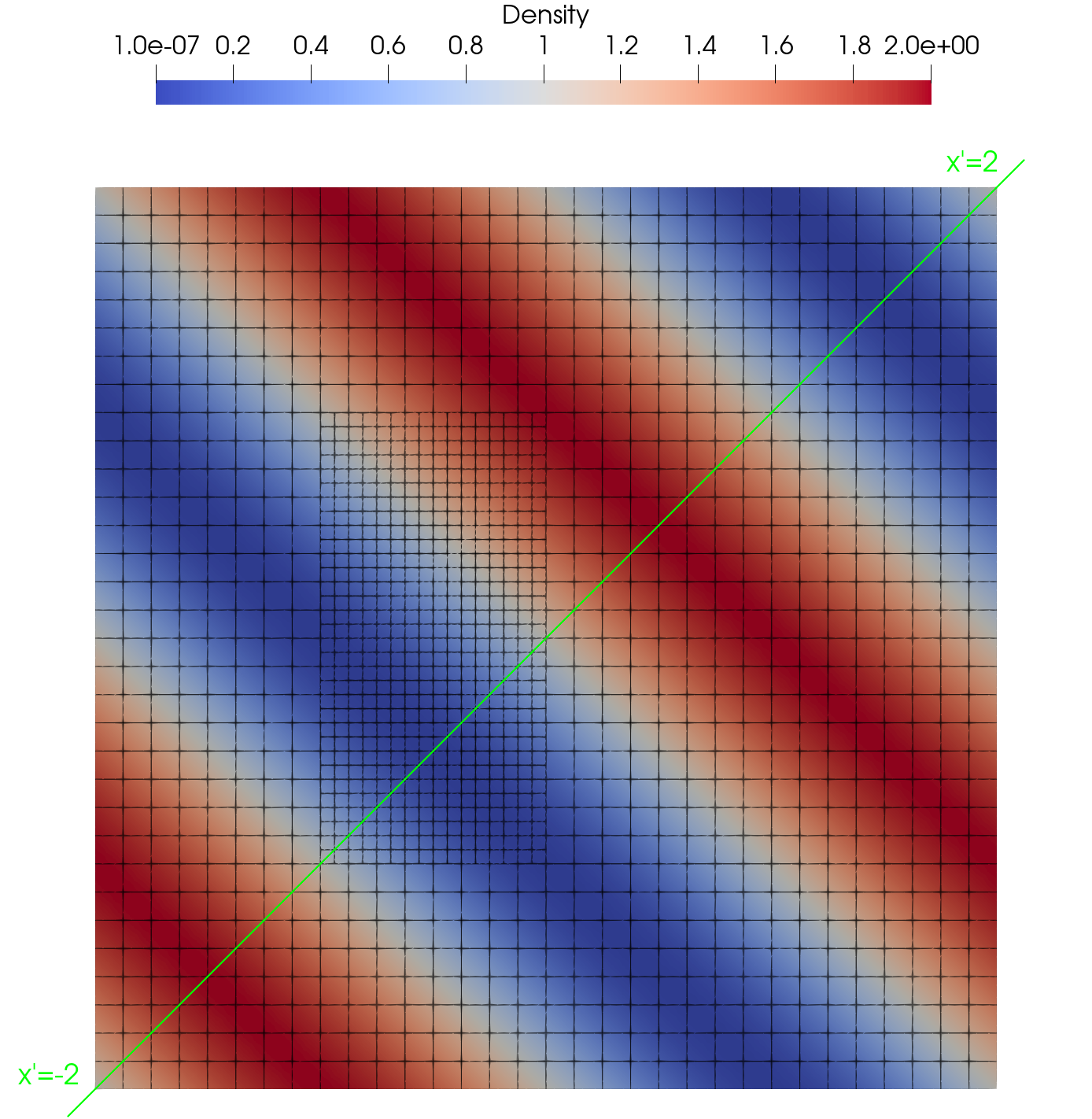}
    \end{subfigure}
    \begin{subfigure}{0.45\textwidth}
        \includegraphics[width=\textwidth]{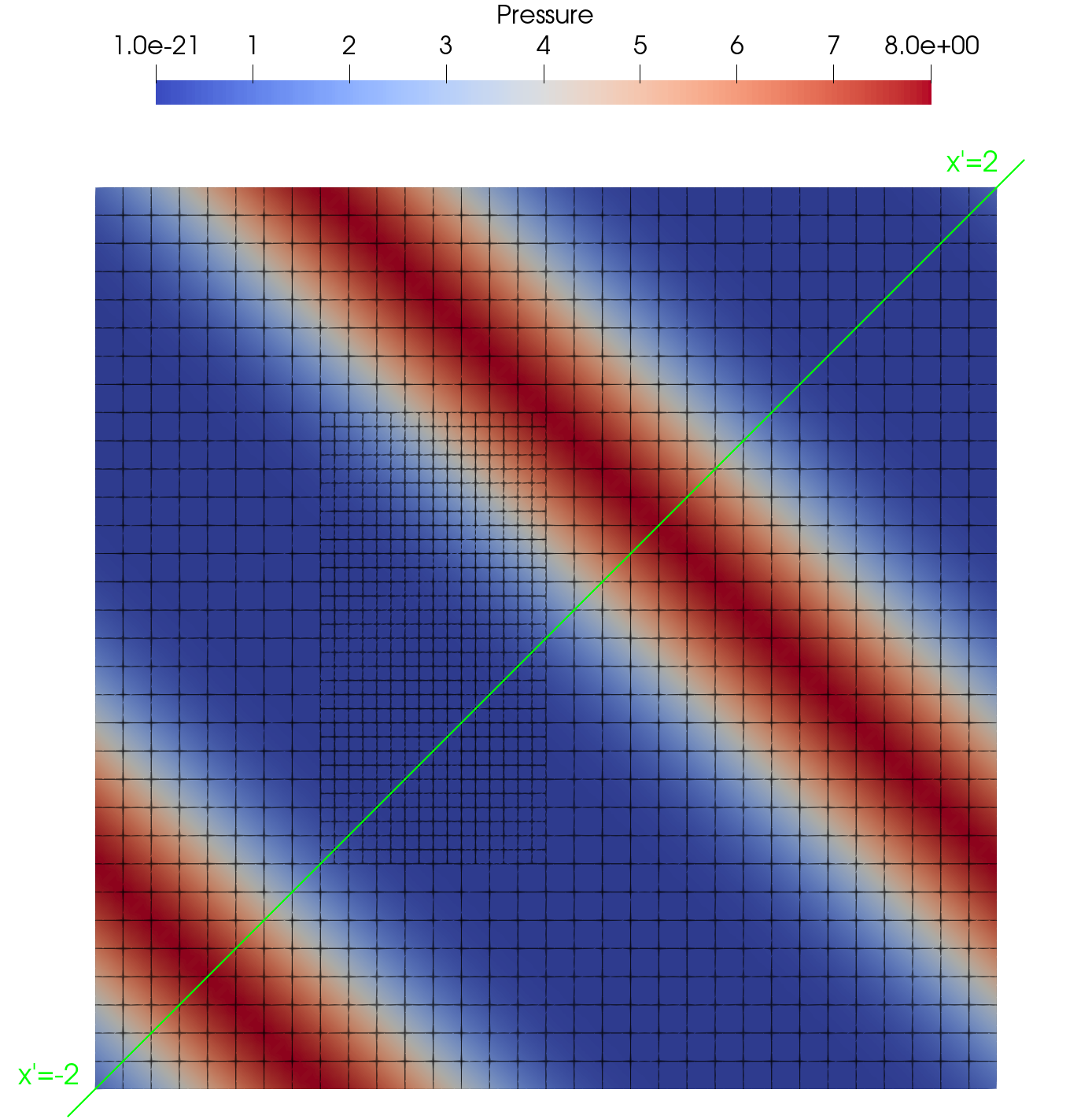}
    \end{subfigure}
    \caption{Density and pressure contours of the initial condition of the isentropic flow setup. A polynomial degree of $N=3$ and a mesh with base refinement level 5 with an additional refined box is used.}
    \label{fig:isentropicflow_initialcondition}
\end{figure}

Again, we want to use a static nonconforming mesh with a more refined region. The refined box from \Cref{fig:RefinedMesh_Box} has mortars in the low-density and low-pressure region when transferred to the new domain. We manually refine a box with $[-1/\sqrt{2}, 0] \times [-1/\sqrt{2}, 1/\sqrt{2}]$.
\Cref{fig:isentropicflow_initialcondition} illustrates the initial condition and the used mesh that contains 896 elements of refinement level 5 and 512 of level 6.

\begin{figure}[pos=htbp]
    \centering
    \begin{subfigure}{0.45\textwidth}
        \includegraphics[width=\textwidth]{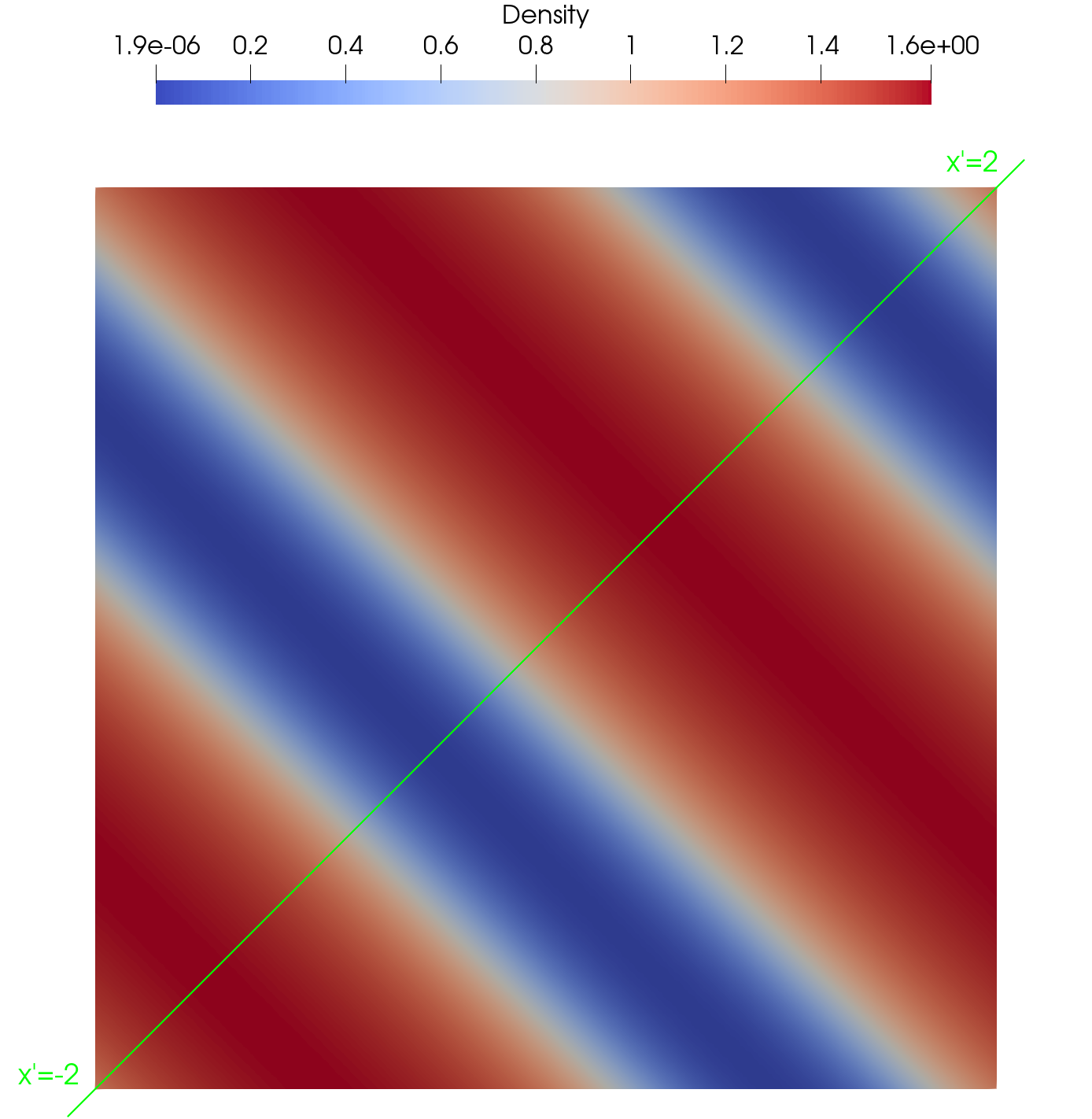}
    \end{subfigure}
    \begin{subfigure}{0.45\textwidth}
        \includegraphics[width=\textwidth]{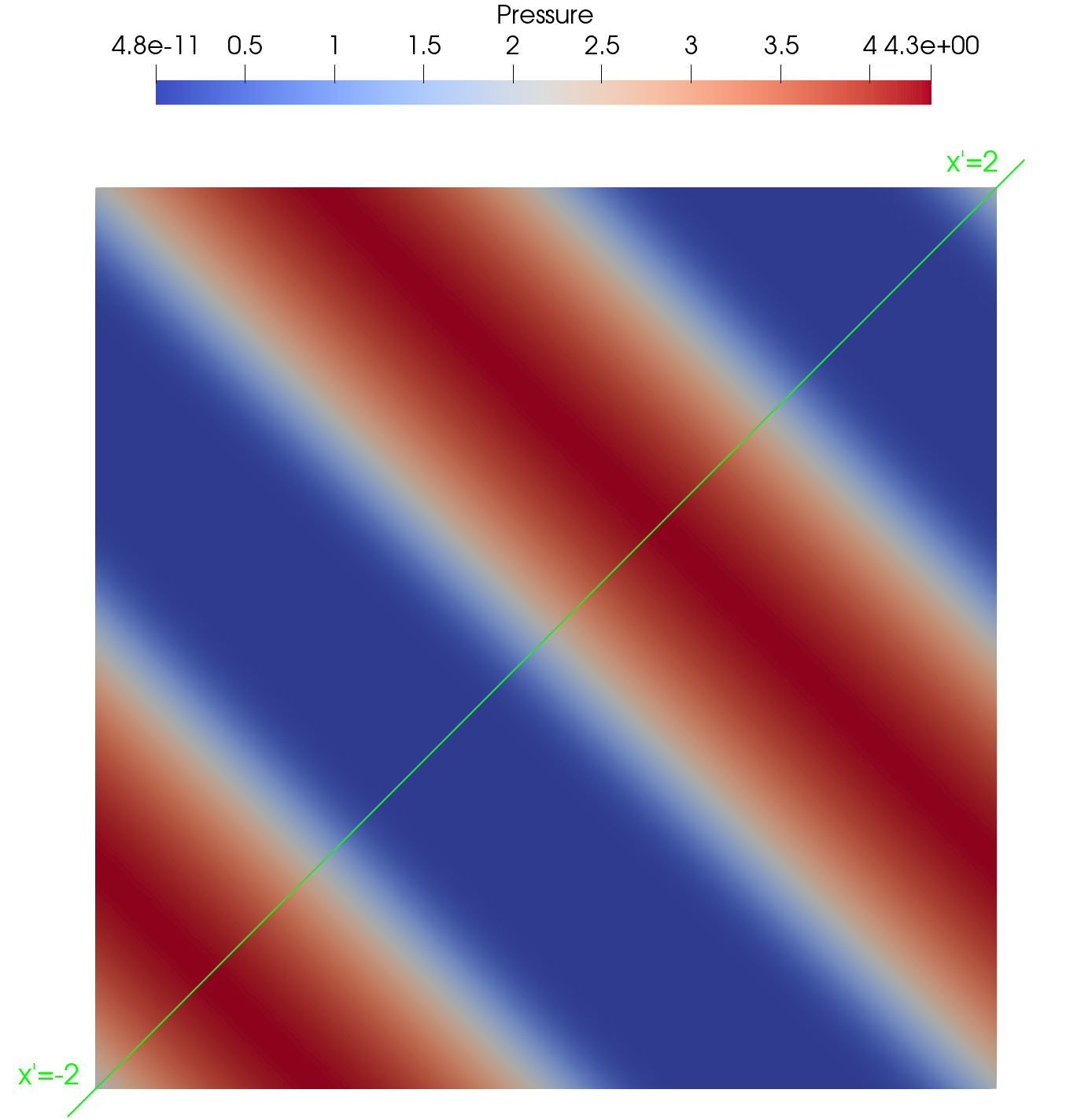}
    \end{subfigure}
    \caption{Density and pressure contours of the solution of the isentropic flow simulation at time $t=0.1$. A polynomial degree of $N=3$ and a mesh with base refinement level 5 with an additional refined box is used.}
    \label{fig:isentropicflow}
\end{figure}

The simulation runs until $t=0.1$ and uses positivity limiting of density and pressure in the volume integral and the mortars. Due to the pure positivity limiting we use the low-order stability CFL condition~\eqref{eq:CFLcondition_lowordermethod} with $\text{CFL}=0.5$ to reduce the number of time steps. \Cref{fig:isentropicflow} shows the resulting density and pressure contours at $t=0.1$.

\begin{figure}
    \centering
    \begin{subfigure}{0.45\linewidth}
        \includegraphics[width=\linewidth]{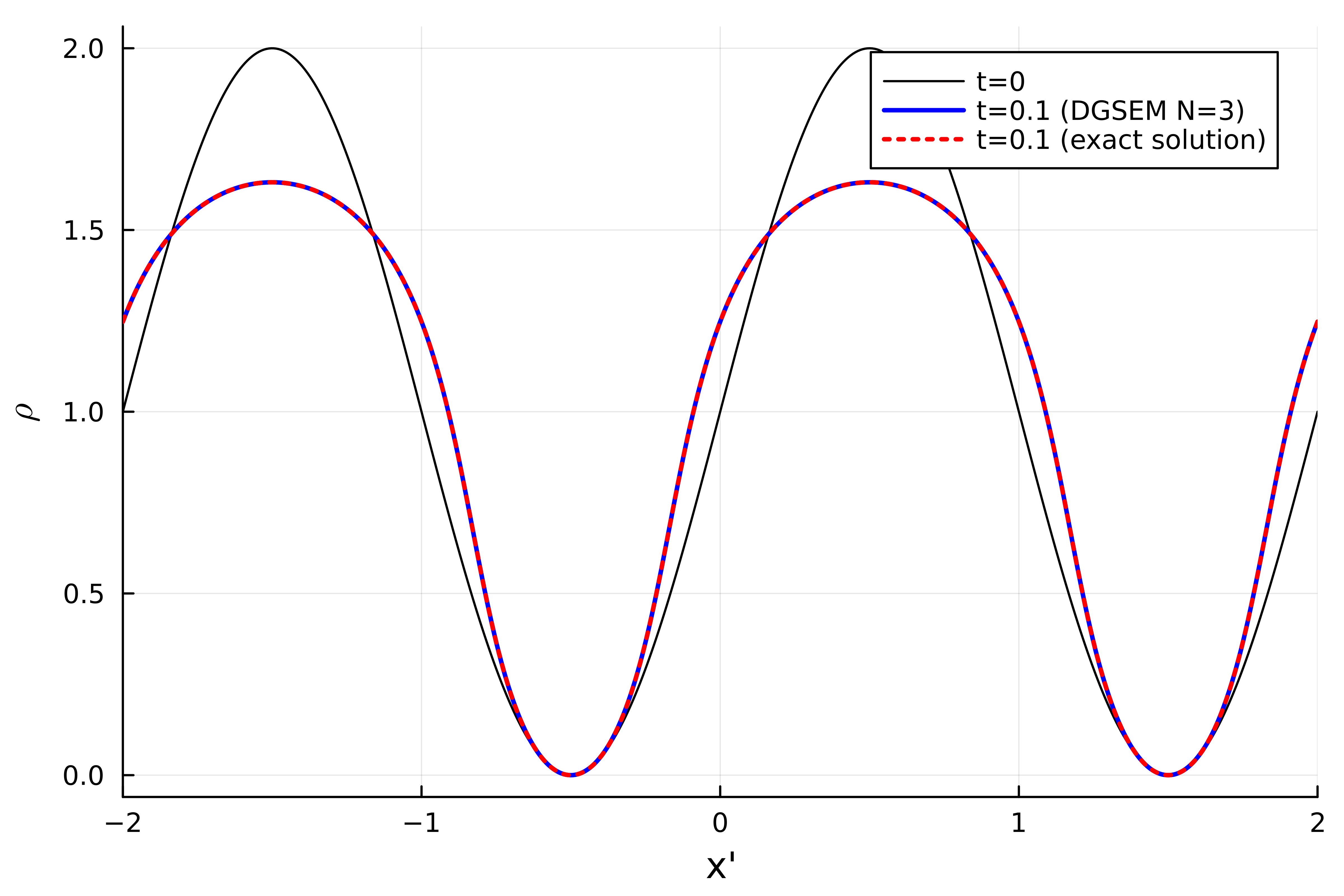}
        \caption{Density}
    \end{subfigure}
    \begin{subfigure}{0.45\linewidth}
        \includegraphics[width=\linewidth]{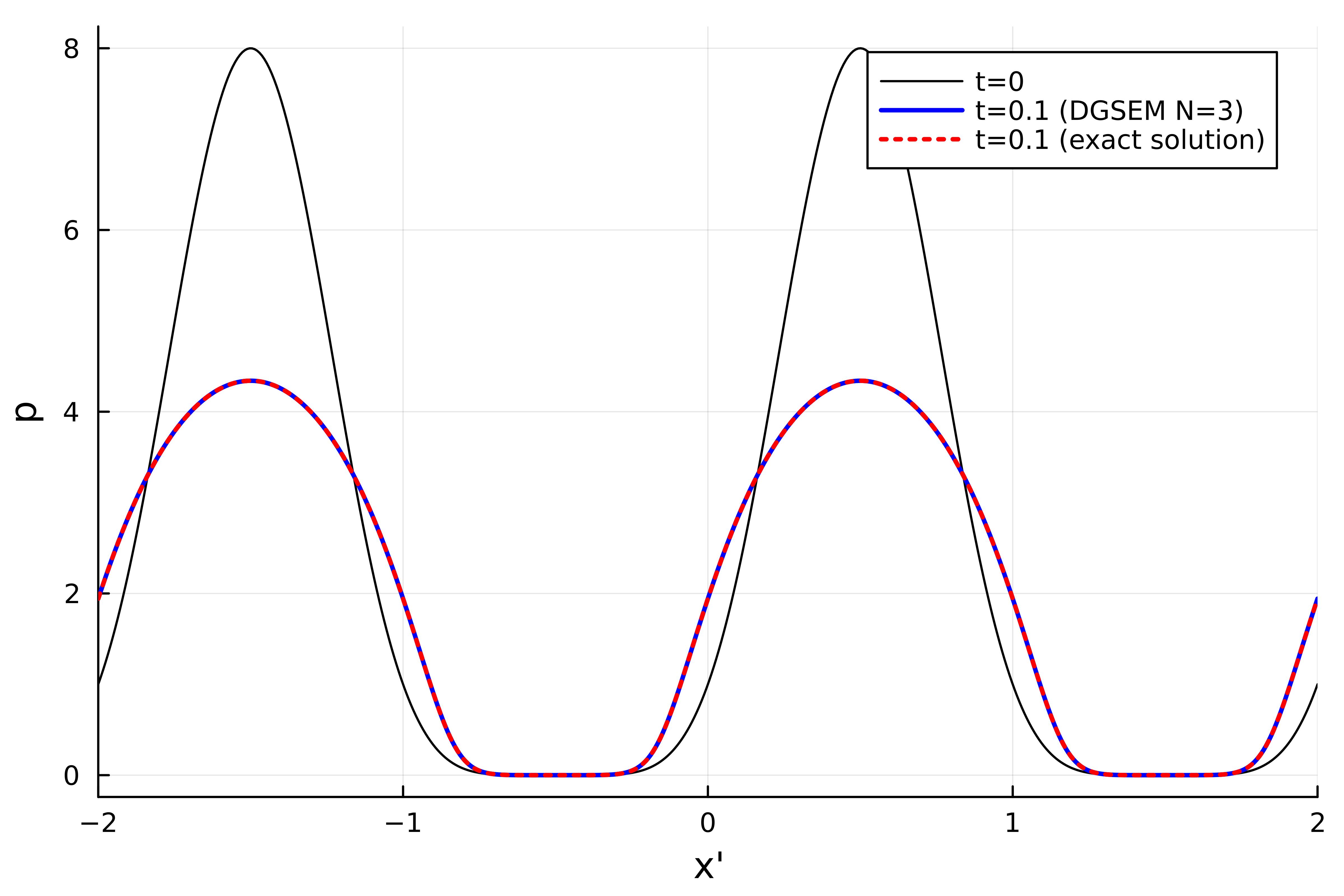}
        \caption{Pressure}
    \end{subfigure}
    \caption{Resulting density and pressure along the diagonal $x'$ slice for the simulation of the isentropic flow at time $t=0.1$.}
    \label{fig:isentropicflow_along_diagonal}
\end{figure}

To test the result qualitatively, we look at the solution on the diagonal $x'$. We compare the resulting values with the exact solution in \Cref{fig:isentropicflow_along_diagonal} and confirm the correctness.

\begin{figure}
    \centering
    \begin{subfigure}[t]{0.49\textwidth}
        \includegraphics[clip, trim=0 0 0 180, width=\textwidth]{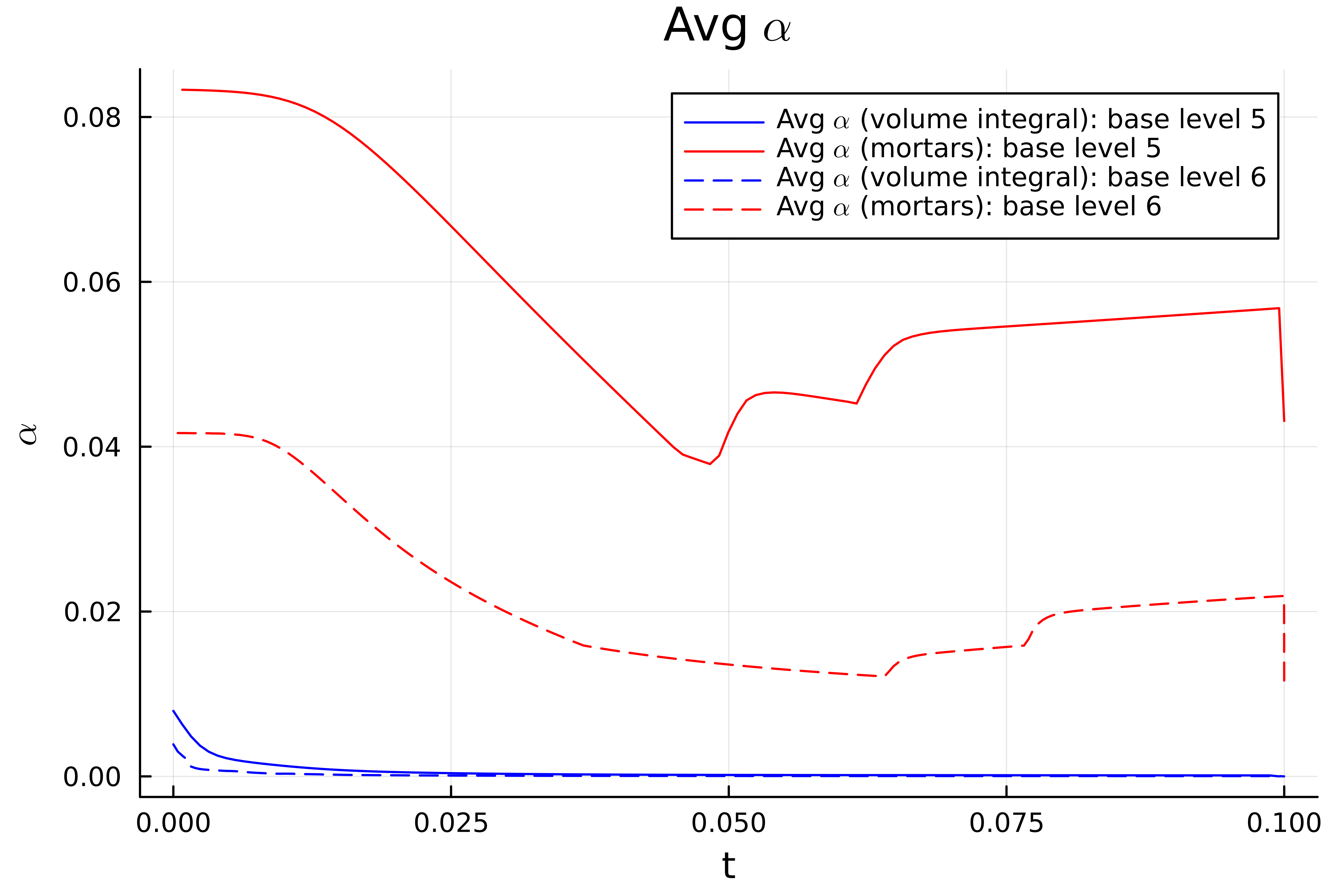}
        \caption{Evolution of limiting factors}
        \label{fig:isentropicflow_limitingfactor_evolution}
    \end{subfigure}
    \begin{subfigure}[t]{0.45\textwidth}
        \includegraphics[trim=850 180 550 180, clip, width=\textwidth]{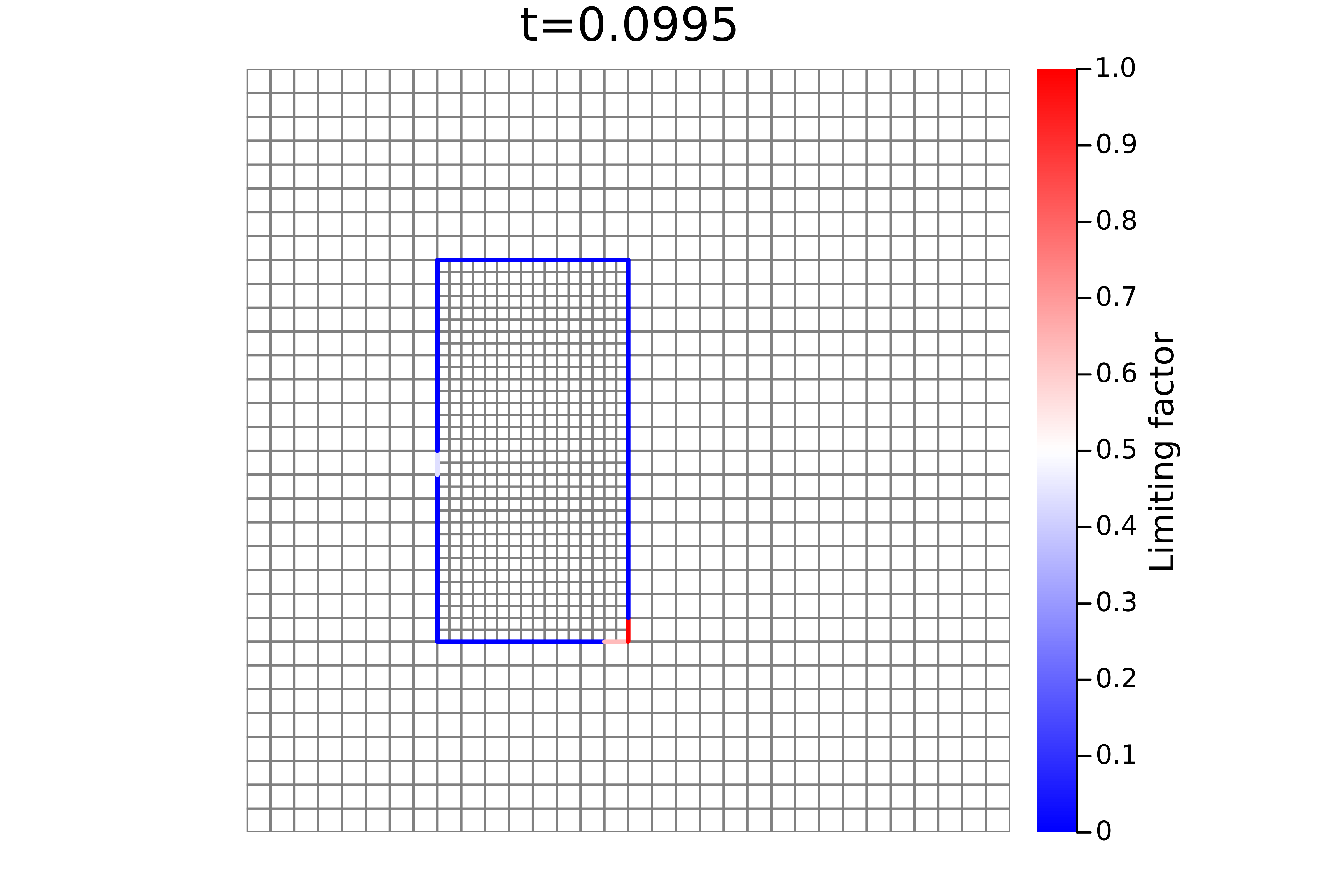}
        \caption{Illustration of the limiting factors at mortars at $t=0.1$.}
    \label{fig:isentropicflow_alpha_mortar}
    \end{subfigure}
    \caption{Analysis of the mortar limiting in the simulation of the isentropic flow.}
\end{figure}

Furthermore, we want to analyze the amount of limiting needed in this simulation. For that we take a look at the evolution of the averaged limiting factor $\alpha$ in the volume integral and at the mortars. We do this once for the given resolution and once with a mesh refined by one additional level ($3584$ elements of refinement level 6 and $2048$ elements of refinement level 7). The results are shown in \Cref{fig:isentropicflow_limitingfactor_evolution}.

We note that the fact that the amount of limiting at the mortars is higher in comparison to the volume integral is expected due to multiple reasons. The first one is that in the volume integral we limit in a more local way on the node level while we have one limiting factor per mortar. Therefore, the limiting factor has to ensure the bounds of all adjacent nodes on the interface.
Secondly, the computation of the high-order mortar flux at the mortars contains an interpolation and a projection step between the different resolutions. This can cause non-admissible numerical states which would break the simulation and therefore yields the full use of the stable low-order mortar flux at the affected mortars. The high-order flux in the volume integral is less affected by this.
Another reason specifically in our implementation is the fact that we first apply the correction in the volume integral and only afterward at the mortars. Because of that, the limiting process in the volume integral starts with the pure low-order solution, while the solution for the mortar limiting already pushes the bounds to their limits.
A fourth reason has to do with the way we compute the number. The visualized limiting factor is the weighted average of all limiting factors for both volume and mortars. For the volume integral, that includes all nodes in the domain. For the mortars, there are only a few interfaces of the domain that are mortars. This increases the average limiting factor in this visualization.

One can see that the average amount of limiting decreases for the more refined mesh. There are two simple reasons for that. First, all elements and mortars are smaller in the more refined mesh which yields limiting in a more local way. Additionally, in our FCT-type limiting approach, the limiting factor directly depends on the time step size which is typically smaller for more refined meshes~\cite{RUEDARAMIREZ2022, ruedaramirez2024}.
Note that the decreased amount of limiting in the last time step is also normal due to the last argument. The last time step of the simulations might be smaller to match the end time and therefore have less limiting.

\Cref{fig:isentropicflow_alpha_mortar} sketches the used limiting factors at the mortars at final time $t=0.1$. Note that it is the final step where the limiting is lower compared to the previous steps. The illustration confirms the local impact of the limiting in the low density and low pressure region. The fact that the limiting is not perfectly symmetric in the corner can be explained since the mesh is not created for this particular setup.

Last, we want to do a convergence test using this challenging setup. Due to the use of the stable low-order fluxes we expect the order of convergence to not reach full high-order $N+1$.
The results of the convergence test with polynomial degree $3$ and $4$ are shown in \Cref{table:ConvergenceTestIsentropicFlow}.

\begin{table}[pos=htbp]
    \begin{subtable}{\textwidth}
    \centering
    \resizebox{0.85\textwidth}{!}{%
    \begin{tabular}{c|cc|cc|cc|cc}
    \toprule
     & \multicolumn{2}{c|}{\(\rho\)}
     & \multicolumn{2}{c|}{\(\rho v_1\)}
     & \multicolumn{2}{c|}{\(\rho v_2\)}
     & \multicolumn{2}{c}{\(\rho e\)} \\
    level & error & EOC & error & EOC & error & EOC & error & EOC \\
    \midrule
    2(+1) & $3.45 \times 10^{-2}$ & -    & $4.09 \times 10^{-2}$ & -    & $4.06 \times 10^{-2}$ & -    & $6.94 \times 10^{-2}$ & -    \\
    3(+1) & $4.31 \times 10^{-3}$ & 3.00 & $4.82 \times 10^{-3}$ & 3.08 & $4.80 \times 10^{-3}$ & 3.08 & $8.13 \times 10^{-3}$ & 3.09 \\
    4(+1) & $5.11 \times 10^{-4}$ & 3.08 & $5.30 \times 10^{-4}$ & 3.19 & $5.27 \times 10^{-4}$ & 3.19 & $7.05 \times 10^{-4}$ & 3.53 \\
    5(+1) & $6.55 \times 10^{-5}$ & 2.96 & $4.53 \times 10^{-5}$ & 3.55 & $4.52 \times 10^{-5}$ & 3.55 & $5.11 \times 10^{-5}$ & 3.79 \\
    6(+1) & $7.81 \times 10^{-6}$ & 3.07 & $3.32 \times 10^{-6}$ & 3.77 & $3.32 \times 10^{-6}$ & 3.77 & $3.44 \times 10^{-6}$ & 3.89 \\
    7(+1) & $8.34 \times 10^{-7}$ & 3.23 & $2.34 \times 10^{-7}$ & 3.83 & $2.34 \times 10^{-7}$ & 3.83 & $2.26 \times 10^{-7}$ & 3.93 \\
    \bottomrule
    \end{tabular}}
    \subcaption{Polynomial degree of 3.}
    \end{subtable}

    \begin{subtable}{\textwidth}
    \centering
    \resizebox{0.85\textwidth}{!}{%
    \begin{tabular}{c|cc|cc|cc|cc}
    \toprule
     & \multicolumn{2}{c|}{\(\rho\)}
     & \multicolumn{2}{c|}{\(\rho v_1\)}
     & \multicolumn{2}{c|}{\(\rho v_2\)}
     & \multicolumn{2}{c}{\(\rho e\)} \\
    level & error & EOC & error & EOC & error & EOC & error & EOC \\
    \midrule
    2(+1) & $1.77 \times 10^{-2}$ & -    & $1.76 \times 10^{-2}$ & -    & $1.70 \times 10^{-2}$ & -    & $2.72 \times 10^{-2}$ & -    \\
    3(+1) & $1.16 \times 10^{-3}$ & 3.93 & $1.31 \times 10^{-3}$ & 3.74 & $1.30 \times 10^{-3}$ & 3.71 & $2.12 \times 10^{-3}$ & 3.68 \\
    4(+1) & $9.76 \times 10^{-5}$ & 3.58 & $7.84 \times 10^{-5}$ & 4.06 & $7.83 \times 10^{-5}$ & 4.06 & $1.05 \times 10^{-4}$ & 4.33 \\
    5(+1) & $7.98 \times 10^{-6}$ & 3.61 & $3.61 \times 10^{-6}$ & 4.44 & $3.61 \times 10^{-6}$ & 4.44 & $3.85 \times 10^{-6}$ & 4.77 \\
    6(+1) & $6.93 \times 10^{-7}$ & 3.52 & $1.45 \times 10^{-7}$ & 4.64 & $1.45 \times 10^{-7}$ & 4.64 & $1.29 \times 10^{-7}$ & 4.90 \\
    7(+1) & $4.91 \times 10^{-8}$ & 3.82 & $5.68 \times 10^{-9}$ & 4.67 & $5.68 \times 10^{-9}$ & 4.67 & $5.61 \times 10^{-9}$ & 4.52 \\
    \bottomrule
    \end{tabular}}
    \subcaption{Polynomial degree of 4.}
    \end{subtable}
    \caption{$L_2$ errors and EOC for convergence test with isentropic flow setup using positivity limiting. Each mesh consists of elements of the given refinement level with an additionally refined box.}
    \label{table:ConvergenceTestIsentropicFlow}
\end{table}

The convergence tables show a decreased order of convergence as we expected. We reach convergence orders between 3 and 4 for a polynomial degree of 3 and between 3.5 and 5 for polynomial degree of 4. For such a challenging setup, the results are as expected and confirm the functionality of the method.

\subsection{Kelvin-Helmholtz Instability}
Next, we want to consider the inviscid and two-dimensional Kelvin-Helmholtz instability (KHI) setup, which was presented for instance in \cite{RR2021,RUEDARAMIREZ2022}.

The initial condition is given by
\begin{align}\label{eq:KHI_IniCond}
\rho (x,y) &= \frac{1}{2} + \frac{3}{4} B,
& p (x,y) &= 1, \nonumber\\
v_{1} (x,y) &= \frac{1}{2} \left( B-1 \right),
& v_{2} (x,y) &= \frac{1}{10} \sin(2 \pi x),
\end{align}
with $B=\tanh \left( 15 y + 7.5 \right) - \tanh(15y-7.5)$ on the periodic domain $\Omega = [-1,1]^2$.

We use a polynomial degree of $N=7$. From now on, we want to use adaptive mesh refinement instead of manually refined fixed meshes. For that, we use an AMR indicator that measures whether an element is going to be refined or coarsened based on different parameters~\cite{ranocha2022adaptive}.
Here, we an AMR indicator adapted from a FEM indicator by Löhner~\cite{LOHNER1987323}. It was also used in the FLASH astrophysical code~\cite{fryxell2000flash}. We apply it for the density.
We tessellate the domain with AMR into elements of refinement levels between 4 and 8 which corresponds to a refinement of between $2^4=16$ and $2^8=256$ elements per spatial dimension. Note that the amount of refinement and therefore the amount of limiting in the volume integral and the mortars very much depends on the AMR setup.

The simulation uses positivity limiting in the volume integral and in the mortars to run the simulation until $t=3.7$.
Therefore, we again use the low-order stability CFL condition~\eqref{eq:CFLcondition_lowordermethod} with a CFL number of $0.2$.

\begin{figure}[pos=htbp]
    \centering
    \includegraphics[trim=1800 440 1800 160, clip, width=0.48\textwidth]{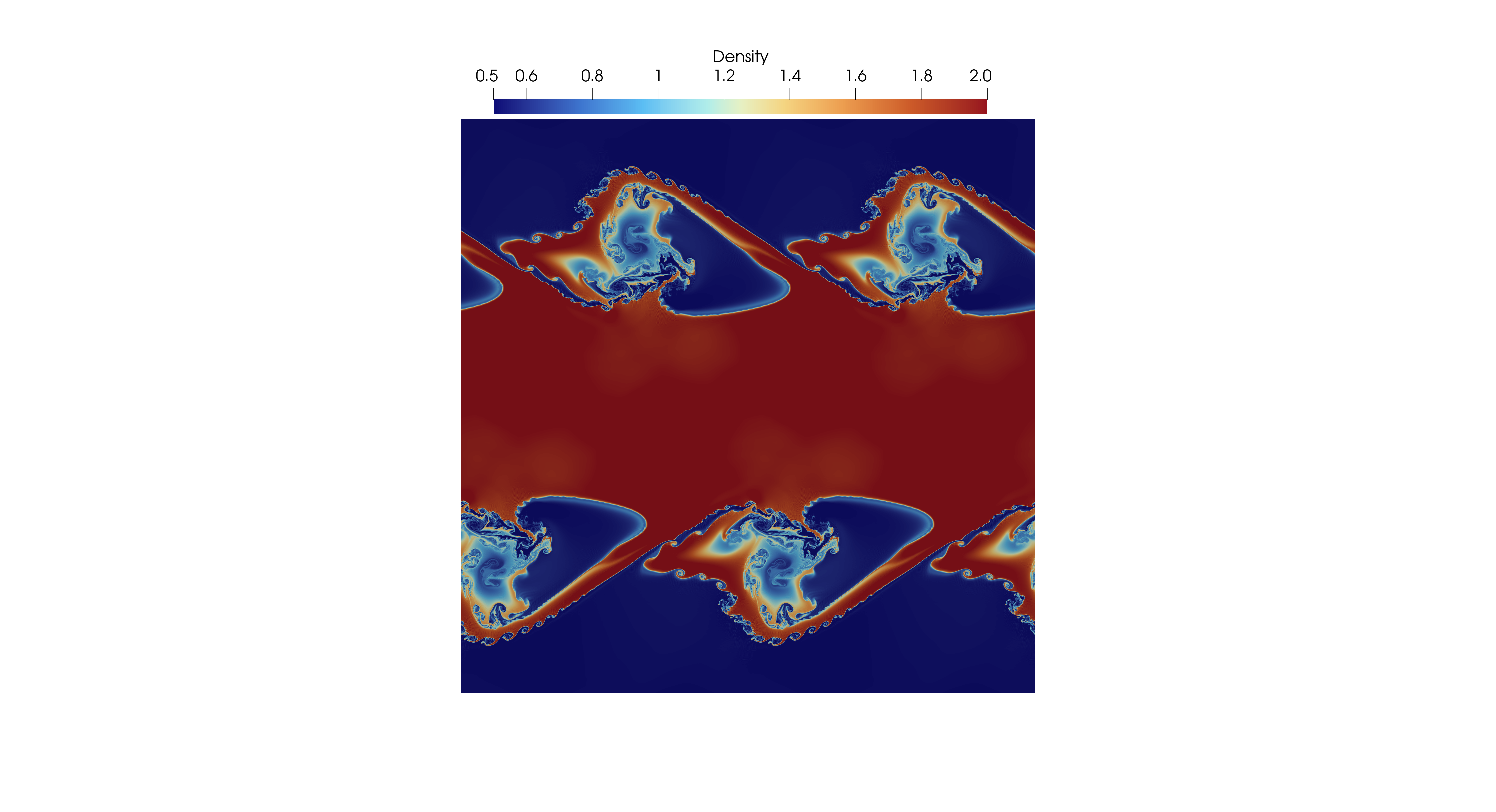}
    \includegraphics[trim=1800 440 1800 160, clip, width=0.48\textwidth]{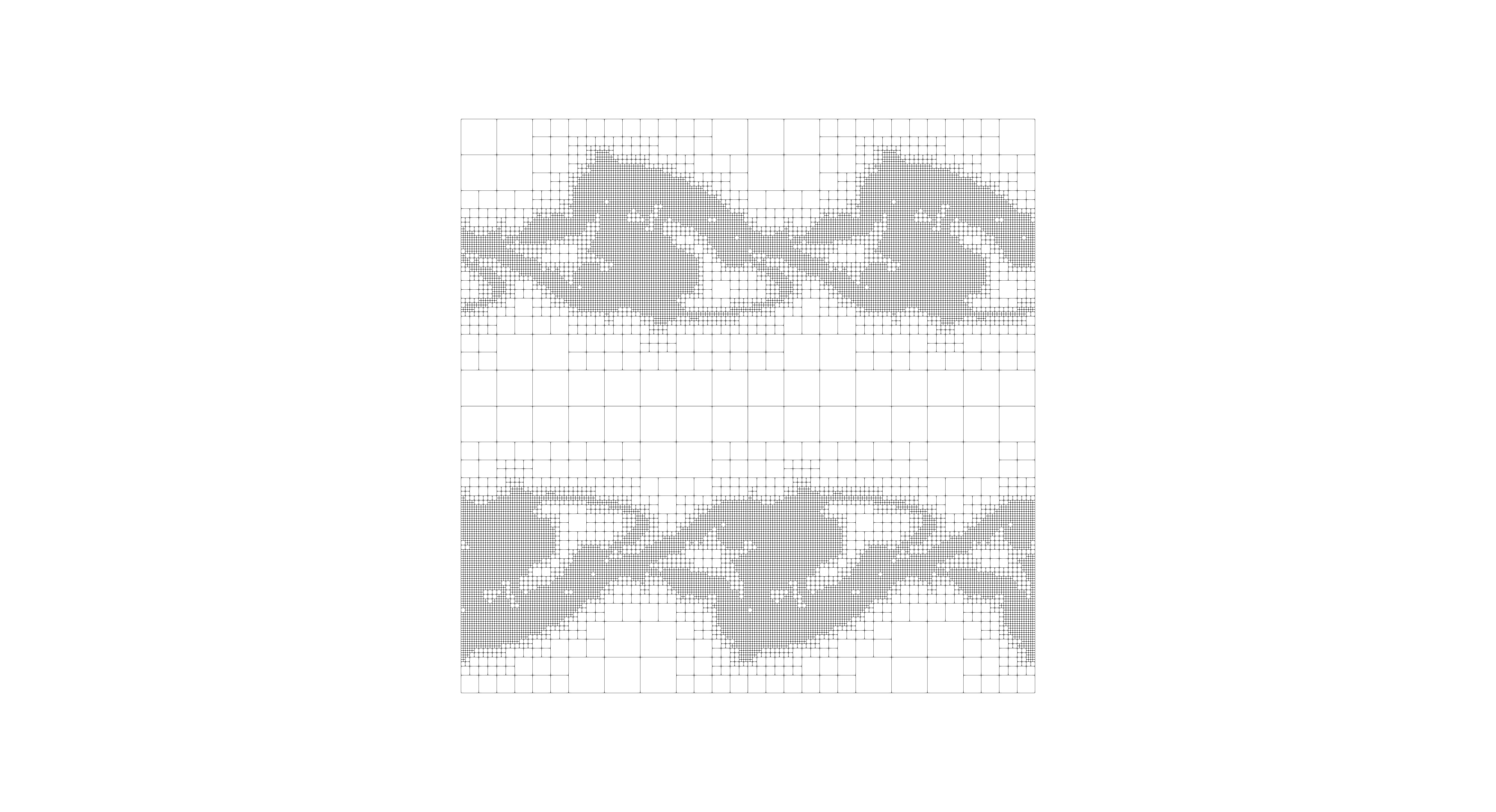}
    \caption{Density contours and the used mesh for the  simulation of the Kelvin-Helmholtz instability with polynomial degree of $7$.}
    \label{fig:khi}
\end{figure}

The resulting density contour and the mesh at $t=3.7$ are shown in \Cref{fig:khi}. The final mesh contains 21508 elements corresponding to 1376512 degrees of freedom. 18640 of those elements are of the highest refinement level.
The general form of the solution fits the results of this setup in the past \cite{rueda2022flux, RUEDARAMIREZ2022}. It shows the many small-scale structures typical for this setup with the minimum amount of added dissipation due to pure positivity limiting. They seem to be generated by either the higher resolution in the finest elements or in general by the use of AMR with mortars.

\begin{figure}
    \centering
    \includegraphics[trim=0 0 0 200, clip, width=0.6\textwidth]{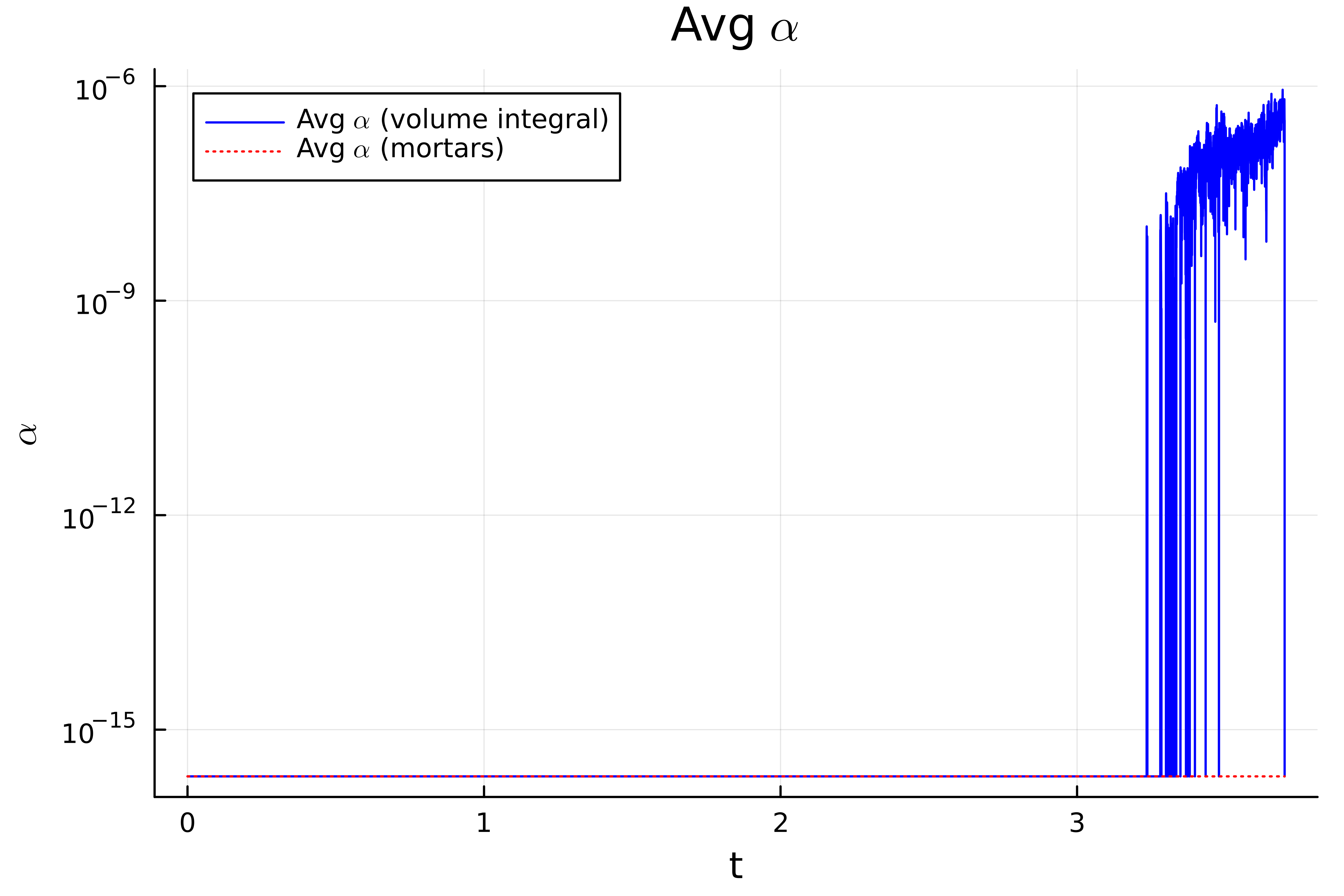}
    \caption{Evolution of the averaged limiting factors for the simulation of the Kelvin-Helmholtz instability.}
    \label{fig:khi_alphas}
\end{figure}

Consistent with previous observations, the setup of the Kelvin-Helmholtz instability only requires a very small amount of positivity limiting \cite{RUEDARAMIREZ2022}. Limiting is only applied after about $t=3.2$ in the volume integral. With this AMR setup, there is now limiting applied at the mortars (see \Cref{fig:khi_alphas}).

\subsection{Sedov Blast Explosion}\label{sec:sedovblast}
The next setup is the very challenging Sedov blast problem with a strong circular shock. We use the initial condition from the FLASH astrophysical code \cite{fryxell2000flash}. Initially, the domain is filled with a gas at rest with $v_1(t=0) = v_2(t=0) =0$ and a constant density $\rho(t=0)=1$. The background atmospheric pressure is $p_0=10^{-5}$. We insert a quantity of dimensionless energy $e=1$ into a small region at the center of the domain $\Omega=[-2,2]^2$ with radius $r_0=0.21875$,
\begin{equation}
    p(t=0) = \begin{cases}
        \frac{(\gamma - 1) e}{\pi r_0^2} & \text{if } r\leq r_0,\\
        p_0 & \text{otherwise},
    \end{cases}
\end{equation}
with $r=\sqrt{x^2+y^2}$.

We use the shock-capturing indicator developed by Hennemann et al.~\cite{hennemann2021} applied as an AMR indicator for density and pressure to add refinement at shocks. The AMR setup employs between $2^4=16$ and $2^7=128$ elements per dimension. The polynomial degree is $N=3$ and the final time is $t=3$.

When only positivity limiting is applied for density and pressure, the simulation remains sufficiently stable and does not crash. However, the resulting solution has numerical artifacts in directions of the four coordinate axes.
Therefore, we also apply local limiting including minimum and maximum principles for density and maximum principle for the modified specific entropy. We also apply the positivity limiter of Zhang and Shu \cite{zhang2011} with a global threshold of $10^{-10}$ to ensure positivity for density and pressure during the solution transfer in the refinement and coarsening steps as described in \Cref{sec:LimitingTransfer}.

\begin{figure}[pos=htbp]
    \centering
    \begin{subfigure}[t]{0.48\textwidth}
        \includegraphics[trim=900 200 900 80, clip, width=\textwidth]{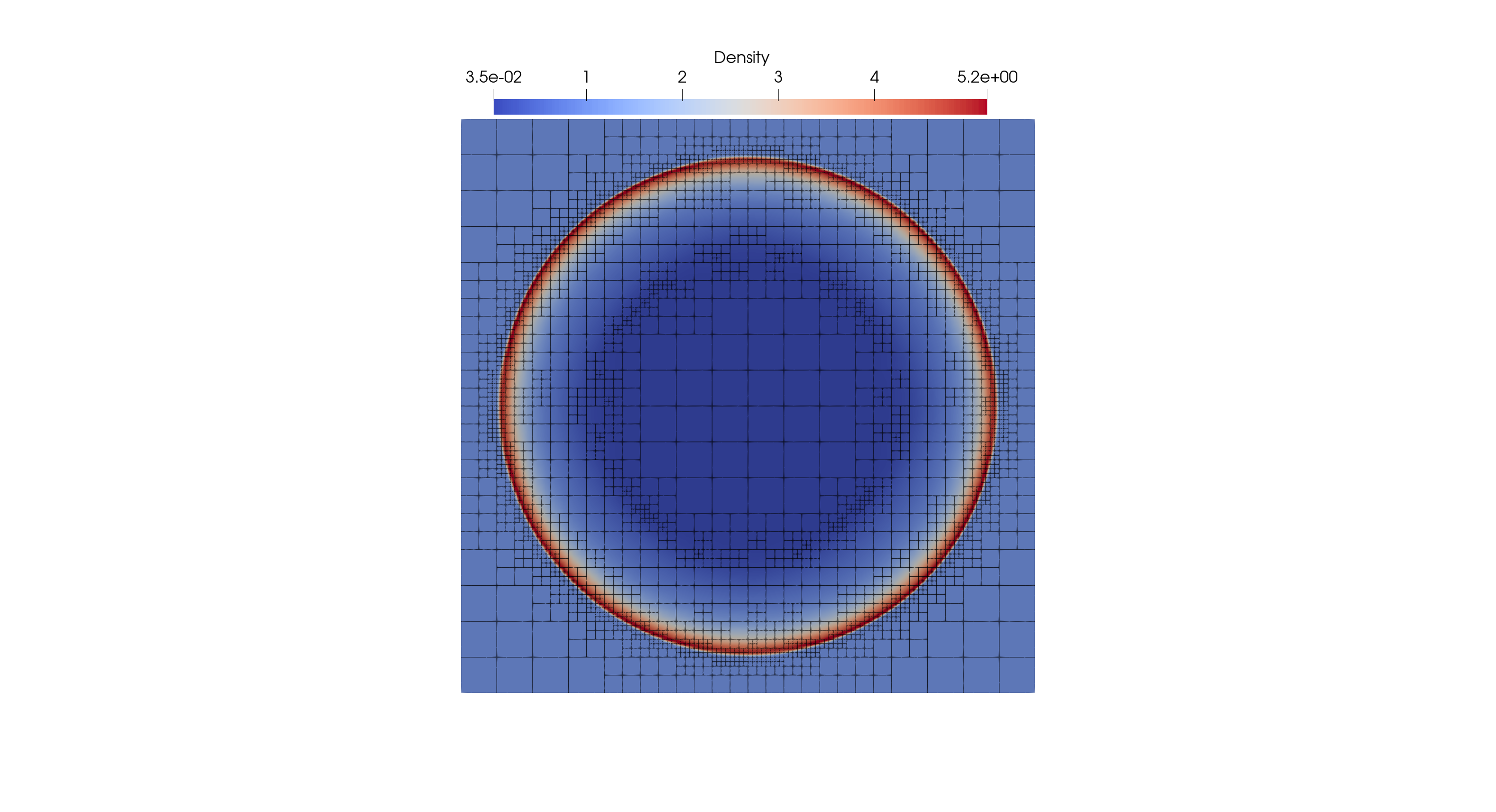}
        \subcaption{Density contours and mesh at $t=3$.}
    \end{subfigure}
    \begin{subfigure}[t]{0.48\textwidth}
        \includegraphics[trim=0 0 0 200, clip, width=\textwidth]{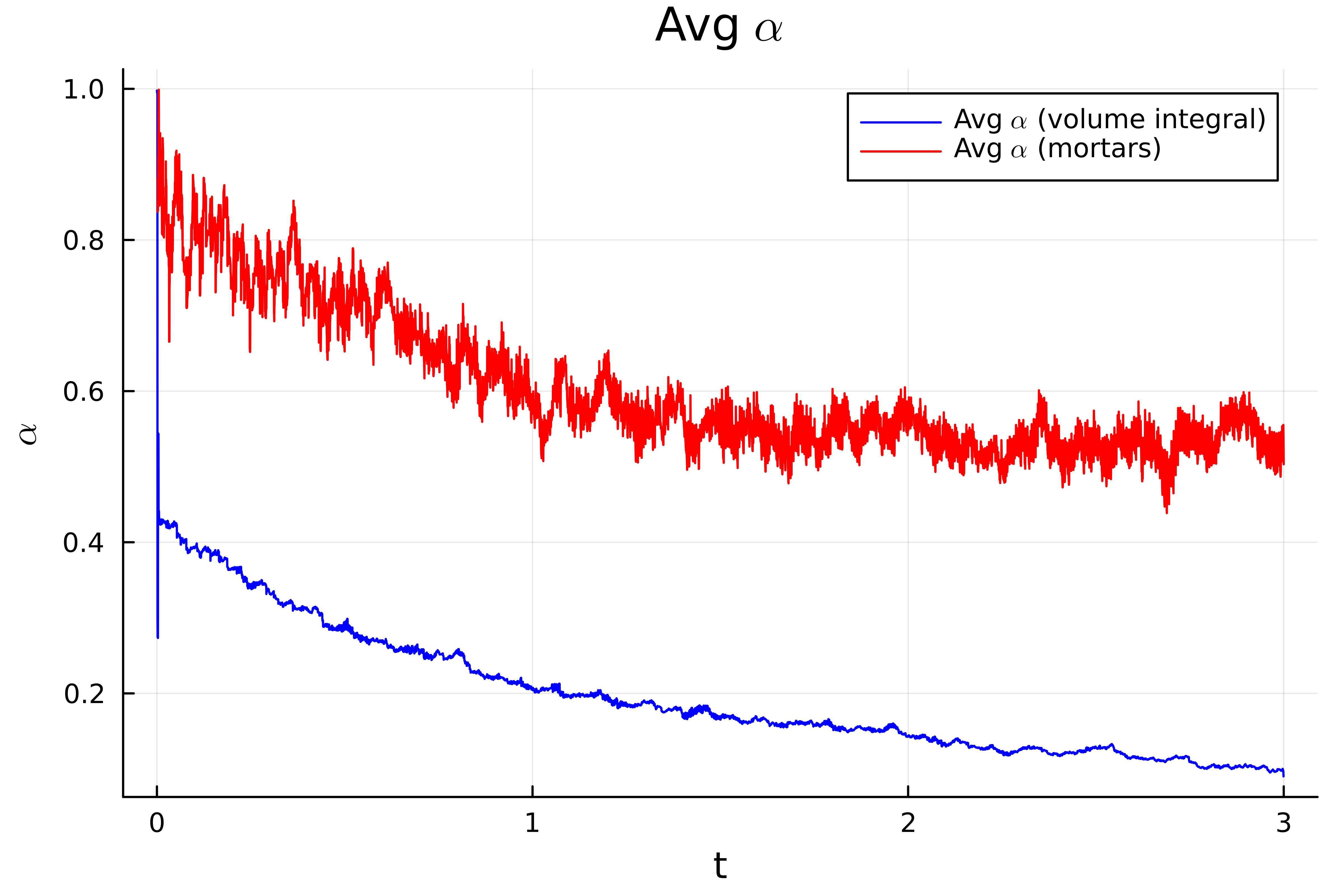}
        \subcaption{Evolution of the limiting factors.}
    \end{subfigure}
    \caption{Results of the simulation of the Sedov blast wave using the bar-state local bounds.}
    \label{fig:sedov_barstates=true}
\end{figure}

As mentioned before, in practice, we can often bypass the restrictive IDP time step restriction~\eqref{eq:timestep_restriction} and use the CFL condition for the stability of the low-order method~\eqref{eq:CFLcondition_lowordermethod}. We do a comparison here. The first simulation uses local bounds based on the bar-states (see \eqref{eq:Bounds_Barstates}) and the IDP CFL condition \eqref{eq:timestep_restriction} with $\text{CFL}=0.95$. The results are shown in \Cref{fig:sedov_barstates=true}.
The second simulation uses local bounds based on the low-order solution,
\begin{equation}\label{eq:Bounds_LowOrder}
    \min_{j\in\NN(i)} \rho_{j}^\FV \leq \rho_i \leq \max_{j\in\NN(i)} \rho_{j}^\FV, \qquad
    \min_{j\in \NN(i)} \eta(\stateL{u}_{j}^\FV) \leq \eta(\stateL{u}_i).
\end{equation}

\begin{figure}[pos=htbp]
    \centering
    \begin{subfigure}[t]{0.48\textwidth}
        \includegraphics[trim=900 200 900 80, clip, width=\textwidth]{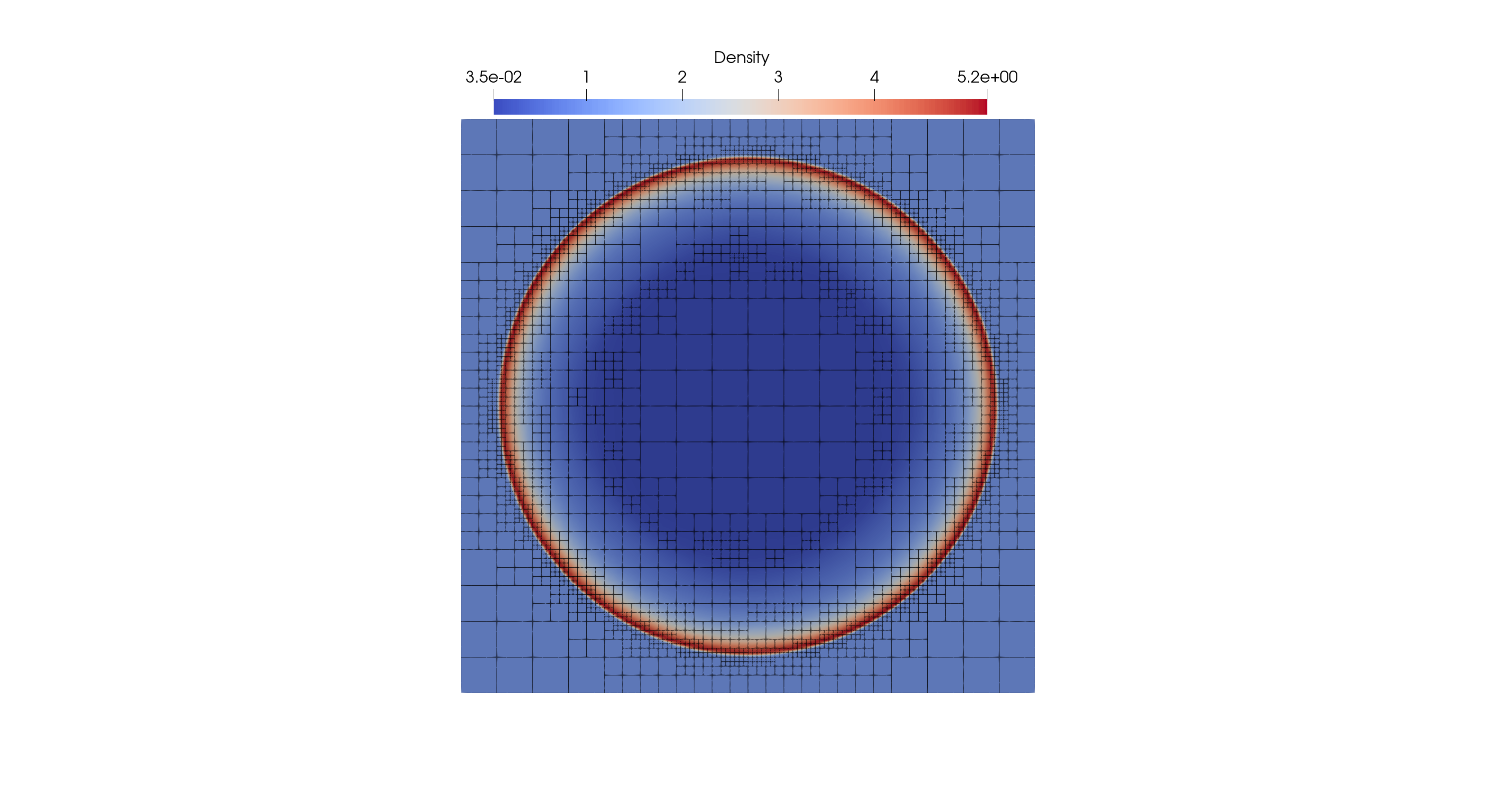}
        \subcaption{Density contours and mesh at $t=3$.}
    \end{subfigure}
    \begin{subfigure}[t]{0.48\textwidth}
        \includegraphics[trim=0 0 0 200, clip, width=\textwidth]{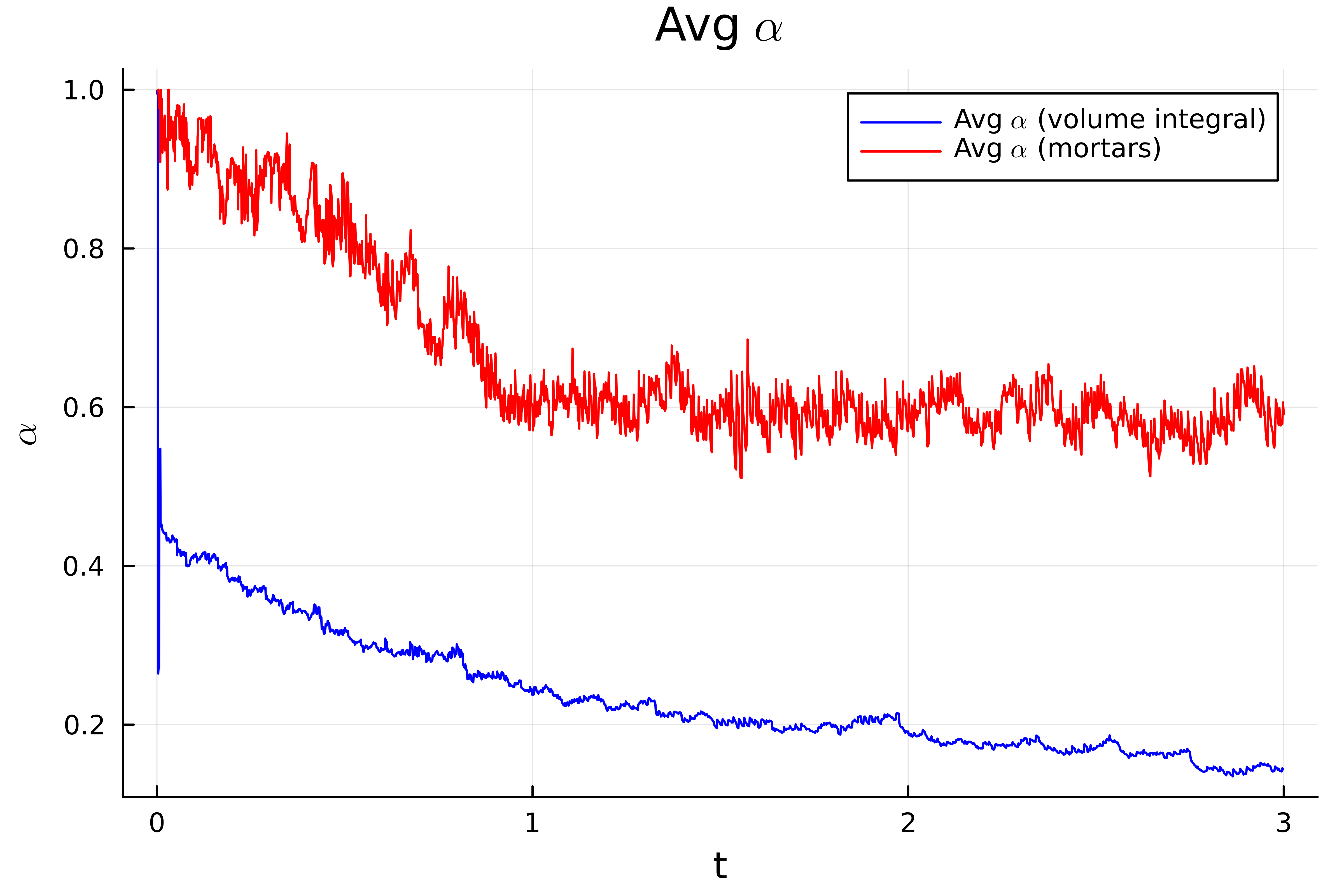}
        \subcaption{Evolution of the limiting factors.}
    \end{subfigure}
    \caption{Results of the simulation of the Sedov blast wave using local bounds based on the low-order solution.}
    \label{fig:sedov}
\end{figure}

These bounds are generally more restrictive. We use the low-order stability CFL condition~\eqref{eq:CFLcondition_lowordermethod} with $\text{CFL} = 0.4$. The results are shown in \Cref{fig:sedov}.

\begin{figure}
    \centering
    \includegraphics[clip, trim=0 0 0 200, width=0.48\textwidth]{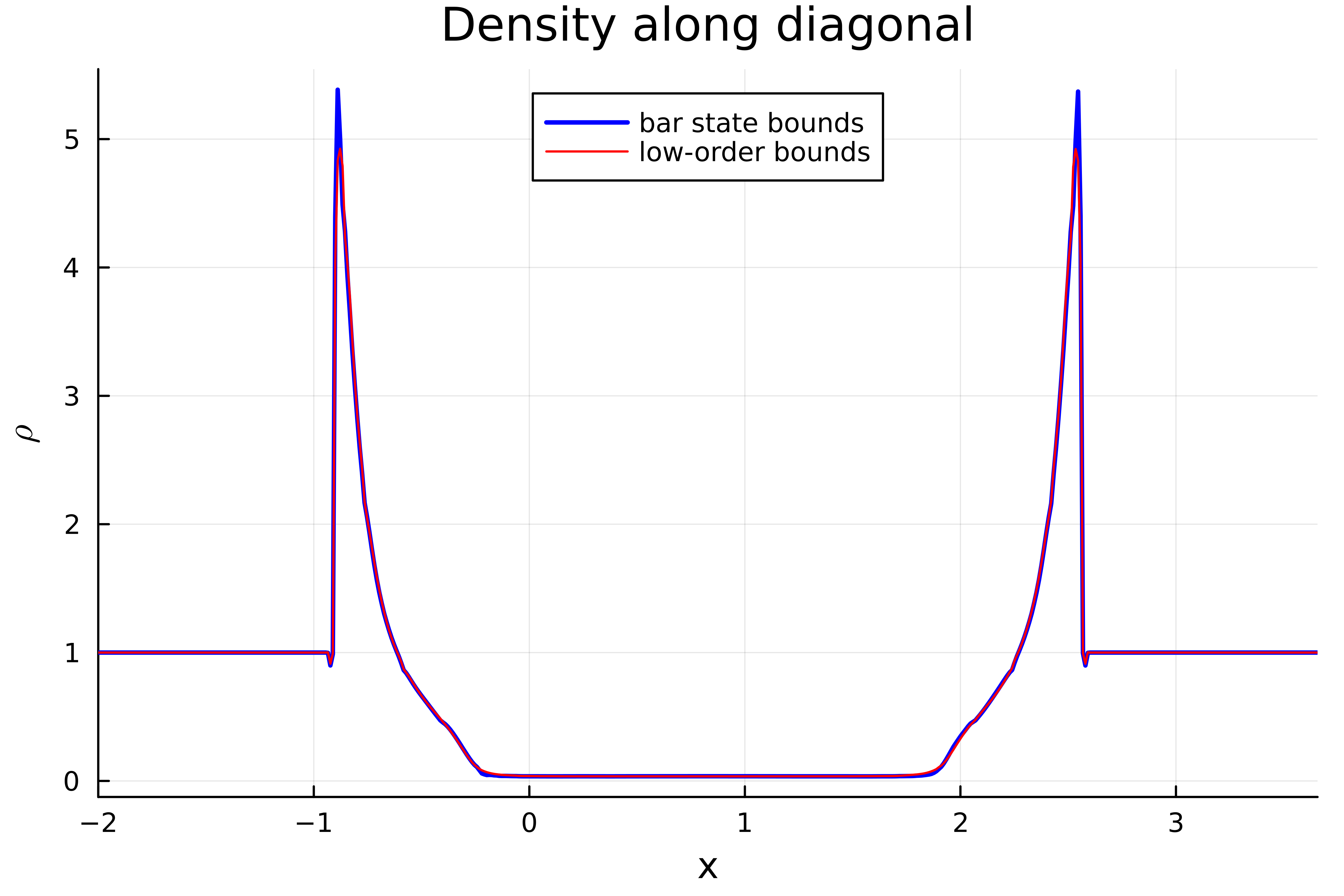}
    \includegraphics[clip, trim=0 0 0 200, width=0.48\textwidth]{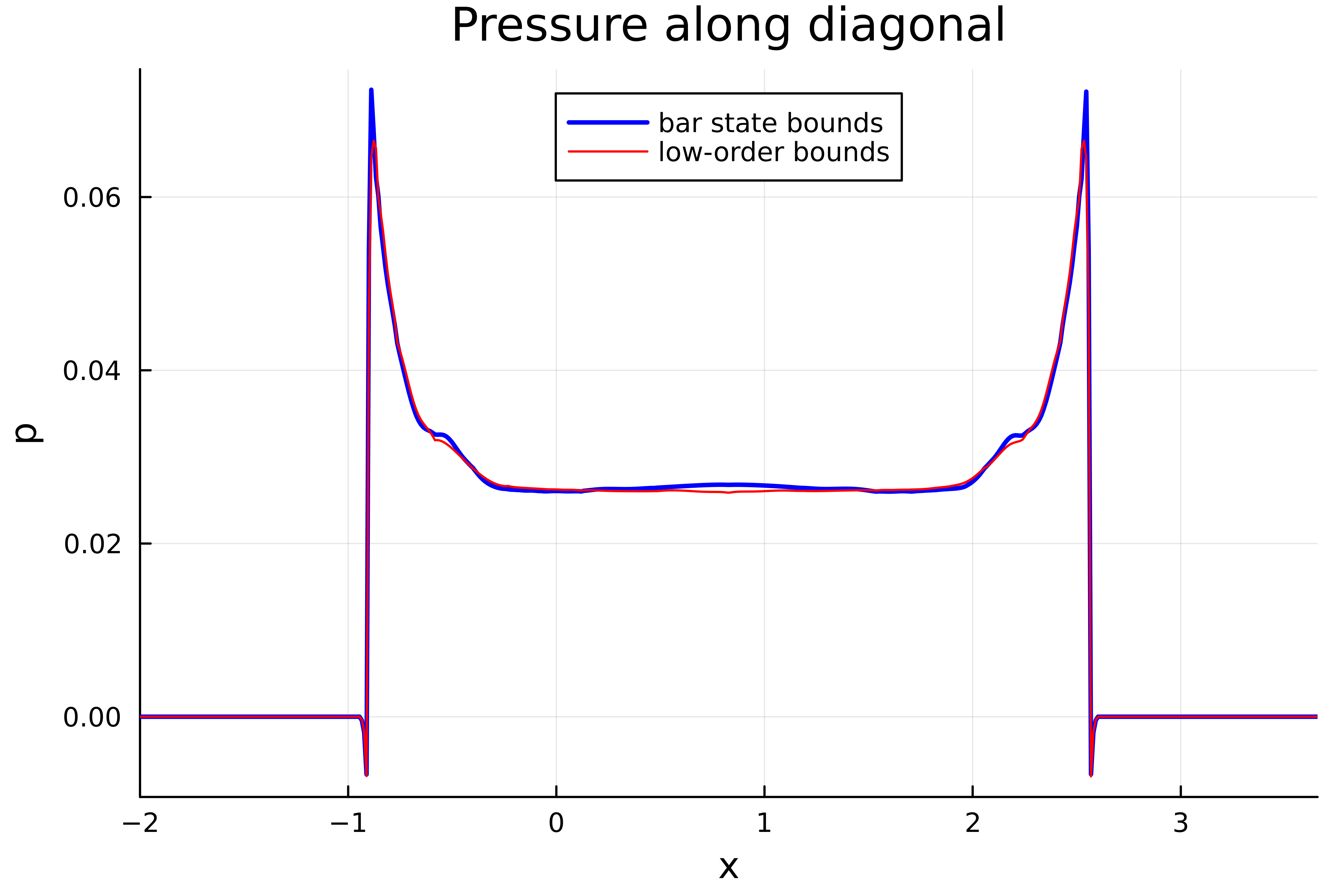}
    \caption{Resulting density and pressure along the diagonal $x=y$.}
    \label{fig:sedov_diagonal}
\end{figure}

The applied local limiting is very sensitive to small deviations of the solution with respect to the bounds. In regions with constant (within the machine's tolerance) solutions (e.g., outside the blast wave), we may have fluxes near zero $\stateL{F}^{L_2}, \tilde{\stateL{F}} \approx \vec{0}$ or a flux-difference of about zero $\stateL{F}^{L_2} - \tilde{\stateL{F}} \approx \vec{0}$. In these cases, the value of the limiting factor is essentially irrelevant. Experience shows that this value is typically rather large, rendering the average value insignificant. This phenomenon has been observed previously, for instance in~\cite{ruedaramirez2024}, and is not specific to mortar limiting.
Here, we only want to compare the given results and observe that the more restrictive low-order bounds lead to slightly larger limiting factors at the mortars while the limiting in the volume integral remains roughly the same.
The density contours are also very similar. To compare the results qualitatively, we therefore plot the solution along the diagonal $x=y$ in \Cref{fig:sedov_diagonal}. This comparison clearly shows that the simulation with low-order bounds is indeed more dissipative, with the shock being represented less sharply.

\begin{figure}
    \centering
    \includegraphics[clip, trim=0 0 0 200, width=0.48\textwidth]{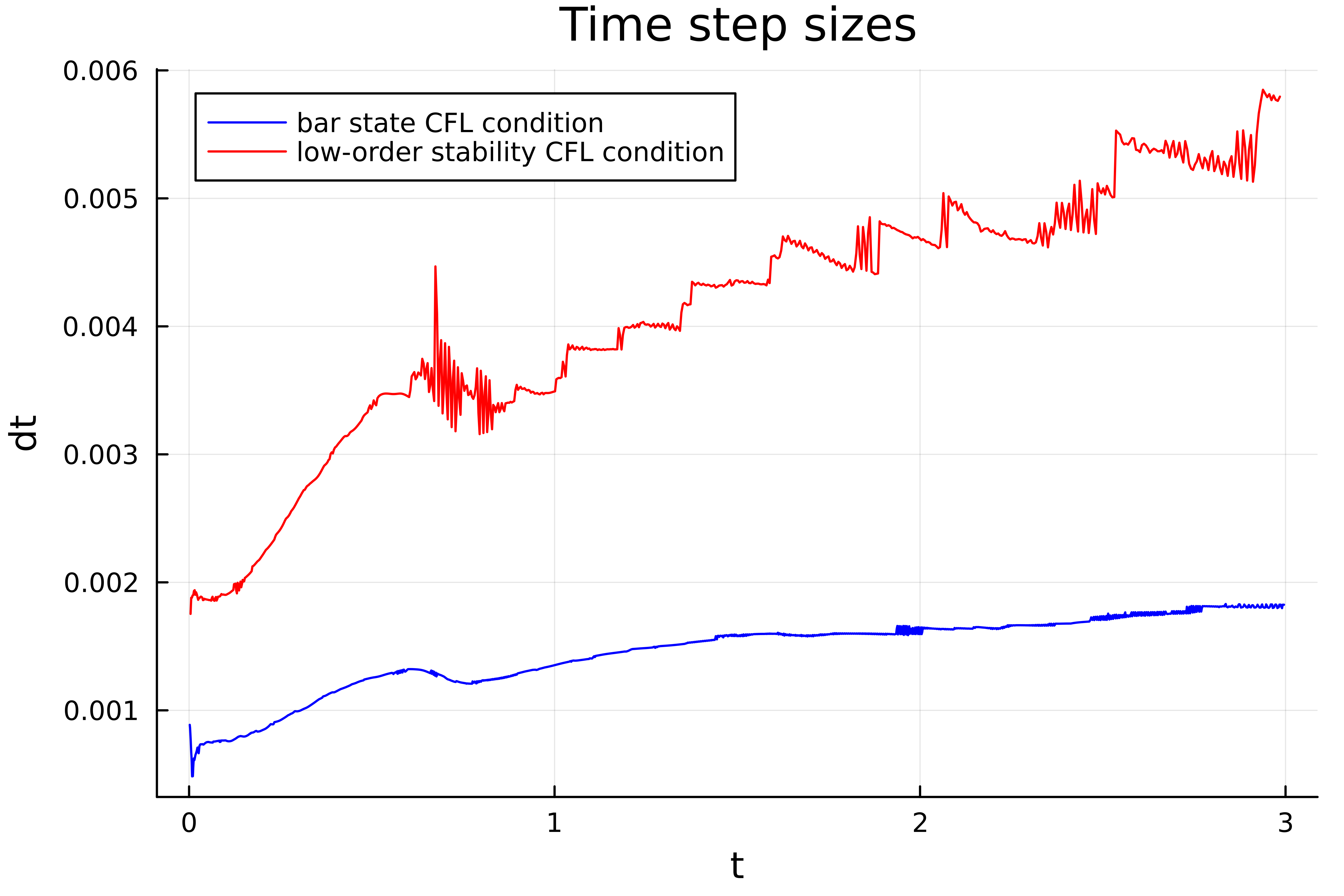}
    \caption{Comparison of time step sizes of simulations with bar-state bounds and low-order bounds.}
    \label{fig:sedov_timesteps}
\end{figure}

As a final comparison, we examine the time steps used in both simulations. \Cref{fig:sedov_timesteps} shows the corresponding time-step sizes.

Here, we observe that the time steps in the simulation with bar-state bounds, and hence with the IDP CFL condition~\eqref{eq:timestep_restriction}, are significantly smaller despite the larger CFL number. These smaller time steps result in a total of 6523 time steps for the bar-state bounds, compared to 2357 for the low-order bounds. In total, the simulation with bar-state bounds took about four times longer.

\subsection{Double Mach Reflection}
The next simulation uses the Double Mach Reflection setup based on \cite{WOODWARD1984115}. It was also used in~\cite{kuzmin2020monolithic}. The initial setup contains a propagating shock of Mach 10 with an angle of 60°. It divides the spatial domain $\Omega=[0,3.25]\times[0,1]$, which we have shortened slightly on the right side compared to the original setup, into a post-shock (left) and pre-shock (right) region.

The initial condition in the post-shock region $\Omega_L=\{(x,y)^T : \, x <1/6 + y/\sqrt{3}\}$ is defined by
\begin{equation}
    \rho_L = 8, \quad \vec{v}_L=\left(\begin{array}{c} 8.25 \cos(\pi/6) \\ -8.25 \sin(\pi/6)\end{array}\right), \quad p_L = 116.5.
\end{equation}
In the pre-shock region, $\Omega_R = \Omega\setminus \Omega_L$, we have
\begin{equation}
    \rho_R = 1.4,\quad \vec{v}_R = (0,0)^T, \quad p_R = 1.
\end{equation}
This is the first setup with a non-periodic domain. In both $\xi_1$ and positive $\xi_2$ direction we consider a characteristic-based inflow/outflow boundary condition. For that, we use the initial condition propagated over time such that the post-shock values (left) are applied for $x<1/6+(y+20t)/\sqrt{3}$ and the pre-shock values otherwise. The boundary in negative $\xi_2$ direction is divided into a characteristic-based boundary for $x<1/6$ and a reflecting slip wall boundary condition for $y\geq 1/6$. The final time is $t=0.2$.

For the spatial discretization we again use AMR. In this case, we use the AMR indicator adapted from a FEM indicator by Löhner~\cite{LOHNER1987323}. We apply it for the density. Since we have no square domain anymore, the base mesh of refinement level 0 uses a resolution of $16\times5$ elements. We use AMR with refinement levels between 0 and 6 corresponding to between $16\times 5$ and $16\cdot 2^6\times 5\cdot 2^6 = 1024\times320$ elements. The polynomial degree is $N=4$. We use local bounds based on the bar-states and a CFL of 0.9 with the IDP time step restriction~\eqref{eq:timestep_restriction}.

\begin{figure}[pos=htbp]
    \centering
    \includegraphics*[trim=200 760 200 360, clip, width=0.8\textwidth]{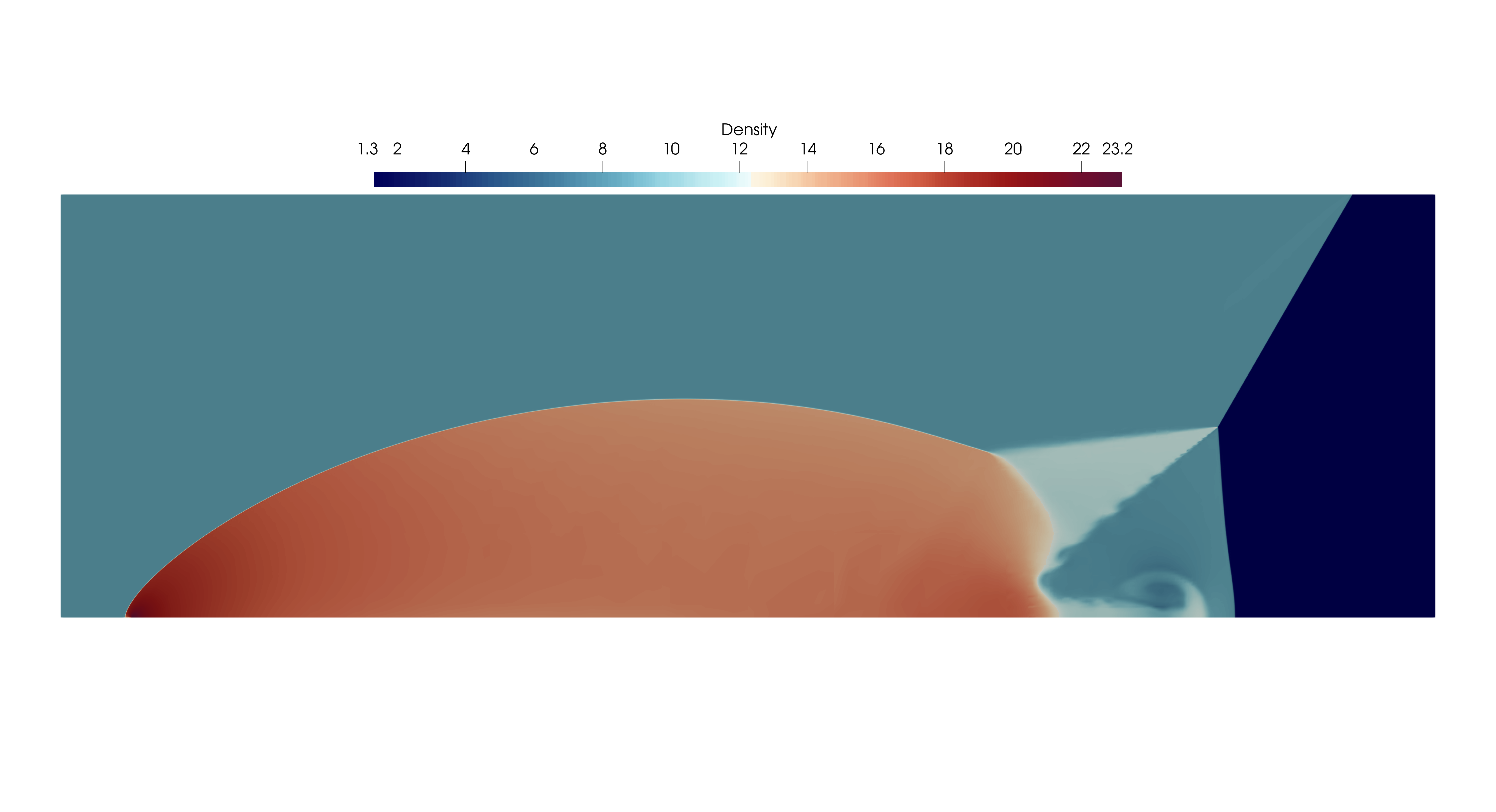}
    \includegraphics*[trim=200 760 200 760, clip, width=0.8\textwidth]{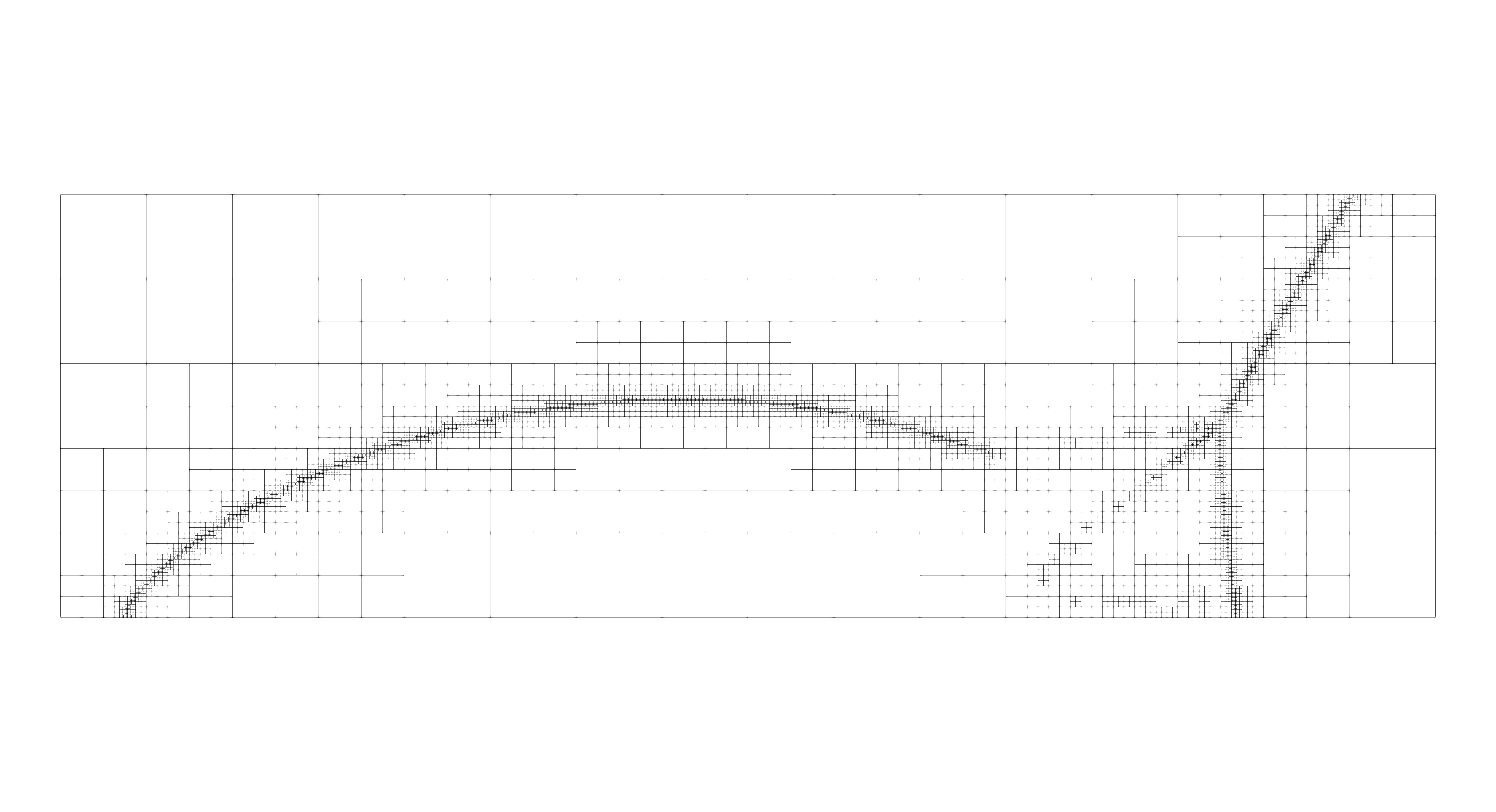}
    \caption{Density contours and resulting mesh of the simulation of the Double Mach reflection at $t=0.2$ with local limiting.}
    \label{fig:doublemach_local}
\end{figure}

For this setup, we first apply the standard local limiting to suppress oscillations and numerical artifacts. The result is shown in \Cref{fig:doublemach_local}. We also apply the positivity limiter of Zhang and Shu \cite{zhang2011} with a global threshold of $10^{-6}$ to ensure positivity for density and pressure during the solution transfer in the refinement and coarsening steps as described in \Cref{sec:LimitingTransfer}.

\begin{figure}[pos=htbp]
    \centering
    \includegraphics[clip, trim=0 0 0 200, width=0.48\linewidth]{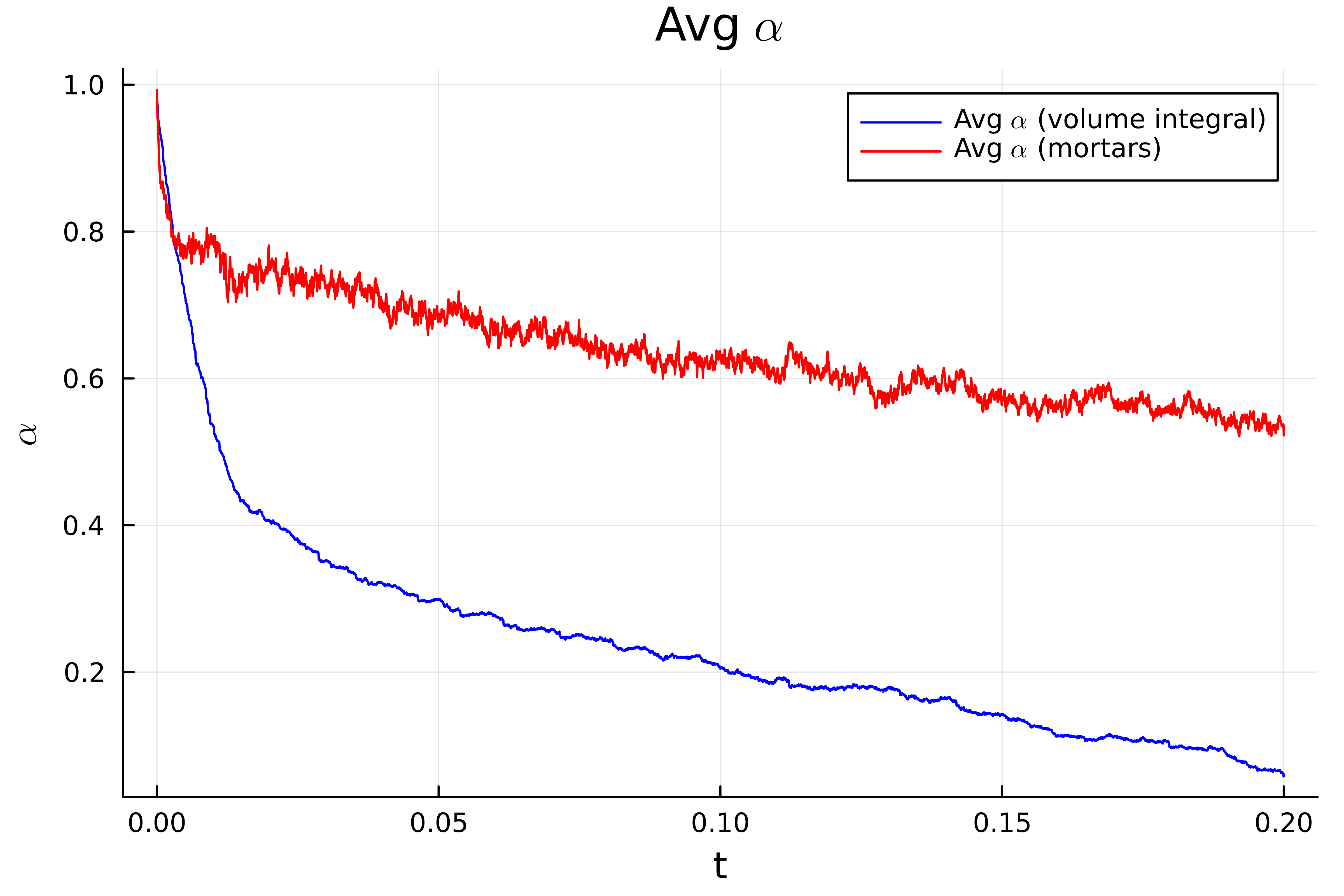}
    \caption{Evolution of the limiting factors.}
    \label{fig:doublemach_limitingfactor}
\end{figure}

The result in \Cref{fig:doublemach_local} shows clean shock fronts which are sharply resolved by AMR. In comparison, the interaction region between the shocks is not as finely resolved. The final mesh contains 6644 elements corresponding to 166100 degrees of freedom.
A detailed view of the interaction region on the left of the right shock front suggests a relatively large amount of limiting due to a very dissipative solution. The visualized evolution of the limiting factors in \Cref{fig:doublemach_limitingfactor} generally supports this theory. Nevertheless, note that there are reasons why the average mortar limiting factor can appear larger when compared to the limiting factor in the volume integral, as described earlier.

We want to slightly modify the setup to get a less dissipative solution, while still rendering the shock very sharply. One option here is to use pure positivity limiting instead of local limiting. In our tests, this provided insufficient limiting. Mesh imprinting becomes very visible and also the shocks are not sharply defined anymore.

We therefore combine both approaches. We separately compute limiting factors for pure positivity limiting and local limiting, $\alpha^\text{posLim}$ and $\alpha^\text{local}$.
In order to combine these, we use the element-wise shock-capturing indicator by Hennemann et al.~\cite{hennemann2021} applied to the pressure as smoothness indicator. This returns a continuous indicator variable, $\alpha_{\mathrm{indicator}}\in[0,1]$, where 0 indicates a smooth solution and 1 indicates a full shock.
So, to apply local limiting where shocks are detected, but use only positivity limiting otherwise, we compute the final limiting factors as follows:
\begin{equation}\label{eq:combination_posLim_local}
    \alpha = (1 - \alpha_\mathrm{indicator}) \, \alpha^\text{posLim} + \alpha_\mathrm{indicator} \, \alpha^\text{local}.
\end{equation}

\begin{figure}[pos=htbp]
    \begin{subfigure}[t]{\textwidth}
        \centering
        \includegraphics*[trim=200 760 200 460, clip, width=0.8\textwidth]{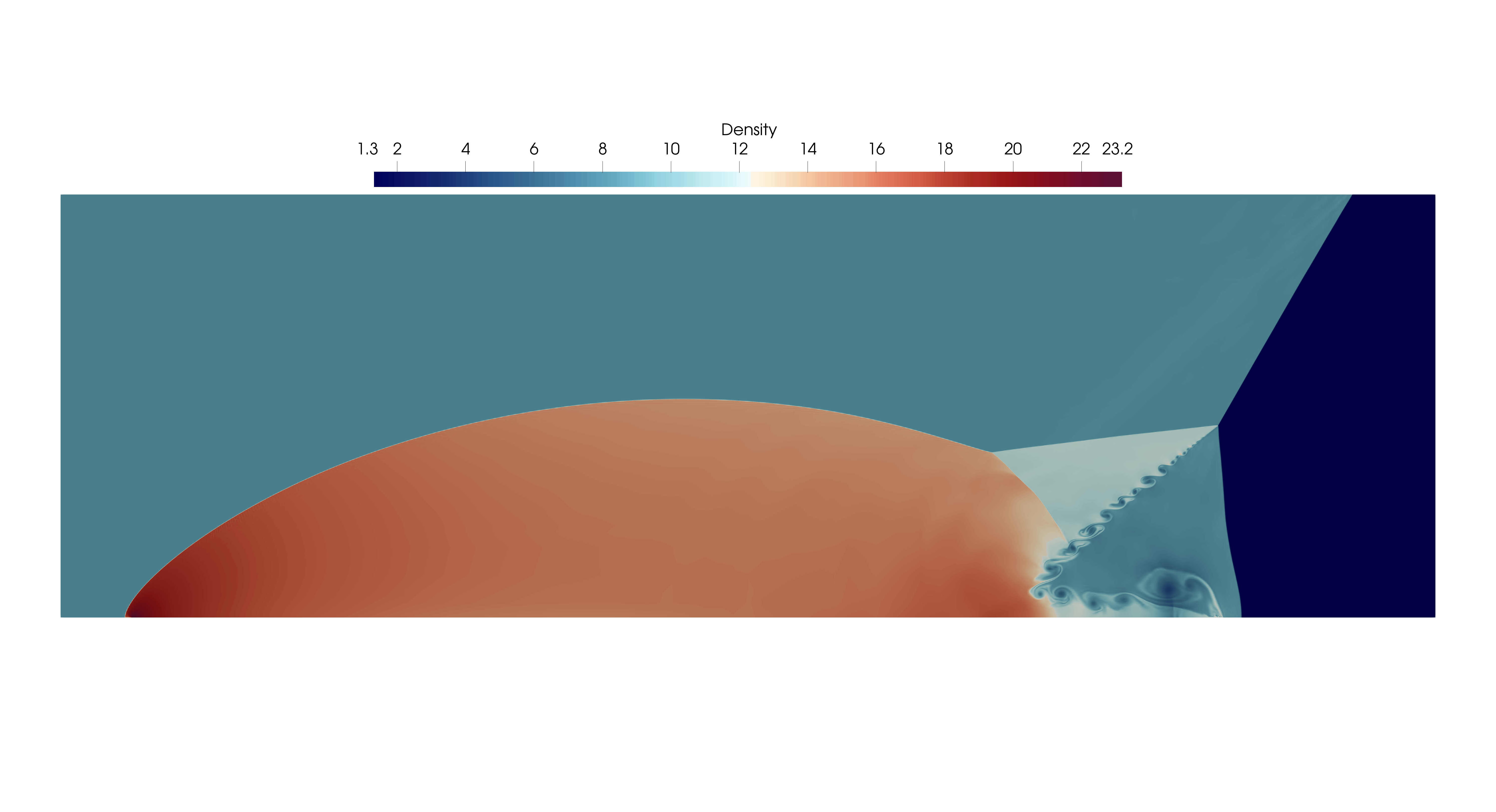}
        \includegraphics*[trim=200 760 200 760, clip, width=0.8\textwidth]{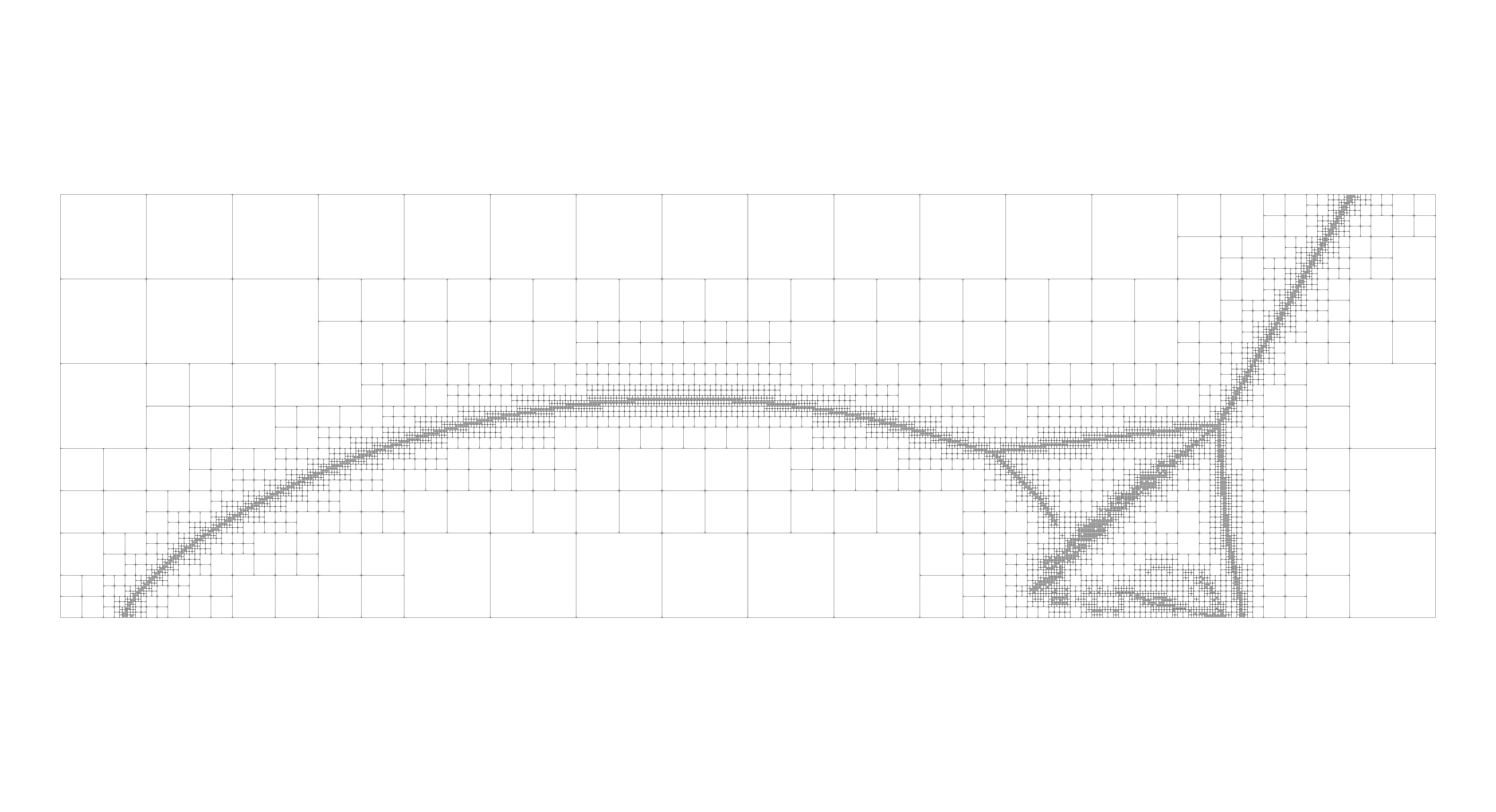}
        \subcaption{Density contours and resulting mesh.}
    \end{subfigure}
    \begin{subfigure}[t]{\textwidth}
        \centering
        \includegraphics*[trim=200 760 200 460, clip, width=0.8\textwidth]{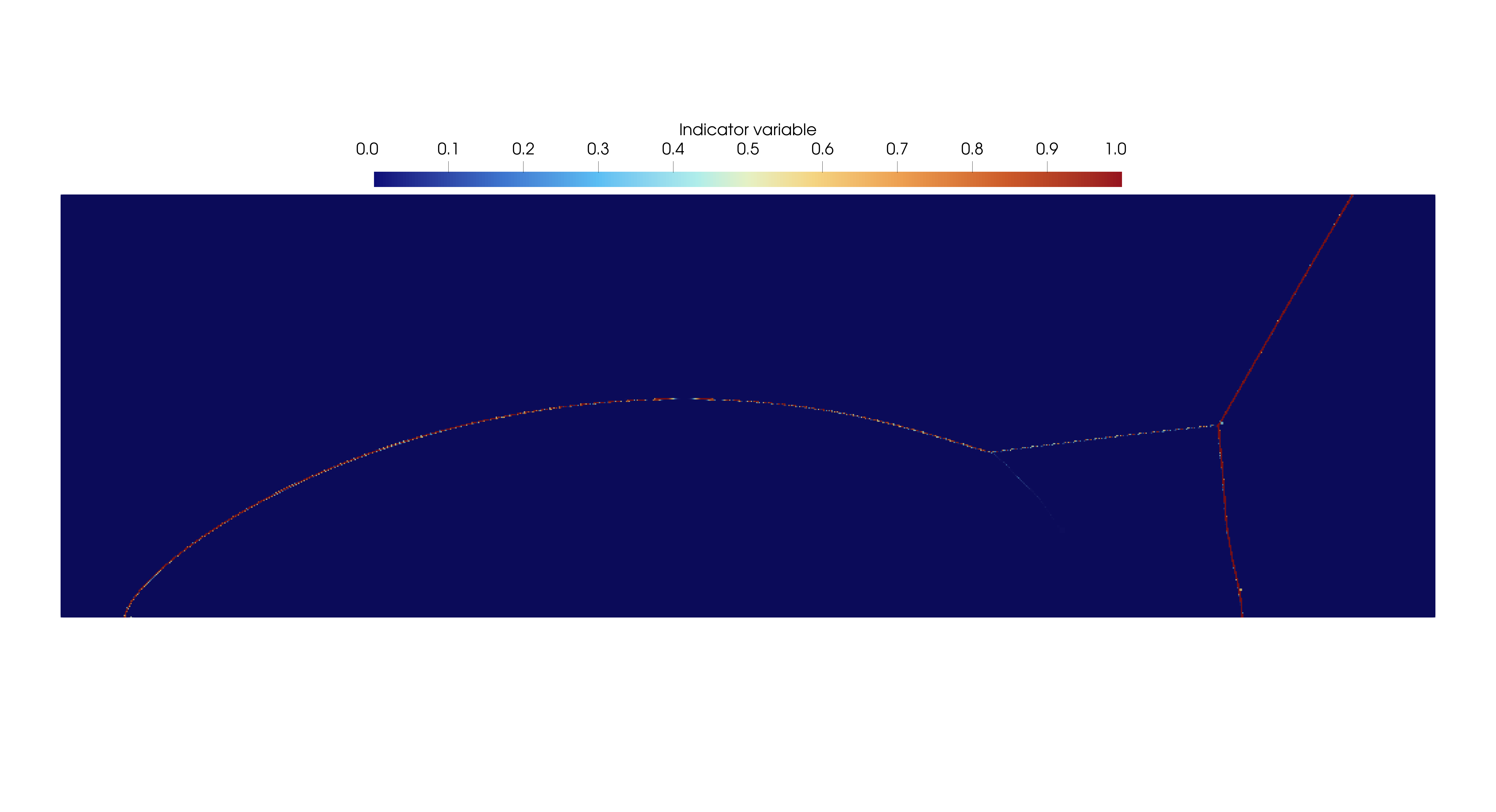}
        \subcaption{Visualization of the indicator variable $\alpha_\mathrm{indicator}$.}
    \end{subfigure}
    \caption{Results of the simulation of the Double Mach reflection at $t=0.2$ using a combination of positivity limiting and local limiting.}
    \label{fig:doublemach_indicator}
\end{figure}

Additionally, we add the new indicator as a second AMR indicator and refine an element to the maximum extent when $\alpha_\mathrm{indicator} \geq 0.5$.
The result with the new setup is shown in \Cref{fig:doublemach_indicator}.

A few things stand out. First, the shocks are still very sharply defined. In addition, the interaction region now contains significantly more small-scale structures.
This region is now also significantly more refined, resulting in a total of 10067 elements at time $t=0.2$. The plot of the indicator variable shows exactly the intended behavior. We are now using positivity limiting in smooth regions and local limiting at the shocks.

\subsection{High-Mach Astrophysical Jet}
A common benchmark to test the robustness of the used scheme is the astrophysical jet with a Mach number of about 2000. The setup was originally proposed in~\cite{ha2005numerical}.
The computational domain $\Omega=[-0.5, 0.5]^2$ contains a mono-atomic gas ($\gamma = 5/3$) and is at rest at the start of the simulation with
\begin{equation}
    \rho(x, y) = 0.5, \qquad
    p(x, y) = 0.4127, \qquad
    v_1(x, y) = 0, \qquad
    v_2(x, y) = 0.
\end{equation}

On the left, there is a hypersonic inflow with
\begin{align}
    \rho(-0.5,y_B) = 5, \qquad
    p(-0.5,y_B) = 0.4127, \qquad
    v_1(-0.5,y_B) = 800, \qquad
    v_2(-0.5,y_B) = 0,
\end{align}
for $y_B \in [-0.05, 0.05]$, which corresponds to a Mach number of $\textrm{Ma}=2156.91$ with respect to the speed of sound of the jet gas, and $\textrm{Ma}=682.08$ with respect to the speed of sound of the ambient gas.
We use periodic boundary conditions at the bottom and top boundaries and characteristic-based inflow/outflow boundary condition on the left and right.

We use the entropy-conserving and kinetic energy preserving flux of Chandrashekar~\cite{Chandrashekar2013} for the volume numerical fluxes.
The simulation again uses the previously described combination of positivity limiting and local limiting~\eqref{eq:combination_posLim_local}. We also apply the positivity limiter of Zhang and Shu \cite{zhang2011} with a global threshold of $10^{-10}$ to ensure positivity for density and pressure during the solution transfer in the refinement and coarsening steps.

We use $N=3$ and tessellate the domain using AMR with the combination of the AMR indicators by Löhner~\cite{LOHNER1987323} and Hennemann et al.~\cite{hennemann2021} for the density. We use refinement levels between 2 and 9 which corresponds to between $2^2=4$ and $2^9=512$ elements per spatial dimension. We use a CFL of 0.9 with the IDP time step restriction~\eqref{eq:timestep_restriction}.
\begin{figure}
    \centering
    \begin{subfigure}[t]{0.48\textwidth}
        \includegraphics[trim=3600 920 3600 360, clip, width=\textwidth]{figures/astrojet/density_mesh.png}
        \subcaption{Density contours and final mesh.}
    \end{subfigure}
    \begin{subfigure}[t]{0.48\textwidth}
        \includegraphics[trim=900 230 900 90, clip, width=\textwidth]{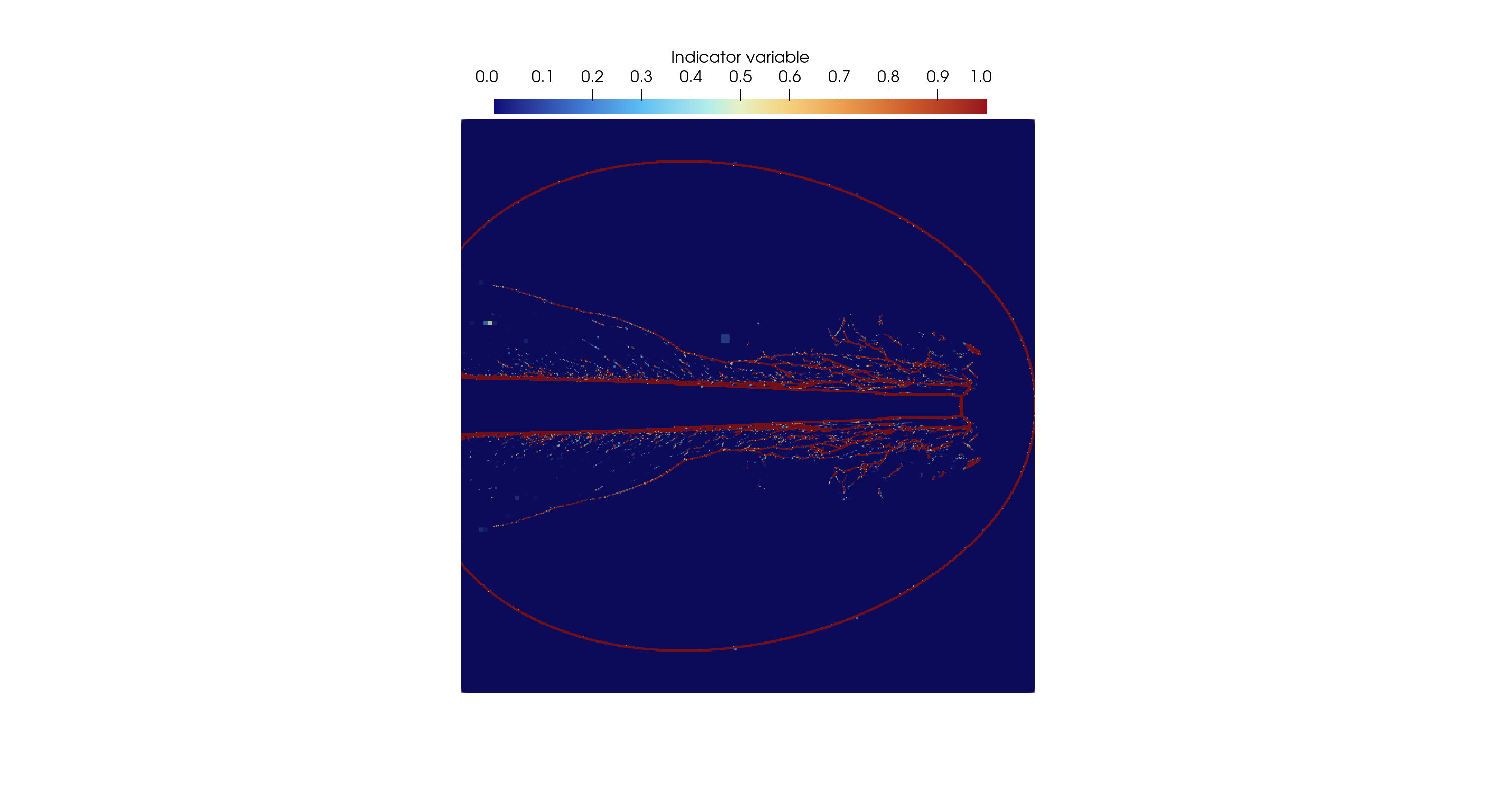}
        \subcaption{Visualization of the element-wise indicator.}
    \end{subfigure}
    \caption{Results of the simulation of the astrophysical jet at $t=0.0015$.}
    \label{fig:astrojet}
\end{figure}

\Cref{fig:astrojet} shows the results of the astrophysical jet simulation at time $t=0.0015$. The final mesh contains $51586$ elements corresponding to $825376$ degrees of freedom. A qualitative assessment is difficult for this setup. What is noticeable, however, is that the visualization of the indicator variable suggests that local limiting is again only used sparingly. Moreover, the resulting density contours show a sharply resolved front of the jet without visible numerical artifacts or carbuncles. The small-scale structures on both sides of the jet are symmetric, except for a few very minor irregularities.

\begin{remark}
    The results presented in this paper are, in some cases, highly sensitive to minor differences in the numerical setup. There are several reasons for this. Primarily, the chosen setups are inherently susceptible to such sensitivities, since calculations often involve very large numbers while also dealing with numbers very close to 0. For the latter, even small absolute errors can have a major impact.
    The design of the limiting strategy further contributes to this sensitivity, especially for nonlinear variables. In these cases, the optimal limiting factor is determined using a combination of Newton’s method and the bisection method.
    This was already the case for our subcell limiting scheme on conforming meshes, but becomes even more pronounced on nonconforming meshes.
    First, mortar limiting is more susceptible in this regard, since it uses one limiting factor per mortar and, as a result, has a less localized effect. Second, the use of AMR with heuristic refinement indicators further amplifies this behavior.
    This effect is particularly relevant for the astrophysical jet, where even minor differences in the setup can lead to significantly different results.

    Generally, the sensitivity is particularly pronounced in simulations that rely solely on positivity limiting due to the small amount of introduced dissipation, which would otherwise help mask the effect of small-scale differences. For this reason, we introduced the limiting setup that uses a combination of the positivity limiting and local limiting.

    In addition to this general sensitivity, we observe asymmetric solutions in some cases, even though both the initial conditions and the grid are symmetric w.r.t. the spatial domain. This occurs mainly with pure positivity limiting, since this method introduces less dissipation by not explicitly suppressing oscillations and therefore tends to produce solutions with more small-scale structures.
    This is known to be a possible problem of low-dissipation schemes, as discussed in \cite{FLEISCHMANN201994}, and is not specific to our method. This problem can be controlled, for example, by not using pure positivity limiting, but rather local limiting or the combination of the two described above.
\end{remark}

\section{Conclusion}\label{sec:conclusion}

We have developed a conservative and invariant-domain-preserving mortar formulation for nonconforming interfaces in Legendre--Gauss--Lobatto Discontinuous Galerkin Spectral Element Methods on Cartesian meshes.
The proposed approach extends graph-viscosity-based low-order schemes and convex limiting techniques to adaptive meshes containing hanging nodes. We adjusted existing limiting techniques to apply stabilizing limiting at the interface level. We introduced a new limiting strategy that combines positivity limiting with local limiting to maintain robustness while reducing numerical dissipation in smooth regions.

Starting from a conservative all-to-all mortar coupling, we derived a sparse interface discretization based on characteristic functions of LGL subcells.
The resulting mortar fluxes preserve conservation, satisfy the structural requirements of invariant-domain-preserving low-order schemes, and naturally reduce to the standard conforming interface formulation when neighboring elements share matching nodes.
The compact stencil obtained through sparsification avoids the excessive diffusion associated with dense mortar couplings and remains independent of the polynomial degree.

The numerical tests verify the expected convergence behavior in smooth regimes and demonstrate robustness for challenging compressible Euler simulations involving shocks, strong rarefactions, and adaptive refinement.
The proposed construction therefore provides a consistent low-order building block for convex limiting in DGSEM on (adaptively) refined meshes.
As a result, it enables the combination of high-order accuracy, robustness, and local mesh adaptivity within a common framework for nonlinear hyperbolic conservation laws.

\section*{Acknowledgments}

Andrés M. Rueda-Ramírez gratefully acknowledges funding from the Spanish Ministry of Science, Innovation, and Universities through the ``Beatriz Galindo'' grant (BG23-00062).
This project has received funding from the European Research Council (ERC) under the European Union's Horizon Europe research and innovation programme (grant agreement No. 101167322 - TRANSDIFFUSE).
Gregor J. Gassner and Andrés M. Rueda-Ramírez acknowledge funding through the German Federal Ministry for Education and Research (BMBF) project ``ICON-DG'' (01LK2315B) of the ``WarmWorld Smarter'' program.
Gregor J. Gassner further acknowledges funding from the German Research Foundation DFG through the research unit “SNuBIC” (DFG-FOR5409) and from Germany’s Excellence Strategy – EXC 3037 – 533607693 (Dynaverse).

\section*{Supplementary information}
All the implementations in the paper are carried out using the open-source code \texttt{Trixi.jl}~\cite{ranocha2022adaptive, schlottkelakemper2021purely, schlottkelakemper2025trixi}. The access to the data and the exact simulation configurations are granted upon request.

\appendix

\section{Illustrative calculation of local weights}\label{sec:local_weights}
In this section, we provide an illustrative computation of the local weights \eqref{eq:local_weights}.
The local weight $\tilde{\omega}_{(i, j)}^S$ from~\eqref{eq:local_weights} connects nodes $i\in\NN_S^-$ and $j\in\NN_{S}^+$.
Both nodes have an associated characteristic function $\psi_i^-$ and $\psi_j^+$. We denote the support of $\psi_i^-$ on $S$ by $[b_L^-, b_U^-]$ and the support of $\psi_j^+$ on $S$ by $[b_L^+, b_U^+]$. The local weight $\tilde\omega_{(i, j)}^S$ is only nonzero if those intervals overlap. We define the overlapping part of those two intervals as $[b_L, b_U]$ with $b_L = \max(b_L^-, b_L^+)$ and $b_U = \min(b_U^-, b_U^+)$. So,
\begin{equation}
    \tilde{\omega}_{(i, j)}^S
    = \int_S \psi_i^- \psi_j^+ \mathrm{d}s
    = \begin{cases}
        b_U - b_L, & \text{if } b_U > b_L,\\
        0, & \text{otherwise}.
    \end{cases}
\end{equation}
Due to symmetry of the quadrature weights, we only need to compute the weights between nodes of the large and one small element.

Since the spatial discretization for all elements is done on the reference element $[-1,1]^2$ and then mapped to the physical space, in practice, the computation of the local weights is done on a reference interface $S$ with 1D coordinates in $[-1,1]$. Also, the quadrature weights are based on the 1D LGL quadrature in $[-1,1]$, while the weights in the small elements on $S$ are then scaled by a factor of $\tfrac{1}{2}$.

We assume that the interface $S$ connects one large element on the left and two smaller elements on the right of the interface in $\xi_1$-direction as sketched \Cref{fig:SketchNonConformingInterface_Connections}. The interface nodes $i\in\NN_S^-$ have corresponding 2D indices $Nk_i$ and equivalently for $j\in\NN_S^+$ with $0k_j$. Assuming the right element is the lower small element, $k_i, k_j\in\{0,\dots,N\}$.
Then, due to the structure of the LGL subgrid, the upper and lower bounds $b_L$ and $b_U$ can be computed as follows,
\begin{equation}\begin{split}
b_L &= \max\left(b_L^-, b_L^+\right)
= -1 + \max\left(\sum_{\ell=0}^{k_i-1} \omega_{\ell, S}^-, \sum_{\ell=0}^{k_j-1} \omega_{\ell, S}^+ \right),\\
b_U &= \min\left(b_U^-, b_U^+\right)
= -1 + \min\left(\sum_{\ell=0}^{k_i} \omega_{\ell, S}^-, \sum_{\ell=0}^{k_j} \omega_{\ell, S}^+ \right).\\
\end{split}\end{equation}

We compute the local weights illustratively for a DG scheme with polynomial degree $N=3$, i.\,e., $N+1=4$ LGL nodes per dimension. We have $\{\omega_{\ell, S}^-\}_{\ell=0}^3 = \{\frac{1}{6}, \frac{5}{6}, \frac{5}{6}, \frac{1}{6}\}$ and $\omega_{j, S}^+ = \frac{1}{2}\,\omega_{j, S}^-$ for $j=0, \dots, 3$.
\begin{table}[pos=htbp]
    \centering
    {\renewcommand{\arraystretch}{1.2}
    {\small
    \begin{tabular}{@{}ll | cccc | cccc@{}} \toprule
        && \multicolumn{8}{c}{Node $j\in\NN_{S}^+$}\\
        && \multicolumn{4}{c|}{lower element}
         & \multicolumn{4}{c}{upper element}\\
        && 0 & 1 & 2 & 3 & 0 & 1 & 2 & 3\\
        \midrule
        \multirow{4}{*}{Node $i\in\NN_S^-$}
        & 0 & $\frac{1}{12}$ & $\frac{1}{12}$ & & & & & & \\
        & 1 & & $\frac{1}{3}$ & $\frac{5}{12}$ & $\frac{1}{12}$ & & & & \\
        & 2 & & & & & $\frac{1}{12}$ & $\frac{5}{12}$ & $\frac{1}{3}$ & \\
        & 3 & & & & & & & $\frac{1}{12}$ & $\frac{1}{12}$\\ \bottomrule
    \end{tabular}}}
    \caption{Local weights $\tilde{\omega}_{(i,j)}^S$ from \eqref{eq:local_weights} for $N=3$.}
    \label{table:local_weights}
\end{table}
The resulting local weights are shown in \Cref{table:local_weights}. This table provides a clear overview of the node connectivity sketched in \Cref{fig:SketchNonConformingInterface_Connections}.
The local weights for $N=4$ and $N=5$ are shown in \Cref{table:local_weights_4_5}.

\begin{table}[pos=htbp]
    \begin{subtable}{\textwidth}
    \centering
    {\renewcommand{\arraystretch}{1.2}
    {\small
    \begin{tabular}[pos=htbp]{@{}ll|ccccc|ccccc@{}} \toprule
        && \multicolumn{10}{c}{Node $j\in\NN_{S}^+$}\\
        && \multicolumn{5}{c|}{lower element} & \multicolumn{5}{c}{upper element}\\
        && 0 & 1 & 2 & 3 & 4 & 0 & 1 & 2 & 3 & 4\\
        \midrule
        \multirow{5}{*}{Node $i\in\NN_S^-$}
        & 0 & $\frac{1}{20}$ & $\frac{1}{20}$ & & & & & & & &\\
        & 1 & & $\frac{2}{9}$ & $\frac{29}{90}$ & & & & & & &\\
        & 2 & & & $\frac{1}{30}$ & $\frac{49}{180}$ & $\frac{1}{20}$ & $\frac{1}{20}$ & $\frac{49}{180}$ & $\frac{1}{30}$ & &\\
        & 3 & & & & & & & & $\frac{29}{90}$ & $\frac{2}{9}$ &\\
        & 4 & & & & & & & & & $\frac{1}{20}$ & $\frac{1}{20}$\\ \bottomrule
    \end{tabular}}}
    \subcaption{$N=4$ with $\{\omega_{\ell, S}^-\}_{\ell=0}^4 = \left\{\frac{1}{10}, \frac{49}{90}, \frac{32}{45}, \frac{49}{90}, \frac{1}{10}\right\}$.}
    \end{subtable}
    \begin{subtable}{\textwidth}
    \centering
    {\renewcommand{\arraystretch}{1.2}
    {\footnotesize
    \begin{tabular}[pos=htbp]{@{}ll|cccccc|cccccc@{}} \toprule
        && \multicolumn{12}{c}{Node $j\in\NN_{S}^+$}\\
        && \multicolumn{6}{c|}{lower element} & \multicolumn{6}{c}{upper element}\\
        && 0 & 1 & 2 & 3 & 4 & 5 & 0 & 1 & 2 & 3 & 4 & 5\\
        \midrule
        \multirow{6}{*}{Node $i\in\NN_S^-$}
        & 0 & $\frac{1}{30}$ & $\frac{1}{30}$ & & & & & & & & & & \\
        & 1 & & $\frac{12-\sqrt{7}}{60}$ & $\frac{16-\sqrt{7}}{60}$ & & & & & & & & & \\
        & 2 & & & $\frac{-1+\sqrt{7}}{30}$ & $\frac{14+\sqrt{7}}{60}$ & $\frac{14-\sqrt{7}}{60}$ & $\frac{1}{30}$ & & & & & & \\
        & 3 & & & & & & & $\frac{1}{30}$ & $\frac{14-\sqrt{7}}{60}$ & $\frac{14+\sqrt{7}}{60}$ & $\frac{-1+\sqrt{7}}{30}$ & & \\
        & 4 & & & & & & & & & & $\frac{16-\sqrt{7}}{60}$ & $\frac{12-\sqrt{7}}{60}$ & \\
        & 5 & & & & & & & & & & & $\frac{1}{30}$ & $\frac{1}{30}$ \\ \bottomrule
    \end{tabular}}}
    \subcaption{$N=5$ with $\{\omega_{\ell, S}^-\}_{\ell=0}^5 = \left\{\frac{1}{15}, \frac{14-\sqrt{7}}{30}, \frac{14+\sqrt{7}}{30}, \frac{14+\sqrt{7}}{30}, \frac{14-\sqrt{7}}{30}, \frac{1}{15}\right\}$.}
    \end{subtable}
    \caption{Local weights $\tilde{\omega}_{(i,j)}^S$ for $N=4$ and $N=5$.}
    \label{table:local_weights_4_5}
\end{table}

\printcredits

\bibliographystyle{model1-num-names}

\bibliography{cas-bibliography}



\end{document}